\documentclass[12pt]{amsart}
\usepackage{amsmath,amssymb,mathrsfs,epsfig,rotate,color,
url,hyperref}
\usepackage[all]{xy}
\usepackage{fancyvrb}
\newenvironment{tobiwrite}[1]
  {\typeout{Writing file #1}\VerbatimOut{#1}}
  {\endVerbatimOut}

\def\ev{\pi}
\def\nf{\qopname\relax o{nf}}
\def\mcg{\mathrm{MCG}}
\def\homeo{\mathrm{Homeo}}
\def\cwl{{\mathrm{wl}}}
\def\zwl{{\mathrm{zwl}}}

\newtheorem{prop}{Proposition}[section]
\newtheorem{lemm}[prop]{Lemma}
\newtheorem{coro}[prop]{Corollary}
\newtheorem{theo}[prop]{Theorem}
\theoremstyle{definition}
\newtheorem{defi}{Definition}[section]
\newtheorem{rem}{Remark}[section]
\newtheorem{exam}{Example}[section]
\author{Ivan Dynnikov}
\title{Counting intersections of normal curves}
\address{\noindent Steklov Mathematical Institute of Russian Academy of Sciences, 8 Gubkina Str., Moscow 119991, Russia}
\email{dynnikov@mech.math.msu.su}
\thanks{\copyright\ 2021. This manuscript version is made available under the CC-BY-NC-ND 4.0 license \url{https://creativecommons.org/licenses/by-nc-nd/4.0/}}

\begin{tobiwrite}{denorm7.eps}
/v {gsave translate
4 0 moveto 0 0 4 0 360 arc gsave
1 setgray fill grestore 0 setgray
2 setlinewidth stroke
grestore} def

gsave
newpath
15 60 moveto
30 80 lineto
33 84 45 87 50 85 curveto
75 75 lineto
80 73 85 65 85 62 curveto
85 38 lineto
85 35 80 27 75 25 curveto
50 15 lineto
45 13 33 16 30 20 curveto
15 40 lineto
12 42 12 58 15 60 curveto
closepath
50 50 moveto
65 60 85 70 85 55 curveto
85 45 moveto
85 30 105 40 120 50 curveto
85 38 moveto
85 33 91 27 95 25 curveto
115 15 lineto
117 14 123 14 125 15 curveto
145 25 lineto
149 27 155 33 155 38 curveto
155 62 lineto
155 67 149 73 145 75 curveto
125 85 lineto
123 86 117 86 115 85 curveto
95 75 lineto
91 73 85 67 85 62 curveto
15 60 moveto
12 56 7 53 5 52 curveto
15 40 moveto
12 44 7 47 5 48 curveto
30 80 moveto
33 84 35 90 37 95 curveto
30 20 moveto
33 16 35 10 37 5 curveto
50 85 moveto
45 87 43 92 40 95 curveto
50 15 moveto
45 13 43 8 40 5 curveto
75 75 moveto
80 73 91 73 95 75 curveto
75 25 moveto
80 27 91 27 95 25 curveto
115 15 moveto
117 14 117 8 118 5 curveto
125 15 moveto
123 14 123 8 122 5 curveto
115 85 moveto
117 86 117 92 118 95 curveto
125 85 moveto
123 86 123 92 122 95 curveto
155 38 moveto
155 35 160 31 163 28 curveto
155 62 moveto
155 65 160 69 163 72 curveto
145 25 moveto
149 27 155 27 160 24 curveto
145 75 moveto
149 73 155 73 160 76 curveto
.5 setlinewidth {2.5 2.5} 0 setdash .7 setgray stroke
grestore

5 20 moveto 50 50 lineto
5 80 lineto
50 50 moveto
165 50 lineto
65 5 moveto 50 50 lineto
65 95 lineto
100 95 moveto 120 50 lineto
100 5 lineto
145 95 moveto
120 50 lineto
145 5 lineto
2 setlinewidth stroke

newpath
50 50 moveto
70 63 83 63 83 50 curveto
83 30 97 25 107 20 curveto
117 15 123 15 133 20 curveto
143 25 154 30 154 50 curveto
154 70 143 75 133 80 curveto
123 85 117 85 107 80 curveto
97 75 87 70 87 50 curveto
87 37 100 37 120 50 curveto
1 setlinewidth stroke

50 50 v
120 50 v

showpage
\begin{tobiwrite}{cleanup3.eps}
newpath
0 50 moveto
5 40 25 30 35 30 curveto
80 50 moveto
75 40 55 30 45 30 curveto stroke
newpath 0 30 moveto 80 0 rlineto 2 setlinewidth stroke
showpage
\end{tobiwrite}
\begin{tobiwrite}{corner2.eps}

100 50 moveto
50 50 50 0 360 arc clip

gsave
newpath
-10 -10 moveto
30 30 lineto
45 45 55 45 70 30 curveto
110 -10 lineto
closepath gsave .9 setgray fill
grestore 2 setlinewidth stroke
grestore

newpath
5 7 100 {0 exch moveto
100 0 rlineto} for
.5 setlinewidth stroke

showpage
\end{tobiwrite}
\begin{tobiwrite}{tight3.eps}
newpath
3 30 moveto
3 30 3 0 360 arc fill
newpath
3 30 moveto
18 60 60 60 80 0 curveto
3 30 moveto
18 0 60 0 80 60 curveto
stroke
showpage
\end{tobiwrite}
\begin{tobiwrite}{triang2track2.eps}
/v {gsave translate newpath 2 0 moveto
0 0 2 0 360 arc
gsave 1 setgray fill grestore
0 setgray 0.7 setlinewidth stroke grestore} def
0.5 setlinewidth
newpath
110 60 moveto
115 80 lineto
120 100 140 100 135 100 curveto
165 100 lineto
175 100 180 95 185 85 curveto
200 55 lineto
205 45 200 20 190 15 curveto
170 5 lineto
160 0 143 -1 135 5 curveto
115 20 lineto
107 26 105 40 110 60 curveto
115 80 moveto
120 100 122 104 120 110 curveto
110 140 lineto
108 146 98 148 94 148 curveto
62 148 lineto
56 148 52 141 48 135 curveto
32 111 lineto
28 105 30 91 35 86 curveto
51 70 lineto
54 67 60 62 65 60 curveto
90 50 lineto
100 46 105 40 110 60 curveto
135 100 moveto
125 100 123 101 120 110 curveto
165 100 moveto
180 100 186 105 188 110 curveto
200 140 lineto
202 145 205 152 200 163 curveto
185 190 lineto
182.5 194.5 175 195 165 195 curveto
135 195 lineto
130 195 123 194 120 185 curveto
112 161 lineto
110 155 108 146 110 140 curveto
185 85 moveto
182 91 186 105 188 110 curveto
200 140 moveto
204 150 210 155 220 155 curveto
240 155 lineto
250 155 260 145 265 135 curveto
275 115 lineto
280 105 277 90 275 85 curveto
265 60 lineto
261 50 255 48 245 48 curveto
220 48 lineto
213 48 205 45 200 55 curveto
115 20 moveto
95 35 100 46 90 50 curveto
65 60 moveto
60 62 58 62 56 60 curveto
51 70 moveto
53 67 53 65 53 62 curveto
35 86 moveto
33 89 31 90 28 91 curveto
32 111 moveto
30 108 28 106 27 105 curveto
48 135 moveto
50 138 52 143 52 148 curveto
62 148 moveto
60 148 58 150 56 152 curveto
94 148 moveto
100 148 110 155 112 161 curveto
120 185 moveto
121 188 122 192 122 195 curveto
135 195 moveto
133 195 129 196 127 199 curveto
165 195 moveto
168 195 170 196 174 199 curveto
185 190 moveto
180 199 182 199 182 199 curveto
200 163 moveto
205 152 210 155 220 155 curveto
240 155 moveto
245 155 247 156 250 158 curveto
265 135 moveto
260 145 257 151 256 155 curveto
275 115 moveto
280 105 283 103 285 102 curveto
275 85 moveto
279 91 282 95 285 96 curveto
265 60 moveto
263 55 262 50 262 48 curveto
245 48 moveto
250 48 252 47 254 45 curveto
220 48 moveto
205 48 210 25 190 15 curveto
170 5 moveto
164 2 163 1 162 0 curveto
135 5 moveto
139 2 140 1 141 0 curveto
70 100 moveto
60 70 51 70 43 78 curveto
150 45 moveto
125 90 115 80 112.5 70 curveto
150 155 moveto
195 150 200 140 194.5 126 curveto
230 100 moveto
210 55 200 55 191 73 curveto
0.7 setgray
[2 2] 0 setdash
stroke
newpath
150 45 moveto
132 90 115 90 112 72 curveto
106 36 104 16 139 2 curveto
144 0 180 0 200 20 curveto
215 35 205 46 240 46 curveto
260 46 262 50 275 85 curveto
280 100 280 105 268 129 curveto
256 153 248 155 222 155 curveto
208 155 205 150 195 120 curveto
185 90 185 85 195 70 curveto
201 61 220 76 230 100 curveto
220 72 201 57 193 69 curveto
181 87 184 95 194 125 curveto
200 143 200 155 150 155 curveto
195 150 192 138 185 117 curveto
180 102 174 102 150 102  curveto
125 102 123 105 116 126 curveto
108 150 104 150 77 149 curveto
55 149 54 147 38 123 curveto
22 99 27 92 41 78 curveto
51 68 62  70 70 100 curveto
60 72 53 70 41 82 curveto
27 96 25 99 41 123 curveto
57 145 58 146 77 146 curveto
104 146 105 146 113 126 curveto
120 102 122 99 150 99 curveto
170 99 175 93 188 67 curveto
200 43 200 18 178 12 curveto
166 6 143 3 127 15 curveto
115 24 113 33 116 57 curveto
119 81 124 80 150 45 curveto
124 86 118 87 114 63 curveto
110 31 109 25 129 10 curveto
141 1 164 3 182 11 curveto
204 22 204 41 188 73 curveto
180 89 181 100 191 125 curveto
197 140 195 153 150 155 curveto
0 setgray
[] 0 setdash
stroke
70 100 v
150 45 v
150 155 v
230 100 v
showpage
\end{tobiwrite}
\begin{tobiwrite}{cleanup10.eps}
newpath
35 50 moveto
40 40 60 40 70 30 curveto
stroke
newpath 0 30 moveto 70 0 rlineto 2 setlinewidth stroke
70 30 translate
newpath
2.5 0 moveto 0 0 2.5 0 360 arc 1 setlinewidth gsave
1 setgray fill grestore 0 setgray stroke
showpage
\end{tobiwrite}
\begin{tobiwrite}{tight10.eps}
newpath
3 30 moveto
3 30 3 0 360 arc fill
newpath
3 30 moveto
30 80 52 33 62 33 curveto
72 33 77 51 92 51 curveto
100 51 110 51 140 0 curveto
3 30 moveto
30 -20 52 27 62 27 curveto
72 27 77 9 92 9 curveto
100 9 110 9 140 60 curveto
stroke
showpage
\end{tobiwrite}
\begin{tobiwrite}{case233a.eps}
20 0 translate
newpath
40 0 moveto
50 10 50 15 50 25 curveto
50 95 lineto
50 105 50 110 60 120 curveto
60 0 moveto
50 10 50 15 50 25 curveto
50 30 moveto
50 40 50 45 40 55 curveto
20 55 moveto
30 65 30 70 30 80 curveto
30 95 lineto
30 105 30 110 20 120 curveto
30 80 moveto
30 95 50 80 50 95 curveto
stroke
showpage
\end{tobiwrite}
\begin{tobiwrite}{or.eps}
5 15 translate
newpath
0 0 moveto
20 15 50 15 50 0 curveto
50 -15 moveto 50 15 lineto
stroke
2.5 0 moveto
0 0 2.5 0 360 arc
gsave 1 setgray fill grestore 0 setgray stroke
-70 rotate
newpath
10 0 moveto
0 0 10 0 140 arc
stroke
140 rotate
newpath
10 3 moveto
2 -5 rlineto
-2 2 rlineto
-2 -2 rlineto
closepath fill
showpage
\end{tobiwrite}
\begin{tobiwrite}{slide2.eps}
newpath
0 30 moveto
100 30 lineto
0 0 moveto
6 5 45 30 90 30 curveto
0 60 moveto
20 40 40 30 50 30 curveto
stroke

showpage
\end{tobiwrite}
\begin{tobiwrite}{univ3.eps}
6 8 translate
/v {gsave translate newpath
4 0 moveto 0 0 4 0 360 arc gsave
1 setgray fill grestore 0 setgray
stroke grestore} def
/r 50 3 sqrt div def

2 setlinewidth
newpath
0 0 moveto
100 0 lineto
50 50 3 sqrt mul lineto
closepath stroke

gsave
newpath
0 0 moveto
25 r 2 div 50 r 2 div 50 r 8 div curveto
50 r translate
0 r neg moveto
0 r 2 mul 3 div neg lineto
0 r 3 div neg
r 3 div dup 30 cos mul exch 30 sin mul
r 2 mul 3 div dup 30 cos mul exch 30 sin mul curveto
r dup 30 cos mul exch 30 sin mul lineto
r 2 mul 3 div dup 30 cos mul exch 30 sin mul moveto
r 3 div dup 30 cos mul exch 30 sin mul
r 3 div dup 150 cos mul exch 150 sin mul
r 2 mul 3 div dup 150 cos mul exch 150 sin mul curveto
r dup 150 cos mul exch 150 sin mul lineto
r 2 mul 3 div dup 150 cos mul exch 150 sin mul moveto
r 3 div dup 150 cos mul exch 150 sin mul
0 r 3 div neg
0 r 2 mul 3 div neg curveto
1 setlinewidth stroke
grestore
0 0 v
100 0 v
50 50 3 sqrt mul v
showpage
\end{tobiwrite}
\begin{tobiwrite}{cleanup8.eps}
newpath
0 50 moveto
5 40 25 40 55 40 curveto
30 50 moveto
35 40 50 40 55 40 curveto
65 40 65 30 75 30 curveto
15 10 moveto
20 20 40 30 50 30 curveto stroke
newpath 0 30 moveto 80 0 rlineto 2 setlinewidth stroke
showpage
\end{tobiwrite}
\begin{tobiwrite}{norm1.eps}
0 111 translate
-90 rotate
6 93 translate
/v {gsave translate newpath
4 0 moveto 0 0 4 0 360 arc gsave
1 setgray fill grestore 0 setgray
stroke grestore} def
newpath
0 0 moveto
30 20 50 20 50 0 curveto
gsave
closepath
0.8 setgray fill
0 setgray grestore
50 -20 70 -20 100 0 curveto
1 setlinewidth
stroke
2 setlinewidth
newpath
0 0 moveto
100 0 lineto
50 50 3 sqrt mul lineto
closepath
50 -50 3 sqrt mul lineto
100 0 lineto
stroke
0 0 v
100 0 v
50 50 3 sqrt mul v
50 -50 3 sqrt mul v
-80 rotate
newpath
10 0 moveto
0 0 10 0 160 arc
1 setlinewidth stroke
160 rotate
newpath
10 3 moveto
2 -5 rlineto
-2 2 rlineto
-2 -2 rlineto
closepath fill
showpage
\end{tobiwrite}
\begin{tobiwrite}{case0.eps}
newpath
40 20 moveto
50 30 50 40 50 50 curveto
50 70 lineto
50 80 50 90 40 100 curveto
60 20 moveto
50 30 50 40 50 50 curveto
50 70 moveto
50 80 50 90 60 100 curveto
stroke
showpage
\end{tobiwrite}
\begin{tobiwrite}{tight2.eps}
newpath
0 0 moveto
25 50 25 15 50 15 curveto
75 15 75 50 100 0 curveto
0 60 moveto
25 10 25 45 50 45 curveto
75 45 75 10 100 60 curveto
stroke
showpage
\end{tobiwrite}
\begin{tobiwrite}{corner3.eps}

100 50 moveto
50 50 50 0 360 arc clip

gsave
newpath
0 0 moveto
50 50 lineto 100 0 lineto
closepath .9 setgray fill
grestore

newpath
0 0 moveto
100 100 lineto
100 0 moveto
0 100 lineto
2 setlinewidth stroke

newpath
5 7 100 {0 moveto
0 100 rlineto} for
.5 setlinewidth stroke

50 50 translate
newpath
0 0 moveto
0 0 4 0 360 arc fill

showpage
\end{tobiwrite}
\begin{tobiwrite}{cleanup2.eps}
newpath 0 30 moveto 80 0 rlineto 2 setlinewidth stroke
newpath
30 10 moveto
35 20 45 40 50 50 curveto {4 2} 0 setdash 1 setlinewidth
30 10 moveto
35 20 55 30 65 30 curveto
50 50 moveto
45 40 25 30 15 30 curveto stroke
showpage
\end{tobiwrite}
\begin{tobiwrite}{case24c.eps}
newpath
10 110 moveto
20 100 20 90 20 80 curveto
20 40 lineto
20 30 20 20 30 10 curveto
40 70 moveto
40 10 90 10 90 70 curveto
90 130 40 130 40 70 curveto
stroke

gsave
newpath
90 70 translate
7 0 moveto
0 0 7 0 360 arc
gsave 1 setgray fill grestore 0 setgray stroke

45 rotate
newpath
7 0 moveto -7 0 lineto
0 7 moveto 0 -7 lineto
stroke
grestore
newpath
20 60 translate
7 0 moveto
0 0 7 0 360 arc
gsave 1 setgray fill grestore 0 setgray stroke

45 rotate
newpath
7 0 moveto -7 0 lineto
0 7 moveto 0 -7 lineto
stroke
showpage
\end{tobiwrite}
\begin{tobiwrite}{tight5.eps}
newpath
0 0 moveto
30 60 60 60 80 30 curveto
100 0 130 0 160 60 curveto
0 60 moveto
30 0 60 0 80 30 curveto
100 60 130 60 160 0 curveto
stroke
showpage
\end{tobiwrite}
\begin{tobiwrite}{case25.eps}
newpath
40 30 moveto
50 40 50 50 50 60 curveto
50 80 lineto
50 90 50 100 40 110 curveto
60 30 moveto
50 40 50 50 50 60 curveto
50 80 moveto
50 90 50 100 60 110 curveto
70 120 90 110 90 60 curveto
90 10 70 5 40 5 curveto
20 5 30 20 40 30 curveto
stroke
showpage
\end{tobiwrite}
\begin{tobiwrite}{case231a.eps}
newpath
40 0 moveto
50 10 50 15 50 25 curveto
50 30 lineto
50 30 moveto
50 40 50 45 40 55 curveto
40 65 moveto
50 75 50 80 50 90 curveto
50 115 lineto
40 75 moveto
50 85 50 90 50 100 curveto
80 0 moveto
70 10 70 15 70 25 curveto
70 80 lineto
70 95 50 95 50 110 curveto
50 25 moveto
50 40 70 40 70 55 curveto
stroke
showpage
\end{tobiwrite}
\begin{tobiwrite}{split3.eps}
newpath
0 10 moveto
20 30 40 30 60 30 curveto
96 30 lineto
0 50 moveto
20 30 40 30 60 30 curveto
stroke
newpath

99 30 moveto
96 30 3 0 360 arc
gsave 1 setgray fill
grestore 0 setgray stroke
showpage
\end{tobiwrite}
\begin{tobiwrite}{denorm3.eps}
/v {gsave translate
4 0 moveto 0 0 4 0 360 arc gsave
1 setgray fill grestore 0 setgray
2 setlinewidth stroke
grestore} def

gsave
newpath
15 60 moveto
30 80 lineto
33 84 45 87 50 85 curveto
75 75 lineto
80 73 85 65 85 62 curveto
85 38 lineto
85 35 80 27 75 25 curveto
50 15 lineto
45 13 33 16 30 20 curveto
15 40 lineto
12 42 12 58 15 60 curveto
closepath
50 50 moveto
65 60 85 70 85 55 curveto
85 45 moveto
85 30 105 40 120 50 curveto
85 38 moveto
85 33 91 27 95 25 curveto
115 15 lineto
117 14 123 14 125 15 curveto
145 25 lineto
149 27 155 33 155 38 curveto
155 62 lineto
155 67 149 73 145 75 curveto
125 85 lineto
123 86 117 86 115 85 curveto
95 75 lineto
91 73 85 67 85 62 curveto
15 60 moveto
12 56 7 53 5 52 curveto
15 40 moveto
12 44 7 47 5 48 curveto
30 80 moveto
33 84 35 90 37 95 curveto
30 20 moveto
33 16 35 10 37 5 curveto
50 85 moveto
45 87 43 92 40 95 curveto
50 15 moveto
45 13 43 8 40 5 curveto
75 75 moveto
80 73 91 73 95 75 curveto
75 25 moveto
80 27 91 27 95 25 curveto
115 15 moveto
117 14 117 8 118 5 curveto
125 15 moveto
123 14 123 8 122 5 curveto
115 85 moveto
117 86 117 92 118 95 curveto
125 85 moveto
123 86 123 92 122 95 curveto
155 38 moveto
155 35 160 31 163 28 curveto
155 62 moveto
155 65 160 69 163 72 curveto
145 25 moveto
149 27 155 27 160 24 curveto
145 75 moveto
149 73 155 73 160 76 curveto
.5 setlinewidth {2.5 2.5} 0 setdash .7 setgray stroke
grestore

5 20 moveto 50 50 lineto
5 80 lineto
50 50 moveto
165 50 lineto
65 5 moveto 50 50 lineto
65 95 lineto
100 95 moveto 120 50 lineto
100 5 lineto
145 95 moveto
120 50 lineto
145 5 lineto
2 setlinewidth stroke

newpath
50 50 moveto
70 65 85 65 85 50 curveto
85 35 100 35 120 50 curveto
1 setlinewidth stroke

50 50 v
120 50 v

showpage
\begin{tobiwrite}{slide1.eps}
newpath
0 30 moveto
100 30 lineto
0 0 moveto
10 25 15 30 25 30 curveto
0 60 moveto
20 40 40 30 50 30 curveto
stroke
showpage
\end{tobiwrite}
\begin{tobiwrite}{or2.eps}
5 15 translate
1 -1 scale
newpath
0 0 moveto
20 15 50 15 50 0 curveto
50 -15 moveto 50 15 lineto
stroke
2.5 0 moveto
0 0 2.5 0 360 arc
gsave 1 setgray fill grestore 0 setgray stroke
-70 rotate
newpath
10 0 moveto
0 0 10 0 140 arc
stroke
140 rotate
newpath
10 3 moveto
2 -5 rlineto
-2 2 rlineto
-2 -2 rlineto
closepath fill
showpage
\end{tobiwrite}
\begin{tobiwrite}{case3a.eps}
newpath
40 20 moveto
50 30 50 40 50 50 curveto
50 70 lineto
50 80 50 90 40 100 curveto
70 20 moveto
60 30 60 95 70 110 curveto
50 70 moveto
50 80 50 90 70 110 curveto
stroke
70 110 translate
newpath
3 0 moveto
0 0 3 0 360 arc
gsave 1 setgray fill grestore 0 setgray stroke
showpage
\end{tobiwrite}
\begin{tobiwrite}{thicktrack3.eps}
newpath
0 50 moveto
30 50 lineto
40 50 50 50 100 60 curveto
30 50 moveto
40 50 50 50 100 40 curveto
stroke

/Times-Roman findfont 12 scalefont setfont
15 53 moveto
(7) show
82 61 moveto
(2) show
82 31 moveto
(5) show
showpage
\end{tobiwrite}
\begin{tobiwrite}{norm3.eps}
0 111 translate
-90 rotate
6 93 translate
/v {gsave translate newpath
4 0 moveto 0 0 4 0 360 arc gsave
1 setgray fill grestore 0 setgray
stroke grestore} def
newpath
0 0 moveto
30 20 50 20 50 0 curveto
gsave
closepath
0.8 setgray fill
0 setgray grestore
50 -20
/a 50 3 sqrt mul def
a 0.7 mul 30 cos mul
a 0.7 mul -30 sin mul
a 4 add 30 cos mul
a 4 add -30 sin mul
curveto
1 setlinewidth
stroke
2 setlinewidth
newpath
0 0 moveto
100 0 lineto
50 50 3 sqrt mul lineto
closepath
50 -50 3 sqrt mul lineto
100 0 lineto
stroke
0 0 v
100 0 v
50 50 3 sqrt mul v
50 -50 3 sqrt mul v
-80 rotate
newpath
10 0 moveto
0 0 10 0 160 arc
1 setlinewidth stroke
160 rotate
newpath
10 3 moveto
2 -5 rlineto
-2 2 rlineto
-2 -2 rlineto
closepath fill
showpage
\end{tobiwrite}
\begin{tobiwrite}{foliation2.eps}
/arcs {gsave rotate
newpath
0 0 moveto
50 0 lineto
1 setlinewidth
stroke
newpath
0 8 moveto
50 0 rlineto
0 -16 rlineto
-50 0 rlineto
closepath
clip
newpath
r0 5 mul 2 add r0 50{/r exch def
r 0 moveto
20 r sub 0 r 2 mul 20 sub 0 360 arc
} for
0.5 setlinewidth stroke
newpath
r0 2 add r0 20{/r exch def
r 0 moveto 0 0 r 0 360 arc}for
0.5 setlinewidth stroke
grestore} def
50 50 translate

/r0 4.2 def

0 arcs
95 arcs
170 arcs
-97 arcs

newpath
2.5 0 moveto
0 0 2.5 0 360 arc
gsave 1 setgray
fill
grestore 0 setgray
stroke
showpage
\end{tobiwrite}
\begin{tobiwrite}{cleanup15.eps}
newpath
0 50 moveto
5 40 25 30 35 30 curveto stroke
newpath 0 30 moveto 80 0 rlineto
30 10 moveto
35 20 55 30 65 30 curveto
2 setlinewidth stroke
showpage
\end{tobiwrite}
\begin{tobiwrite}{thicktrack1.eps}
/v {gsave translate newpath
3 0 moveto 0 0 3 0 360 arc gsave
1 setgray fill grestore 0 setgray 1 setlinewidth
stroke grestore} def

gsave

50 50 translate
newpath
0 0 moveto
50 0 rlineto
110 rotate
0 0 moveto
50 0 rlineto

100 rotate
0 0 moveto
50 0 rlineto

70 rotate
0 0 moveto
50 0 rlineto
stroke
grestore
50 50 v

/Times-Roman findfont 12 scalefont setfont

80 39 moveto
(2) show
41 80 moveto
(3) show
17 37 moveto
(1) show
59 13 moveto
(4) show
showpage
\end{tobiwrite}
\begin{tobiwrite}{spiral1.eps}
newpath
0 120 moveto
10 110 10 100 10 90 curveto
10 50 lineto
10 -10 100 -10 100 60 curveto
100 130 10 130 10 80 curveto
10 70 moveto
10 60 10 50 25 35 curveto
stroke
showpage
\end{tobiwrite}
\begin{tobiwrite}{cleanup1.eps}
newpath
0 10 moveto
5 20 25 30 35 30 curveto
80 50 moveto
75 40 55 30 45 30 curveto stroke
newpath 0 30 moveto 80 0 rlineto 2 setlinewidth stroke
showpage
\end{tobiwrite}
\begin{tobiwrite}{case23.eps}
newpath
40 0 moveto
50 10 50 20 50 30 curveto
50 120 lineto
50 40 moveto
50 50 50 60 40 70 curveto
60 0 moveto
50 10 50 20 50 30 curveto
40 80 moveto
50 90 50 100 50 110 curveto
stroke
showpage
\end{tobiwrite}
\begin{tobiwrite}{tight1.eps}
newpath
0 0 moveto
30 60 70 60 100 0 curveto
0 60 moveto
30 0 70 0 100 60 curveto
stroke
showpage
\end{tobiwrite}
\begin{tobiwrite}{cleanup12.eps}
newpath
15 50 moveto
20 40 40 40 50 40 curveto
stroke
newpath 0 30 moveto 50 0 rlineto
60 30 75 38 80 45 curveto
50 30 moveto
60 30 75 22 80 15 curveto
2 setlinewidth stroke
newpath
50 40 moveto
50 40 75 38 80 45 curveto
50 40 moveto
70 40 75 22 80 15 curveto
{4 2} 0 setdash 1 setlinewidth stroke
showpage
\end{tobiwrite}
\begin{tobiwrite}{cleanup6.eps}
newpath
0 50 moveto
5 40 25 40 55 40 curveto
30 50 moveto
35 40 50 40 55 40 curveto
65 40 65 30 75 30 curveto stroke
newpath 0 30 moveto 80 0 rlineto 2 setlinewidth stroke
showpage
\end{tobiwrite}
\begin{tobiwrite}{tight6.eps}
newpath
0 0 moveto
27 54 15 11 48 11 curveto
73 11 70 27 80 27 curveto
90 27 112 -25 160 60 curveto
0 60 moveto
27 6 15 49 48 49 curveto
73 49 70 33 80 33 curveto
90 33 112 85 160 0 curveto
stroke
showpage
\end{tobiwrite}
\begin{tobiwrite}{case232a.eps}
newpath
40 0 moveto
50 10 50 15 50 25 curveto
50 120 lineto
50 30 moveto
50 40 50 45 40 55 curveto
40 65 moveto
50 75 50 80 50 90 curveto

80 0 moveto
70 10 70 15 70 25 curveto
70 80 lineto
70 95 50 95 50 110 curveto
80 50 moveto
70 60 70 70 70 80 curveto
stroke
showpage
\end{tobiwrite}
\begin{tobiwrite}{denorm5.eps}
/v {gsave translate
4 0 moveto 0 0 4 0 360 arc gsave
1 setgray fill grestore 0 setgray
2 setlinewidth stroke
grestore} def

gsave
newpath
15 60 moveto
30 80 lineto
33 84 45 87 50 85 curveto
75 75 lineto
80 73 85 65 85 62 curveto
85 38 lineto
85 35 80 27 75 25 curveto
50 15 lineto
45 13 33 16 30 20 curveto
15 40 lineto
12 42 12 58 15 60 curveto
closepath
50 50 moveto
65 60 85 70 85 55 curveto
85 45 moveto
85 30 105 40 120 50 curveto
85 38 moveto
85 33 91 27 95 25 curveto
115 15 lineto
117 14 123 14 125 15 curveto
145 25 lineto
149 27 155 33 155 38 curveto
155 62 lineto
155 67 149 73 145 75 curveto
125 85 lineto
123 86 117 86 115 85 curveto
95 75 lineto
91 73 85 67 85 62 curveto
15 60 moveto
12 56 7 53 5 52 curveto
15 40 moveto
12 44 7 47 5 48 curveto
30 80 moveto
33 84 35 90 37 95 curveto
30 20 moveto
33 16 35 10 37 5 curveto
50 85 moveto
45 87 43 92 40 95 curveto
50 15 moveto
45 13 43 8 40 5 curveto
75 75 moveto
80 73 91 73 95 75 curveto
75 25 moveto
80 27 91 27 95 25 curveto
115 15 moveto
117 14 117 8 118 5 curveto
125 15 moveto
123 14 123 8 122 5 curveto
115 85 moveto
117 86 117 92 118 95 curveto
125 85 moveto
123 86 123 92 122 95 curveto
155 38 moveto
155 35 160 31 163 28 curveto
155 62 moveto
155 65 160 69 163 72 curveto
145 25 moveto
149 27 155 27 160 24 curveto
145 75 moveto
149 73 155 73 160 76 curveto
.5 setlinewidth {2.5 2.5} 0 setdash .7 setgray stroke
grestore

5 20 moveto 50 50 lineto
5 80 lineto
50 50 moveto
165 50 lineto
65 5 moveto 50 50 lineto
65 95 lineto
100 95 moveto 120 50 lineto
100 5 lineto
145 95 moveto
120 50 lineto
145 5 lineto
2 setlinewidth stroke

newpath
50 50 moveto
70 63 83 63 83 50 curveto
85 30 80 30 60 20 curveto
50 15 34 15 24 30 curveto
14 45 14 55 24 70 curveto
34 85 50 85 60 80 curveto
70 75 87 70 87 50 curveto
87 37 100 37 120 50 curveto
1 setlinewidth stroke

50 50 v
120 50 v

showpage
\begin{tobiwrite}{norm6.eps}
0 111 translate
-90 rotate
6 93 translate
/v {gsave translate newpath
4 0 moveto 0 0 4 0 360 arc gsave
1 setgray fill grestore 0 setgray
stroke grestore} def
gsave
-60 rotate
newpath
0 0 moveto
30 20 70 20 100 0 curveto
1 setlinewidth
stroke
grestore
2 setlinewidth
newpath
0 0 moveto
100 0 lineto
50 50 3 sqrt mul lineto
closepath
50 -50 3 sqrt mul lineto
100 0 lineto
stroke
0 0 v
100 0 v
50 50 3 sqrt mul v
50 -50 3 sqrt mul v
-80 rotate
newpath
10 0 moveto
0 0 10 0 160 arc
1 setlinewidth stroke
160 rotate
newpath
10 3 moveto
2 -5 rlineto
-2 2 rlineto
-2 -2 rlineto
closepath fill
showpage
\end{tobiwrite}
\begin{tobiwrite}{circle2.eps}
50 50 translate
newpath
0 0 moveto
30 20 40 20 40 0 curveto
0 0 40 0 180 arcn
2.5 0 42.5 180 -6 arcn
stroke
2 setlinewidth
newpath
0 0 moveto
50 0 lineto
70 rotate
0 0 moveto
50 0 lineto
80 rotate
0 0 moveto
50 0 lineto
60 rotate
0 0 moveto
50 0 lineto
85 rotate
0 0 moveto
50 0 lineto
2 setlinewidth stroke
newpath
4 0 moveto
0 0 4 0 360 arc
gsave 1 setgray fill
grestore 0 setgray stroke
showpage
\end{tobiwrite}
\begin{tobiwrite}{cleanup14.eps}
newpath
0 50 moveto
5 40 25 50 45 36.5 curveto stroke
newpath 0 30 moveto 80 0 rlineto
30 50 moveto
35 40 55 30 65 30 curveto
2 setlinewidth stroke
showpage
\end{tobiwrite}
\begin{tobiwrite}{narc1.eps}
6 8 translate

/v {gsave translate newpath
4 0 moveto 0 0 4 0 360 arc gsave
1 setgray fill grestore 0 setgray
stroke grestore} def

2 setlinewidth

newpath
0 0 moveto
100 0 lineto
50 50 3 sqrt mul lineto
closepath stroke

0 0 v
100 0 v
50 50 3 sqrt mul v

newpath
-10 cos 40 mul
-10 sin 40 mul moveto
0 0 40 -10 70 arc
1 setlinewidth

stroke

showpage
\end{tobiwrite}
\begin{tobiwrite}{case234a.eps}
-30 0 translate
newpath
40 15 moveto
50 25 50 30 50 40 curveto
50 45 lineto
50 55 50 60 60 70 curveto
60 80  moveto
50 90 50 95 50 105 curveto
50 125 90 125 90 60 curveto
90 10 80 5 50 5 curveto
40 5 30 5 40 15 curveto
60 15 moveto
50 25 50 30 50 40 curveto
50 45 moveto
50 55 50 60 40 70 curveto
40 80 moveto
50 90 50 95 50 105 curveto
stroke
showpage
\end{tobiwrite}
\begin{tobiwrite}{norm8.eps}
0 111 translate
-90 rotate
6 93 translate
/v {gsave translate newpath
4 0 moveto 0 0 4 0 360 arc gsave
1 setgray fill grestore 0 setgray
stroke grestore} def
gsave
-60 rotate
newpath
0 0 moveto
30 20 50 20 50 -4 curveto
1 setlinewidth
stroke
grestore
2 setlinewidth
newpath
0 0 moveto
100 0 lineto
50 50 3 sqrt mul lineto
closepath
50 -50 3 sqrt mul lineto
100 0 lineto
stroke
0 0 v
100 0 v
50 50 3 sqrt mul v
50 -50 3 sqrt mul v
-80 rotate
newpath
10 0 moveto
0 0 10 0 160 arc
1 setlinewidth stroke
160 rotate
newpath
10 3 moveto
2 -5 rlineto
-2 2 rlineto
-2 -2 rlineto
closepath fill
showpage
\end{tobiwrite}
\begin{tobiwrite}{cleanup13.eps}
newpath
0 50 moveto
5 40 25 30 35 30 curveto stroke
newpath 0 30 moveto 80 0 rlineto
30 50 moveto
35 40 55 30 65 30 curveto
2 setlinewidth stroke
showpage
\end{tobiwrite}
\begin{tobiwrite}{triang2track1.eps}
/v {gsave translate newpath 2 0 moveto
0 0 2 0 360 arc
gsave 1 setgray fill grestore
0 setgray 0.7 setlinewidth stroke grestore} def
newpath
70 100 moveto
150 45 lineto
150 155 lineto
230 100 lineto
150 45 lineto
150 155 moveto
70 100 lineto
10 50 rlineto
70 100 moveto
-40 30 rlineto
70 100 moveto
-40 -30 rlineto
70 100 moveto
10 -50 rlineto
230 100 moveto
0 60 rlineto
230 100 moveto
50 30 rlineto
230 100 moveto
0 -60 rlineto
230 100 moveto
50 -30 rlineto
150 45 moveto
-30 -40 rlineto
150 45 moveto
30 -40 rlineto
150 155 moveto
-50 25 rlineto
150 155 moveto
50 25 rlineto
150 155 moveto
0 45 rlineto
0.5 setlinewidth
[5 5] 2 setdash
stroke
newpath
110 60 moveto
115 80 lineto
120 100 140 100 135 100 curveto
165 100 lineto
175 100 180 95 185 85 curveto
200 55 lineto
205 45 200 20 190 15 curveto
170 5 lineto
160 0 143 -1 135 5 curveto
115 20 lineto
107 26 105 40 110 60 curveto
115 80 moveto
120 100 122 104 120 110 curveto
110 140 lineto
108 146 98 148 94 148 curveto
62 148 lineto
56 148 52 141 48 135 curveto
32 111 lineto
28 105 30 91 35 86 curveto
51 70 lineto
54 67 60 62 65 60 curveto
90 50 lineto
100 46 105 40 110 60 curveto
135 100 moveto
125 100 123 101 120 110 curveto
165 100 moveto
180 100 186 105 188 110 curveto
200 140 lineto
202 145 205 152 200 163 curveto
185 190 lineto
182.5 194.5 175 195 165 195 curveto
135 195 lineto
130 195 123 194 120 185 curveto
112 161 lineto
110 155 108 146 110 140 curveto
185 85 moveto
182 91 186 105 188 110 curveto
200 140 moveto
204 150 210 155 220 155 curveto
240 155 lineto
250 155 260 145 265 135 curveto
275 115 lineto
280 105 277 90 275 85 curveto
265 60 lineto
261 50 255 48 245 48 curveto
220 48 lineto
213 48 205 45 200 55 curveto
115 20 moveto
95 35 100 46 90 50 curveto
65 60 moveto
60 62 58 62 56 60 curveto
51 70 moveto
53 67 53 65 53 62 curveto
35 86 moveto
33 89 31 90 28 91 curveto
32 111 moveto
30 108 28 106 27 105 curveto
48 135 moveto
50 138 52 143 52 148 curveto
62 148 moveto
60 148 58 150 56 152 curveto
94 148 moveto
100 148 110 155 112 161 curveto
120 185 moveto
121 188 122 192 122 195 curveto
135 195 moveto
133 195 129 196  127 199 curveto
165 195 moveto
168 195 170 196 174 199 curveto
185 190 moveto
180 199 182 199 182 199 curveto
200 163 moveto
205 152 210 155 220 155 curveto
240 155 moveto
245 155 247 156 250 158 curveto
265 135 moveto
260 145 257 151 256 155 curveto
275 115 moveto
280 105 283 103 285 102 curveto
275 85 moveto
279 91 282 95 285 96 curveto
265 60 moveto
263 55 262 50 262 48 curveto
245 48 moveto
250 48 252 47 254 45 curveto
220 48 moveto
205 48 210 25 190 15 curveto
170 5 moveto
164 2 163 1 162 0 curveto
135 5 moveto
139 2 140 1 141 0 curveto
70 100 moveto
60 70 51 70 43 78 curveto
150 45 moveto
125 90 115 80 112.5 70 curveto
150 155 moveto
195 150 200 140 194.5 126 curveto
230 100 moveto
210 55 200 55 191 73 curveto
[] 0 setdash
stroke
70 100 v
150 45 v
150 155 v
230 100 v
showpage
\end{tobiwrite}
\begin{tobiwrite}{norm4.eps}
0 111 translate
-90 rotate
6 93 translate
/v {gsave translate newpath
4 0 moveto 0 0 4 0 360 arc gsave
1 setgray fill grestore 0 setgray
stroke grestore} def
newpath
0 0 moveto
/a 50 3 sqrt mul def
a 4 add 30 cos mul
a 4 add -30 sin mul
lineto
1 setlinewidth
stroke
2 setlinewidth
newpath
0 0 moveto
100 0 lineto
50 50 3 sqrt mul lineto
closepath
50 -50 3 sqrt mul lineto
100 0 lineto
stroke
0 0 v
100 0 v
50 50 3 sqrt mul v
50 -50 3 sqrt mul v
-80 rotate
newpath
10 0 moveto
0 0 10 0 160 arc
1 setlinewidth stroke
160 rotate
newpath
10 3 moveto
2 -5 rlineto
-2 2 rlineto
-2 -2 rlineto
closepath fill
showpage
\end{tobiwrite}
\begin{tobiwrite}{norm5.eps}
0 111 translate
-90 rotate
6 93 translate
/v {gsave translate newpath
4 0 moveto 0 0 4 0 360 arc gsave
1 setgray fill grestore 0 setgray
stroke grestore} def
newpath
0 0 moveto
30 20 50 20 50 0 curveto
gsave
closepath
0.8 setgray fill
0 setgray grestore
50 -50 3 sqrt mul lineto
1 setlinewidth
stroke
2 setlinewidth
newpath
0 0 moveto
100 0 lineto
50 50 3 sqrt mul lineto
closepath
50 -50 3 sqrt mul lineto
100 0 lineto
stroke
0 0 v
100 0 v
50 50 3 sqrt mul v
50 -50 3 sqrt mul v
-80 rotate
newpath
10 0 moveto
0 0 10 0 160 arc
1 setlinewidth stroke
160 rotate
newpath
10 3 moveto
2 -5 rlineto
-2 2 rlineto
-2 -2 rlineto
closepath fill
showpage
\end{tobiwrite}
\begin{tobiwrite}{sing1.eps}

newpath
50 50 moveto
50 50 2 0 360 arc fill

newpath
50 50 moveto
50 0 rlineto stroke

50 50 translate
newpath 0 50 moveto
0 0 50 0 360 arc clip

newpath
1.1 1.1 7{
/c exch def
/b 50 c c mul add sqrt c mul 2 mul def
/d 50 b b mul c div c div 3 div sub def
50 b moveto
d b 3 div
d b 3 div neg
50 b neg curveto
} for
stroke
showpage
\end{tobiwrite}
\begin{tobiwrite}{case233.eps}
newpath
40 0 moveto
50 10 50 15 50 25 curveto
50 95 lineto
50 105 50 110 40 120 curveto
60 0 moveto
50 10 50 15 50 25 curveto
50 30 moveto
50 40 50 45 40 55 curveto
40 65 moveto
50 75 50 80 50 90 curveto
50 95 moveto
50 105 50 110 60 120 curveto
stroke
showpage
\end{tobiwrite}
\begin{tobiwrite}{denorm1.eps}
/v {gsave translate
4 0 moveto 0 0 4 0 360 arc gsave
1 setgray fill grestore 0 setgray
2 setlinewidth stroke
grestore} def

gsave
newpath
15 60 moveto
30 80 lineto
33 84 45 87 50 85 curveto
75 75 lineto
80 73 85 65 85 62 curveto
85 38 lineto
85 35 80 27 75 25 curveto
50 15 lineto
45 13 33 16 30 20 curveto
15 40 lineto
12 42 12 58 15 60 curveto
closepath
50 50 moveto
65 60 85 70 85 55 curveto
85 45 moveto
85 30 105 40 120 50 curveto
85 38 moveto
85 33 91 27 95 25 curveto
115 15 lineto
117 14 123 14 125 15 curveto
145 25 lineto
149 27 155 33 155 38 curveto
155 62 lineto
155 67 149 73 145 75 curveto
125 85 lineto
123 86 117 86 115 85 curveto
95 75 lineto
91 73 85 67 85 62 curveto
15 60 moveto
12 56 7 53 5 52 curveto
15 40 moveto
12 44 7 47 5 48 curveto
30 80 moveto
33 84 35 90 37 95 curveto
30 20 moveto
33 16 35 10 37 5 curveto
50 85 moveto
45 87 43 92 40 95 curveto
50 15 moveto
45 13 43 8 40 5 curveto
75 75 moveto
80 73 91 73 95 75 curveto
75 25 moveto
80 27 91 27 95 25 curveto
115 15 moveto
117 14 117 8 118 5 curveto
125 15 moveto
123 14 123 8 122 5 curveto
115 85 moveto
117 86 117 92 118 95 curveto
125 85 moveto
123 86 123 92 122 95 curveto
155 38 moveto
155 35 160 31 163 28 curveto
155 62 moveto
155 65 160 69 163 72 curveto
145 25 moveto
149 27 155 27 160 24 curveto
145 75 moveto
149 73 155 73 160 76 curveto
.5 setlinewidth {2.5 2.5} 0 setdash .7 setgray stroke
grestore

5 20 moveto 50 50 lineto
5 80 lineto
50 50 moveto
165 50 lineto
65 5 moveto 50 50 lineto
65 95 lineto
100 95 moveto 120 50 lineto
100 5 lineto
145 95 moveto
120 50 lineto
145 5 lineto
2 setlinewidth stroke

newpath
50 50 moveto
70 60 100 60 120 50 curveto
1 setlinewidth stroke

50 50 v
120 50 v

showpage
\begin{tobiwrite}{case21.eps}
newpath
40 0 moveto
50 10 50 20 50 30 curveto
50 90 lineto
50 100 50 110 40 120 curveto
60 0 moveto
50 10 50 20 50 30 curveto
50 90 moveto
50 100 50 110 60 120 curveto
50 60 moveto
50 60 50 70 40 80 curveto
stroke
showpage
\end{tobiwrite}
\begin{tobiwrite}{case232.eps}
newpath
40 0 moveto
50 10 50 15 50 25 curveto
50 120 lineto
60 0 moveto
50 10 50 15 50 25 curveto
50 30 moveto
50 40 50 45 40 55 curveto
40 65 moveto
50 75 50 80 50 90 curveto
60 85 moveto
50 95 50 100 50 110 curveto
stroke
showpage
\end{tobiwrite}
\begin{tobiwrite}{narc3.eps}

6 8 translate

/v {gsave translate newpath
4 0 moveto 0 0 4 0 360 arc gsave
1 setgray fill grestore 0 setgray 2 setlinewidth
stroke grestore} def

2 setlinewidth

newpath
0 0 moveto
100 0 lineto
50 50 3 sqrt mul lineto
closepath stroke

1 setlinewidth

newpath
0 0 moveto
40 25 60 25 100 0 curveto
stroke

0 0 v
100 0 v
50 50 3 sqrt mul v

showpage
\end{tobiwrite}
\begin{tobiwrite}{spiral3.eps}
newpath
0 120 moveto
10 110 10 100 10 90 curveto
10 50 lineto
10 -10 100 -10 100 60 curveto
100 130 10 130 10 80 curveto
10 70 moveto
10 0 90 0 90 60 curveto
90 120 20 120 20 70 curveto
20 60 20 50 35 35 curveto
stroke
showpage
\end{tobiwrite}
\begin{tobiwrite}{split5.eps}
newpath
0 0 moveto
20 30 80 30 100 0 curveto
0 60 moveto
20 30 80 30 100 60 curveto
0 0 moveto
20 30 80 30 100 60 curveto
0 0 moveto
20 30 80 30 100 60 curveto
0 60 moveto
20 30 80 30 100 0 curveto
stroke
showpage
\end{tobiwrite}
\begin{tobiwrite}{thicktrack2.eps}
/v {gsave translate newpath
3 0 moveto 0 0 3 0 360 arc gsave
1 setgray fill grestore 0 setgray 1 setlinewidth
stroke grestore} def

gsave

50 50 translate
newpath
0 0 moveto
10 2
20 2
50 2
50 2 curveto
0 0 moveto
10 -2
20 -2
50 -2
50 -2 curveto

110 rotate
0 0 moveto
10 -4
20 -4
50 -4
50 -4 curveto
0 0 moveto
10 4
20 4
50 4
50 4 curveto
0 0 moveto
50 0 rlineto

100 rotate
0 0 moveto
50 0 rlineto

70 rotate
0 0 moveto
10 6
20 6
50 6
50 6 curveto
0 0 moveto
10 2
20 2
50 2
50 2 curveto
0 0 moveto
10 -2
20 -2
50 -2
50 -2 curveto
0 0 moveto
10 -6
20 -6
50 -6
50 -6 curveto
stroke
grestore
50 50 v
showpage
\end{tobiwrite}
\begin{tobiwrite}{arctypes.eps}
1.5 1.5 scale
6 8 translate

/v {gsave translate newpath
3 0 moveto 0 0 3 0 360 arc gsave
1 setgray fill grestore 0 setgray 1.33 setlinewidth
stroke grestore} def

/a {newpath
0 0 moveto
3 sqrt 47 mul
47 lineto
stroke
newpath
-10 cos 35 mul
-10 sin 35 mul moveto
0 0 35 -10 70 arc
0.67 setlinewidth
stroke} def

1.33 setlinewidth
newpath
0 0 moveto
100 0 lineto
50 50 3 sqrt mul lineto
closepath stroke

0.67 setlinewidth
a
gsave
100 0 translate
120 rotate
a
100 0 translate
120 rotate
a
grestore

0 0 v
100 0 v
50 50 3 sqrt mul v
showpage
\end{tobiwrite}
\begin{tobiwrite}{norm2.eps}
0 111 translate
-90 rotate
6 93 translate
/v {gsave translate newpath
4 0 moveto 0 0 4 0 360 arc gsave
1 setgray fill grestore 0 setgray
stroke grestore} def
newpath
0 0 moveto
30 -20 70 -20 100 0 curveto
1 setlinewidth
stroke
2 setlinewidth
newpath
0 0 moveto
100 0 lineto
50 50 3 sqrt mul lineto
closepath
50 -50 3 sqrt mul lineto
100 0 lineto
stroke
0 0 v
100 0 v
50 50 3 sqrt mul v
50 -50 3 sqrt mul v
-80 rotate
newpath
10 0 moveto
0 0 10 0 160 arc
1 setlinewidth stroke
160 rotate
newpath
10 3 moveto
2 -5 rlineto
-2 2 rlineto
-2 -2 rlineto
closepath fill
showpage
\end{tobiwrite}
\begin{tobiwrite}{2bigons1.eps}
newpath
99 50 moveto
50 50 49 0 360 arc
129 50 moveto
80 50 49 0 360 arc
stroke
showpage
\end{tobiwrite}
\begin{tobiwrite}{spiral4.eps}
0 120 moveto
10 110 10 100 10 90 curveto
10 50 lineto
10 -10 110 -10 110 50 curveto
110 90 lineto
110 140 10 140 10 90 curveto
10 70 moveto
10 60 10 50 30 30 curveto
120 120 moveto
110 110 110 100 110 90 curveto
110 70 moveto
110 60 110 50 90 30 curveto
stroke
showpage
\end{tobiwrite}
\begin{tobiwrite}{corner1.eps}

100 50 moveto
50 50 50 0 360 arc clip

gsave
newpath
0 0 moveto
50 50 lineto 100 0 lineto
closepath .9 setgray fill
grestore

newpath
0 0 moveto
100 100 lineto
100 0 moveto
0 100 lineto
2 setlinewidth stroke

newpath
5 7 100 {0 exch moveto
100 0 rlineto} for
.5 setlinewidth stroke

50 50 translate
newpath
0 0 moveto
0 0 4 0 360 arc fill

showpage
\end{tobiwrite}
\begin{tobiwrite}{case234.eps}
newpath
40 15 moveto
50 25 50 30 50 40 curveto
50 105 lineto
50 125 90 125 90 60 curveto
90 10 80 5 50 5 curveto
40 5 30 5 40 15 curveto
60 15 moveto
50 25 50 30 50 40 curveto
50 45 moveto
50 55 50 60 40 70 curveto
40 80 moveto
50 90 50 95 50 105 curveto

stroke
newpath
90 60 translate
7 0 moveto
0 0 7 0 360 arc
gsave 1 setgray fill grestore 0 setgray stroke
45 rotate
newpath
7 0 moveto -7 0 lineto
0 7 moveto 0 -7 lineto
stroke
showpage
\end{tobiwrite}
\begin{tobiwrite}{univ2.eps}
6 8 translate

/v {gsave translate newpath
4 0 moveto 0 0 4 0 360 arc gsave
1 setgray fill grestore 0 setgray
stroke grestore} def
/r 50 3 sqrt div def

2 setlinewidth
newpath
0 0 moveto
100 0 lineto
50 50 3 sqrt mul lineto
closepath stroke

newpath
57 0 moveto
-14 4 rlineto
5 -4 rlineto
-5 -4 rlineto
closepath fill
gsave
newpath
100 0 moveto
50 r translate
25 r 2 div neg
0 r 1.2 mul neg
0 r 2 mul 3 div neg curveto
0 r 3 div neg
r 3 div dup 30 cos mul exch 30 sin mul
r 2 mul 3 div dup 30 cos mul exch 30 sin mul curveto
r dup 30 cos mul exch 30 sin mul lineto
r 2 mul 3 div dup 30 cos mul exch 30 sin mul moveto
r 3 div dup 30 cos mul exch 30 sin mul
r 3 div dup 150 cos mul exch 150 sin mul
r 2 mul 3 div dup 150 cos mul exch 150 sin mul curveto
r dup 150 cos mul exch 150 sin mul lineto
r 2 mul 3 div dup 150 cos mul exch 150 sin mul moveto
r 3 div dup 150 cos mul exch 150 sin mul
0 r 3 div neg
0 r 2 mul 3 div neg curveto
1 setlinewidth stroke
grestore
0 0 v
100 0 v
50 50 3 sqrt mul v
showpage
\end{tobiwrite}
\begin{tobiwrite}{spiral2.eps}
newpath
0 120 moveto
10 110 10 100 10 90 curveto
10 50 lineto
10 -10 100 -10 100 60 curveto
100 130 20 130 20 80 curveto
20 70 lineto
20 60 20 50 35 35 curveto
10 55 moveto
10 65 20 65 20 75 curveto
stroke
showpage
\end{tobiwrite}
\begin{tobiwrite}{split2.eps}
newpath
0 0 moveto
20 30 80 30 100 0 curveto
0 60 moveto
20 30 80 30 100 60 curveto
0 0 moveto
20 30 80 30 100 60 curveto
stroke
showpage
\end{tobiwrite}
\begin{tobiwrite}{tight8.eps}
newpath
3 30 moveto
3 30 3 0 360 arc fill
newpath
3 30 moveto
30 90 70 9 90 9 curveto
100 9 110 9 140 60 curveto
3 30 moveto
30 -30 70 51 90 51 curveto
100 51 110 51 140 0 curveto
stroke
showpage
\end{tobiwrite}
\begin{tobiwrite}{cleanup4.eps}
newpath
0 50 moveto
5 40 40 40 50 40 curveto
60 40 75 40 80 50 curveto stroke
newpath 0 30 moveto 80 0 rlineto 2 setlinewidth stroke
newpath
10 30 moveto
30 30 30 40 50 40 curveto
{4 2} 0 setdash 1 setlinewidth stroke
showpage
\end{tobiwrite}
\begin{tobiwrite}{split4.eps}
newpath
0 10 moveto
10 20 85 27 96 30 curveto
0 50 moveto
10 40 85 33 96 30 curveto
stroke
newpath

99 30 moveto
96 30 3 0 360 arc
gsave 1 setgray fill
grestore 0 setgray stroke
showpage
\end{tobiwrite}
\begin{tobiwrite}{case231.eps}
newpath
40 0 moveto
50 10 50 15 50 25 curveto
50 120 lineto
60 0 moveto
50 10 50 15 50 25 curveto
50 30 moveto
50 40 50 45 40 55 curveto
40 65 moveto
50 75 50 80 50 90 curveto
40 85 moveto
50 95 50 100 50 110 curveto
stroke
showpage
\end{tobiwrite}
\begin{tobiwrite}{tight4.eps}
newpath
3 30 moveto
3 30 3 0 360 arc fill
newpath 3 30 moveto
40 90 65 -5 80 60 curveto
3 30 moveto
40 -30 65 65 80 0 curveto
stroke
showpage
\end{tobiwrite}
\begin{tobiwrite}{case24b.eps}
newpath
10 120 moveto
20 110 20 100 20 90 curveto
20 30 lineto
20 20 20 10 30 0 curveto
20 60 moveto
20 20 90 10 90 70 curveto
90 130 20 120 20 80 curveto
stroke

gsave
newpath
90 70 translate
7 0 moveto
0 0 7 0 360 arc
gsave 1 setgray fill grestore 0 setgray stroke

45 rotate
newpath
7 0 moveto -7 0 lineto
0 7 moveto 0 -7 lineto
stroke
grestore
newpath
20 30 translate
7 0 moveto
0 0 7 0 360 arc
gsave 1 setgray fill grestore 0 setgray stroke

45 rotate
newpath
7 0 moveto -7 0 lineto
0 7 moveto 0 -7 lineto
stroke
showpage
\end{tobiwrite}
\begin{tobiwrite}{sing2.eps}
/f {/y exch def /x exch def
/r x x mul y y mul add 1 3 div exp def
y x atan 2 mul 3 div dup cos r mul exch sin r mul} def

/s {0 0 moveto
70 0 rlineto
80 80 400{/d exch def
-400 d f moveto
-399 10 400 {d f lineto} for
} for} def

newpath
50 50 moveto
50 50 2 0 360 arc fill

50 50 translate
newpath 0 50 moveto
0 0 50 0 360 arc clip

newpath
gsave
s
120 rotate
s
120 rotate
s stroke
showpage
\end{tobiwrite}
\begin{tobiwrite}{spiral5.eps}
0 120 moveto
10 110 10 100 10 90 curveto
10 50 lineto
10 -10 110 -10 110 50 curveto
110 90 lineto
110 140 10 140 10 90 curveto
10 70 moveto
10 60 10 50 30 30 curveto
105 90 moveto
105 50 lineto
105 -5 15 -5 15 50 curveto
15 90 lineto
15 135 105 135 105 90 curveto
120 120 moveto
105 105 105 100 105 90 curveto
105 70 moveto
105 60 105 45 90 30 curveto
stroke
showpage
\end{tobiwrite}
\begin{tobiwrite}{case3.eps}
newpath
40 20 moveto
50 30 50 40 50 50 curveto
50 70 lineto
50 80 50 90 40 100 curveto
60 20 moveto
50 30 50 40 50 50 curveto
50 70 moveto
50 80 50 90 70 110 curveto
stroke
70 110 translate
newpath
3 0 moveto
0 0 3 0 360 arc
gsave 1 setgray fill grestore 0 setgray stroke
showpage
\end{tobiwrite}
\begin{tobiwrite}{case233b.eps}
20 0 translate
newpath
40 0 moveto
50 10 50 15 50 25 curveto
50 95 lineto
50 105 50 110 60 120 curveto
60 0 moveto
50 10 50 15 50 25 curveto
50 30 moveto
50 40 50 45 40 55 curveto
20 55 moveto
30 65 30 70 30 80 curveto
30 95 lineto
30 105 30 110 20 120 curveto
50 80 moveto
50 95 30 80 30 95 curveto
stroke
showpage
\end{tobiwrite}
\begin{tobiwrite}{denorm6.eps}
/v {gsave translate
4 0 moveto 0 0 4 0 360 arc gsave
1 setgray fill grestore 0 setgray
2 setlinewidth stroke
grestore} def

gsave
newpath
15 60 moveto
30 80 lineto
33 84 45 87 50 85 curveto
75 75 lineto
80 73 85 65 85 62 curveto
85 38 lineto
85 35 80 27 75 25 curveto
50 15 lineto
45 13 33 16 30 20 curveto
15 40 lineto
12 42 12 58 15 60 curveto
closepath
50 50 moveto
65 60 85 70 85 55 curveto
85 45 moveto
85 30 105 40 120 50 curveto
85 38 moveto
85 33 91 27 95 25 curveto
115 15 lineto
117 14 123 14 125 15 curveto
145 25 lineto
149 27 155 33 155 38 curveto
155 62 lineto
155 67 149 73 145 75 curveto
125 85 lineto
123 86 117 86 115 85 curveto
95 75 lineto
91 73 85 67 85 62 curveto
15 60 moveto
12 56 7 53 5 52 curveto
15 40 moveto
12 44 7 47 5 48 curveto
30 80 moveto
33 84 35 90 37 95 curveto
30 20 moveto
33 16 35 10 37 5 curveto
50 85 moveto
45 87 43 92 40 95 curveto
50 15 moveto
45 13 43 8 40 5 curveto
75 75 moveto
80 73 91 73 95 75 curveto
75 25 moveto
80 27 91 27 95 25 curveto
115 15 moveto
117 14 117 8 118 5 curveto
125 15 moveto
123 14 123 8 122 5 curveto
115 85 moveto
117 86 117 92 118 95 curveto
125 85 moveto
123 86 123 92 122 95 curveto
155 38 moveto
155 35 160 31 163 28 curveto
155 62 moveto
155 65 160 69 163 72 curveto
145 25 moveto
149 27 155 27 160 24 curveto
145 75 moveto
149 73 155 73 160 76 curveto
.5 setlinewidth {2.5 2.5} 0 setdash .7 setgray stroke
grestore

5 20 moveto 50 50 lineto
5 80 lineto
50 50 moveto
165 50 lineto
65 5 moveto 50 50 lineto
65 95 lineto
100 95 moveto 120 50 lineto
100 5 lineto
145 95 moveto
120 50 lineto
145 5 lineto
2 setlinewidth stroke

newpath
50 50 moveto
110 10 125 10 120 50 curveto
1 setlinewidth stroke

50 50 v
120 50 v

showpage
\begin{tobiwrite}{denorm4.eps}
/v {gsave translate
4 0 moveto 0 0 4 0 360 arc gsave
1 setgray fill grestore 0 setgray
2 setlinewidth stroke
grestore} def

gsave
newpath
15 60 moveto
30 80 lineto
33 84 45 87 50 85 curveto
75 75 lineto
80 73 85 65 85 62 curveto
85 38 lineto
85 35 80 27 75 25 curveto
50 15 lineto
45 13 33 16 30 20 curveto
15 40 lineto
12 42 12 58 15 60 curveto
closepath
50 50 moveto
65 60 85 70 85 55 curveto
85 45 moveto
85 30 105 40 120 50 curveto
85 38 moveto
85 33 91 27 95 25 curveto
115 15 lineto
117 14 123 14 125 15 curveto
145 25 lineto
149 27 155 33 155 38 curveto
155 62 lineto
155 67 149 73 145 75 curveto
125 85 lineto
123 86 117 86 115 85 curveto
95 75 lineto
91 73 85 67 85 62 curveto
15 60 moveto
12 56 7 53 5 52 curveto
15 40 moveto
12 44 7 47 5 48 curveto
30 80 moveto
33 84 35 90 37 95 curveto
30 20 moveto
33 16 35 10 37 5 curveto
50 85 moveto
45 87 43 92 40 95 curveto
50 15 moveto
45 13 43 8 40 5 curveto
75 75 moveto
80 73 91 73 95 75 curveto
75 25 moveto
80 27 91 27 95 25 curveto
115 15 moveto
117 14 117 8 118 5 curveto
125 15 moveto
123 14 123 8 122 5 curveto
115 85 moveto
117 86 117 92 118 95 curveto
125 85 moveto
123 86 123 92 122 95 curveto
155 38 moveto
155 35 160 31 163 28 curveto
155 62 moveto
155 65 160 69 163 72 curveto
145 25 moveto
149 27 155 27 160 24 curveto
145 75 moveto
149 73 155 73 160 76 curveto
.5 setlinewidth {2.5 2.5} 0 setdash .7 setgray stroke
grestore

5 20 moveto 50 50 lineto
5 80 lineto
50 50 moveto
165 50 lineto
65 5 moveto 50 50 lineto
65 95 lineto
100 95 moveto 120 50 lineto
100 5 lineto
145 95 moveto
120 50 lineto
145 5 lineto
2 setlinewidth stroke

newpath
50 50 moveto
30 105 90 70 120 50 curveto
1 setlinewidth stroke

50 50 v
120 50 v

showpage
\begin{tobiwrite}{norm7.eps}
0 111 translate
-90 rotate
6 93 translate
/v {gsave translate newpath
4 0 moveto 0 0 4 0 360 arc gsave
1 setgray fill grestore 0 setgray
stroke grestore} def
newpath
0 0 moveto
30 20 50 20 50 0 curveto
gsave
closepath
0.8 setgray fill
0 setgray grestore
0 0 50 0 -64 arcn
1 setlinewidth
stroke
2 setlinewidth
newpath
0 0 moveto
100 0 lineto
50 50 3 sqrt mul lineto
closepath
50 -50 3 sqrt mul lineto
100 0 lineto
stroke
0 0 v
100 0 v
50 50 3 sqrt mul v
50 -50 3 sqrt mul v
-80 rotate
newpath
10 0 moveto
0 0 10 0 160 arc
1 setlinewidth stroke
160 rotate
newpath
10 3 moveto
2 -5 rlineto
-2 2 rlineto
-2 -2 rlineto
closepath fill
showpage
\end{tobiwrite}
\begin{tobiwrite}{sing3.eps}
50 30 translate
newpath
0 0 moveto
0 0 2 0 360 arc fill

newpath
-50 0 moveto
100 0 rlineto
1.5 setlinewidth
stroke
1 setlinewidth

newpath 50 0 moveto
0 0 50 0 180 arc closepath clip

1 12 100 {gsave
newpath
0 exch translate
-50 0 moveto
-50 3 div -50
50 3 div -50
50 0 curveto stroke
grestore} for
showpage
\end{tobiwrite}
\begin{tobiwrite}{case22.eps}
newpath
40 0 moveto
50 10 50 20 50 30 curveto
50 120 lineto
60 0 moveto
50 10 50 20 50 30 curveto
50 110 moveto
50 100 50 90 60 80 curveto
50 60 moveto
50 60 50 70 40 80 curveto
stroke
showpage
\end{tobiwrite}
\begin{tobiwrite}{cleanup16.eps}
newpath
45 50 moveto
50 40 70 30 80 30 curveto stroke
newpath 0 30 moveto 80 0 rlineto
30 10 moveto
35 20 55 30 65 30 curveto
2 setlinewidth stroke
showpage
\end{tobiwrite}
\begin{tobiwrite}{case24.eps}
newpath
40 30 moveto
50 40 50 50 50 60 curveto
50 80 lineto
50 90 50 100 40 110 curveto
60 30 moveto
50 40 50 50 50 60 curveto
50 80 moveto
50 90 50 100 60 110 curveto
70 120 90 110 90 60 curveto
90 10 70 5 40 5 curveto
20 5 30 20 40 30 curveto
stroke

newpath
90 60 translate
7 0 moveto
0 0 7 0 360 arc
gsave 1 setgray fill grestore 0 setgray stroke

45 rotate
newpath
7 0 moveto -7 0 lineto
0 7 moveto 0 -7 lineto
stroke
showpage
\end{tobiwrite}
\begin{tobiwrite}{anarc.eps}
6 8 translate
/v {gsave translate newpath
4 0 moveto 0 0 4 0 360 arc gsave
1 setgray fill grestore 0 setgray
stroke grestore} def
newpath
0 0 moveto
30 15 45 15 45 -2 curveto
1 setlinewidth
stroke
2 setlinewidth
newpath
0 0 moveto
100 0 lineto
50 50 3 sqrt mul lineto
closepath stroke
0 0 v
100 0 v
50 50 3 sqrt mul v
showpage
\end{tobiwrite}
\begin{tobiwrite}{denorm2.eps}
/v {gsave translate
4 0 moveto 0 0 4 0 360 arc gsave
1 setgray fill grestore 0 setgray
2 setlinewidth stroke
grestore} def

gsave
newpath
15 60 moveto
30 80 lineto
33 84 45 87 50 85 curveto
75 75 lineto
80 73 85 65 85 62 curveto
85 38 lineto
85 35 80 27 75 25 curveto
50 15 lineto
45 13 33 16 30 20 curveto
15 40 lineto
12 42 12 58 15 60 curveto
closepath
50 50 moveto
65 60 85 70 85 55 curveto
85 45 moveto
85 30 105 40 120 50 curveto
85 38 moveto
85 33 91 27 95 25 curveto
115 15 lineto
117 14 123 14 125 15 curveto
145 25 lineto
149 27 155 33 155 38 curveto
155 62 lineto
155 67 149 73 145 75 curveto
125 85 lineto
123 86 117 86 115 85 curveto
95 75 lineto
91 73 85 67 85 62 curveto
15 60 moveto
12 56 7 53 5 52 curveto
15 40 moveto
12 44 7 47 5 48 curveto
30 80 moveto
33 84 35 90 37 95 curveto
30 20 moveto
33 16 35 10 37 5 curveto
50 85 moveto
45 87 43 92 40 95 curveto
50 15 moveto
45 13 43 8 40 5 curveto
75 75 moveto
80 73 91 73 95 75 curveto
75 25 moveto
80 27 91 27 95 25 curveto
115 15 moveto
117 14 117 8 118 5 curveto
125 15 moveto
123 14 123 8 122 5 curveto
115 85 moveto
117 86 117 92 118 95 curveto
125 85 moveto
123 86 123 92 122 95 curveto
155 38 moveto
155 35 160 31 163 28 curveto
155 62 moveto
155 65 160 69 163 72 curveto
145 25 moveto
149 27 155 27 160 24 curveto
145 75 moveto
149 73 155 73 160 76 curveto
.5 setlinewidth {2.5 2.5} 0 setdash .7 setgray stroke
grestore

5 20 moveto 50 50 lineto
5 80 lineto
50 50 moveto
165 50 lineto
65 5 moveto 50 50 lineto
65 95 lineto
100 95 moveto 120 50 lineto
100 5 lineto
145 95 moveto
120 50 lineto
145 5 lineto
2 setlinewidth stroke

newpath
50 50 moveto
70 40 100 40 120 50 curveto
1 setlinewidth stroke

50 50 v
120 50 v

showpage
\begin{tobiwrite}{narc2.eps}

6 8 translate

/v {gsave translate newpath
4 0 moveto 0 0 4 0 360 arc gsave
1 setgray fill grestore 0 setgray 2 setlinewidth
stroke grestore} def

2 setlinewidth

newpath
0 0 moveto
100 0 lineto
50 50 3 sqrt mul lineto
closepath stroke

1 setlinewidth

newpath
0 0 moveto
3 sqrt 47 mul
47 lineto
stroke

0 0 v
100 0 v
50 50 3 sqrt mul v

showpage
\end{tobiwrite}
\begin{tobiwrite}{corner4.eps}

100 50 moveto
50 50 50 0 360 arc clip

gsave
newpath
-10 -10 moveto
30 30 lineto
45 45 55 45 70 30 curveto
110 -10 lineto
closepath gsave .9 setgray fill
grestore 2 setlinewidth stroke
grestore

newpath
5 7 100 {0 moveto
0 100 rlineto} for
.5 setlinewidth stroke

showpage
\end{tobiwrite}
\begin{tobiwrite}{or3.eps}
6 8 translate
/v {gsave translate newpath
4 0 moveto 0 0 4 0 360 arc gsave
1 setgray fill grestore 0 setgray
stroke grestore} def
/r 50 3 sqrt div def
2 setlinewidth
newpath
0 0 moveto
100 0 lineto
50 50 3 sqrt mul lineto
closepath stroke
newpath
57 0 moveto
-14 4 rlineto
5 -4 rlineto
-5 -4 rlineto
closepath fill
gsave
newpath
100 0 moveto
50 r translate
25 r 2 div neg
0 r 1.2 mul neg
0 r 2 mul 3 div neg curveto
0 r 3 div neg
r 3 div dup 30 cos mul exch 30 sin mul
r 2 mul 3 div dup 30 cos mul exch 30 sin mul curveto
r dup 30 cos mul exch 30 sin mul lineto
r 2 mul 3 div dup 30 cos mul exch 30 sin mul moveto
r 3 div dup 30 cos mul exch 30 sin mul
r 3 div dup 150 cos mul exch 150 sin mul
r 2 mul 3 div dup 150 cos mul exch 150 sin mul curveto
r dup 150 cos mul exch 150 sin mul lineto
r 2 mul 3 div dup 150 cos mul exch 150 sin mul moveto
r 3 div dup 150 cos mul exch 150 sin mul
0 r 3 div neg
0 r 2 mul 3 div neg curveto
1 setlinewidth stroke
grestore
0 0 v
100 0 v
50 50 3 sqrt mul v
100 0 translate
70 rotate
newpath
10 0 moveto
0 0 10 0 140 arc
1 setlinewidth stroke
140 rotate
newpath
10 3 moveto
2 -5 rlineto
-2 2 rlineto
-2 -2 rlineto
closepath fill
showpage
\end{tobiwrite}
\begin{tobiwrite}{case21b.eps}
newpath
40 0 moveto
50 10 50 20 50 30 curveto
50 90 lineto
50 100 50 110 40 120 curveto
50 40 moveto
50 40 50 50 40 60 curveto
80 0 moveto
70 10 70 20 70 30 curveto
70 90 lineto
70 100 70 110 80 120 curveto
70 50 moveto
70 70 50 70 50 90 curveto
stroke
showpage
\end{tobiwrite}
\begin{tobiwrite}{sing4.eps}
50 30 translate
newpath
0 0 moveto
0 0 2 0 360 arc fill

newpath
-50 0 moveto
100 0 rlineto
1.5 setlinewidth
stroke
1 setlinewidth

newpath 50 0 moveto
0 0 50 0 180 arc closepath clip

-32 10 100 {gsave
newpath
0 exch translate
-50 0 moveto
-50 3 div 50
50 3 div 50
50 0 curveto stroke
grestore} for
showpage
\end{tobiwrite}
\begin{tobiwrite}{split1.eps}
newpath
0 0 moveto
10 15 20 30 35 30 curveto
65 30 lineto
80 30 90 15 100 0 curveto
0 60 moveto
10 45 20 30 35 30 curveto
65 30 moveto
80 30 90 45 100 60 curveto
stroke
showpage
\end{tobiwrite}
\begin{tobiwrite}{cleanup5.eps}
newpath
0 50 moveto
5 40 25 30 35 30 curveto
30 50 moveto
35 40 55 30 65 30 curveto stroke
newpath 0 30 moveto 80 0 rlineto 2 setlinewidth stroke
showpage
\end{tobiwrite}
\begin{tobiwrite}{case24a.eps}
newpath
10 110 moveto
20 100 20 90 20 80 curveto
20 40 lineto
20 0 90 -20 90 60 curveto
90 120 40 120 40 80 curveto
40 70 40 60 50 50 curveto
40 80 moveto
40 60 20 60 20 40 curveto
stroke

newpath
90 60 translate
7 0 moveto
0 0 7 0 360 arc
gsave 1 setgray fill grestore 0 setgray stroke

45 rotate
newpath
7 0 moveto -7 0 lineto
0 7 moveto 0 -7 lineto
stroke
showpage
\end{tobiwrite}
\begin{tobiwrite}{tight7.eps}
160 0 translate
-1 1 scale
newpath
0 0 moveto
27 54 15 11 48 11 curveto
73 11 70 27 80 27 curveto
90 27 112 -25 160 60 curveto
0 60 moveto
27 6 15 49 48 49 curveto
73 49 70 33 80 33 curveto
90 33 112 85 160 0 curveto
stroke
showpage
\end{tobiwrite}
\begin{tobiwrite}{cleanup7.eps}
newpath
0 50 moveto
5 40 25 30 35 30 curveto
30 50 moveto
35 40 55 30 65 30 curveto
15 10 moveto
20 20 40 30 50 30 curveto stroke
newpath 0 30 moveto 80 0 rlineto 2 setlinewidth stroke
showpage
\end{tobiwrite}
\begin{tobiwrite}{case21a.eps}
newpath
40 0 moveto
50 10 50 20 50 30 curveto
50 90 lineto
50 100 50 110 40 120 curveto
50 60 moveto
50 60 50 70 40 80 curveto
80 0 moveto
70 10 70 20 70 30 curveto
70 90 lineto
70 100 70 110 80 120 curveto
50 30 moveto
50 50 70 50 70 70 curveto
stroke
showpage
\end{tobiwrite}
\begin{tobiwrite}{case22a.eps}
newpath
40 0 moveto
50 10 50 20 50 30 curveto
50 120 lineto
50 50 moveto
50 50 50 60 40 70 curveto
50 110 moveto
50 90 70 90 70 70 curveto
70 30 lineto
70 20 70 10 80 0 curveto
70 70 moveto
70 60 70 50 80 40 curveto
stroke
showpage
\end{tobiwrite}
\begin{tobiwrite}{foliation1.eps}

0 25 translate

newpath
0 25 moveto
50 25 lineto
60 25
80 30
97 47 curveto
50 25 moveto
60 25
80 20
97 3 curveto
stroke

newpath
0 17 moveto
45 17 lineto
55 17 75 15 90 0 curveto
100 10 lineto
90 20 80 24 77 25 curveto
80 26 90 30 100 40 curveto
90 50 lineto
75 35 55 33 45 33 curveto
0 33 lineto
closepath
clip

0.5 setlinewidth

newpath
2 4 50 {0 moveto 0 50 rlineto} for
2.1 2 25{1.15 exp /t exch def
50 t add 0 moveto
50 t 1.6 mul add 20
50 t 1.6 mul add 30
50 t add 50 curveto
} for
stroke
showpage
\end{tobiwrite}
\begin{tobiwrite}{diverge.eps}
newpath
0 0 moveto
20 22 40 30 50 30 curveto
100 30 lineto
0 60 moveto
20 40 40 30 50 30 curveto
stroke
showpage
\end{tobiwrite}
\begin{tobiwrite}{2bigons2.eps}
newpath
3 30 moveto
3 30 3 0 360 arc
117 30 moveto
117 30 3 0 360 arc fill
newpath
3 30 moveto
40 100 80 -40 117 30 curveto
3 30 moveto
40 -40 80 100 117 30 curveto
stroke
showpage
\end{tobiwrite}
\begin{tobiwrite}{univ1.eps}
6 8 translate
/v {gsave translate newpath
4 0 moveto 0 0 4 0 360 arc gsave
1 setgray fill grestore 0 setgray
stroke grestore} def
/r 50 3 sqrt div def

2 setlinewidth
newpath
0 0 moveto
100 0 lineto
50 50 3 sqrt mul lineto
closepath stroke
0 0 v
100 0 v
50 50 3 sqrt mul v

50 r translate
newpath
0 r neg moveto
0 r 2 mul 3 div neg lineto
0 r 3 div neg
r 3 div dup 30 cos mul exch 30 sin mul
r 2 mul 3 div dup 30 cos mul exch 30 sin mul curveto
r dup 30 cos mul exch 30 sin mul lineto
r 2 mul 3 div dup 30 cos mul exch 30 sin mul moveto
r 3 div dup 30 cos mul exch 30 sin mul
r 3 div dup 150 cos mul exch 150 sin mul
r 2 mul 3 div dup 150 cos mul exch 150 sin mul curveto
r dup 150 cos mul exch 150 sin mul lineto
r 2 mul 3 div dup 150 cos mul exch 150 sin mul moveto
r 3 div dup 150 cos mul exch 150 sin mul
0 r 3 div neg
0 r 2 mul 3 div neg curveto
1 setlinewidth stroke
showpage
\end{tobiwrite}
\begin{tobiwrite}{thicktrack4.eps}
newpath
0 52 moveto
30 52 lineto
40 52 50 52 100 62 curveto
0 56 moveto
30 56 lineto
40 56 50 56 100 66 curveto
0 48 moveto
30 48 lineto
40 48 50 48 100 38 curveto
0 44 moveto
30 44 lineto
40 44 50 44 100 34 curveto
0 40 moveto
30 40 lineto
40 40 50 40 100 30 curveto
0 36 moveto
30 36 lineto
40 36 50 36 100 26 curveto
0 32 moveto
30 32 lineto
40 32 50 32 100 22 curveto
stroke
showpage
\end{tobiwrite}
\begin{tobiwrite}{split2a.eps}
0 60 translate
1 -1 scale
newpath
0 0 moveto
20 30 80 30 100 0 curveto
0 60 moveto
20 30 80 30 100 60 curveto
0 0 moveto
20 30 80 30 100 60 curveto
stroke
showpage
\end{tobiwrite}
\begin{tobiwrite}{circle1.eps}
50 50 translate
newpath
0 0 moveto
30 20 45 20 45 -4 curveto
stroke
2 setlinewidth
newpath
0 0 moveto
50 0 lineto
70 rotate
0 0 moveto
50 0 lineto
80 rotate
0 0 moveto
50 0 lineto
60 rotate
0 0 moveto
50 0 lineto
85 rotate
0 0 moveto
50 0 lineto
2 setlinewidth stroke
newpath
4 0 moveto
0 0 4 0 360 arc
gsave 1 setgray fill
grestore 0 setgray stroke
showpage
\end{tobiwrite}
\begin{tobiwrite}{cleanup11.eps}
newpath
0 50 moveto
5 40 25 30 35 30 curveto
stroke
newpath 0 30 moveto 50 0 rlineto
60 30 75 38 80 45 curveto
50 30 moveto
60 30 75 22 80 15 curveto
2 setlinewidth stroke
showpage
\end{tobiwrite}
\begin{tobiwrite}{cleanup9.eps}
newpath
10 50 moveto
15 40 35 30 45 30 curveto
stroke
newpath 0 30 moveto 70 0 rlineto 2 setlinewidth stroke
70 30 translate
newpath
2.5 0 moveto 0 0 2.5 0 360 arc 1 setlinewidth gsave
1 setgray fill grestore 0 setgray stroke
showpage
\end{tobiwrite}
\begin{tobiwrite}{2edges.eps}
6 8 translate
/v {gsave translate newpath
4 0 moveto 0 0 4 0 360 arc gsave
1 setgray fill grestore 0 setgray
stroke grestore} def
newpath
100 90 moveto
100 80 lineto
100 60 80 40 60 30 curveto
40 20 40 10 100 0 curveto
100 80 moveto
100 60 120 40 140 30 curveto
200 0 lineto
60 30 moveto
80 40 120 40 140 30 curveto
stroke
newpath
0 0 moveto 200 0 lineto
200 120 0 120 0 0 curveto
2 setlinewidth stroke
0 0 v 100 0 v 200 0 v
newpath
57 0 moveto
-14 4 rlineto
5 -4 rlineto
-5 -4 rlineto
closepath fill
newpath
157 0 moveto
-14 4 rlineto
5 -4 rlineto
-5 -4 rlineto
closepath fill
showpage
\end{tobiwrite}
\begin{tobiwrite}{tight9.eps}
newpath
3 30 moveto
3 30 3 0 360 arc fill
newpath
3 30 moveto
30 80 52 33 62 33 curveto
72 33 77 51 92 51 curveto
122 51 110 5 140 60 curveto
3 30 moveto
30 -20 52 27 62 27 curveto
72 27 77 9 92 9 curveto
122 9 110 55 140 0 curveto
stroke
showpage
\end{tobiwrite}

\begin{document}

\noindent\hbox to\textwidth{\hss\it To the memory of Patrick Dehornoy}

\vskip1cm

\maketitle

\begin{abstract}
A fast algorithm for counting intersections of two normal curves on a triangulated
surface is proposed. It yields a convenient way for treating mapping class groups
of punctured surfaces by presenting mapping classes by matrices, and the composition by an
exotic matrix multiplication. An efficient solution of the word problem
for mapping class groups of punctured surfaces is proposed, with efficiency
understood in a more restrictive way than the most common one.
\end{abstract}

\section{Introduction}
Among all finitely generated
groups there are those that come with a naturally defined geometry due to
their geometric or topological origin. By a geometry here we mean a quasi-isometry
class of a left-invariant (or right-invariant) metric on the group. When the efficiency
of an algorithm solving some decision problem for such a group is discussed, it is natural to
evaluate it in terms of the accompanying geometry, which may be quite different
from the word-length geometry used as the default option for abstract finitely generated
groups.

Basic examples to look at are the groups~$\mathrm{GL}(n,\mathbb Z)$ of integral invertible
matrices, in which the `natural'
geometry is defined by the distance~$d(x,y)=\log\|x^{-1}y\|$, where~$\|\,\|$ stands
for an operator norm. The \emph{complexity}~$c(x)$ of an element~$x\in\mathrm{GL}(n,\mathbb Z)$
defined by~$c(x)=d(1,x)$ is asymptotically comparable to the amount of space needed for recording~$x$
in the conventional way. If~$x$ is conjugate to a Jordan block, then~$c(x^n)$ grows
with~$n$ as~$\log(n)$ whereas  the word length of~$x^n$ grows as~$n$.
This means that an algorithm operating with elements of~$\mathrm{GL}(n,\mathbb Z)$
which is polynomial-time with respect to the word-length geometry,
may appear to be exponential-time in worst cases with respect to the `natural' geometry. Such
a divergence is unavoidable if
the algorithm uses a presentation of the group elements as decompositions
into a product of generators and reads such presentations on input letter by letter.

Mapping class groups~$\mcg(M,\mathscr P)$ of punctured surfaces are similar to~$\mathrm{GL}(n,\mathbb Z)$
in this respect. There are several equivalent natural ways to define a geometry on them,
and the `natural' complexity of the $n$th power of a Dehn twist grows logarithmically with~$n$.
A polynomial-time solution to the word problem for these
groups with respect to the word-length geometry is given
by Lee Mosher in~\cite{mosher}, but it is not equally efficient in the worst case
with respect to the `natural' geometry, since the algorithm is based on finite state automata.

The present paper proposes a viewpoint on the groups~$\mcg(M,\mathscr P)$
from which these groups appear very much like integral matrix groups, and mimic,
in a certain sense, orthogonal groups.
The group elements are presented by specific integral~$N\times N$ matrices with~$N$ depending on the surface
and the number of punctures. The geometry on~$\mcg(M,\mathscr P)$ arising from this
presentation coincides with the one coming from the action on the thick part of the respective
Teichm\"uller space.

To define the matrix presentation,
a triangulation~$T$ is fixed on~$M$, and the elements of a certain subset~$L\subset{}\mathbb Z^N$
are interpreted as
multiple curves on the surface encoded by their normal coordinates with respect to~$T$.
For any two multiple curves~$\gamma_1,\gamma_2$,
one defines their geometric intersection index denoted
by~$\langle\gamma_1,\gamma_2\rangle$, which becomes a function on~$L\times L$ once the triangulation~$T$ has been fixed.

The matrices representing elements of~$\mcg(M,\mathscr P)$ are `orthogonal' with
respect to~$-\langle\,,\,\rangle$in the sense that the rows (equivalently, columns)
of those matrices form an orthonormal family provided that~$-\langle\,,\,\rangle$
is used instead of the
standard scalar product. The $ij$th element of the matrix
representing the product~$xy$ is equal to~$\langle r_i,c_j\rangle$,
where~$r_i$ is the $i$th row of the matrix representing~$x$, and~$c_j$ is the $j$th
column of the matrix representing~$y$. (However, not all integral matrices
satisfying the above mentioned `orthogonality' condition represent elements of~$\mcg(M,\mathscr P)$. What they do represent is
the set of isotopy classes of all triangulations of~$M$ with vertices at~$\mathscr P$.)

So, one can operate
with elements of~$\mcg(M,\mathscr P)$ as efficiently as with
those of~$\mathrm{GL}(n,\mathbb Z)$ provided that the geometric
intersection index can be computed as efficiently as the
standard scalar product. And this is exactly what the technical
part of the paper is devoted to---an algorithm for computing
the geometric intersection index~$\langle\gamma_1,\gamma_2\rangle$
of two multiple curves~$\gamma_1,\gamma_2$ represented
by their normal coordinates. The asymptotic running time of the algorithm, for fixed~$M$ and~$\mathscr P$,
is~$O(|\gamma_1|\cdot|\gamma_2|)$ where~$|\gamma|$ stands for a complexity measure
comparable to the amount of space needed for writing the normal coordinates of~$\gamma$
(see Proposition~\ref{compute-index-prop}).

This is used to construct an algorithm solving the word problem for~$\mcg(M,\mathscr P)$
that accepts as the input \emph{zipped} words, and whose running time is quadratic in the size of the input
provided that the generating set satisfies certain conditions (see Theorem~\ref{maintheo}). By saying `zipped' here we mean that powers of generators~$a^k$
are encoded as pairs~$(a,k)$ with~$k$ written in a positional numeral system, so
the contribution of such a term to the size of the input is~$O(\log k)$.

However,  from the practical point of view, the matrix presentation of~$\mcg(M,\mathscr P)$ has an advantage over
the zipped-word presentation in the following respects:
\begin{itemize}
\item
the asymptotic time for computing the matrix presentation~$g_1g_2$ can be made \emph{linear} in
the size of each of the matrix presentations of~$g_1$ and~$g_2$ (see Corollary~\ref{main-coro});
\item
inverting a group element requires constant
time (provided that the result is returned by reference), since it amounts just to transposing the respective matrix;
\item
the dependence of the running time of these algorithms on the complexity of the punctured surface~$(M,\mathscr P)$
is polynomial (see Theorem~\ref{n-dep-theo}).
\end{itemize}

This is important because, for computations involving a non-abelian group, one might prefer not only to have
a solution of the word problem for this group, but also an algorithm for computing a normal form of any element, which would
allow to avoid operating with unreasonably long presentations of group elements without knowing that they can
be simplified. When a normal form is defined for all group elements, the key question is how efficiently
one can compute the normal form of~$g_1g_2$ from the normal forms of~$g_1$ and~$g_2$,
and the normal form of~$g^{-1}$ from the normal form of~$g$.

The author is unaware of any approach based on decompositions of group elements into products of generators
that yields a solution of these problems with the above mentioned properties.

Many ideas we use are pretty well known to date (such as representing
curves by measured train tracks and simplifying them by a
procedure similar in nature to the accelerated Euclidean algorithm).
The feature of the method proposed here is that
to compute~$\langle\gamma_1,\gamma_2\rangle$ we simplify the presentation of both~$\gamma_1$
and~$\gamma_2$ simultaneously to the extent in which their simplifications go in
parallel, which makes all intersections detectable without fully untangling any of~$\gamma_1$
and~$\gamma_2$. In cases when only few simplification steps are needed (which
seem to be typical in a sense), this allows to benefit from fast multiplication
algorithms for integers, which are known since the work of A.\,Karatsuba~\cite{kar1,kar2}.

Another possible way to compute~$\langle\gamma_1,\gamma_2\rangle$ could be
by simplifying the presentation of~$\gamma_1$ as much as possible at the expense
of possibly getting the presentation of~$\gamma_2$ more complicated.
This idea is realized in~\cite{bell,bellwebb} where a method is suggested to change the triangulation
of the surface so that normal coordinates of~$\gamma_1$ become small. Another way
of simplifying the presentation of a normal curve on a surface is given in~\cite{en} by means
of constructing a special cell decomposition of the surface called the street complex.
The computational efficiency of these approaches, for a fixed surface,
would be comparable, in worst cases, to the one proposed here, but in good
cases, no acceleration due to Karatsuba type algorithms can be achieved.

It is worth noting that in order to have a polynomial bound for the running time of the algorithms simplifying
the presentation of a normal curve it is important to devlope an analogue of the \emph{accelerated} version
of the Euclidean algorithm (the one that uses Euclidean division instead of subtraction). This means that
the algorithm must somehow detect `large spirals' in the given normal curve and untwist each of them in a single step.
Without such a feature the running time of the algorithm would be, in worst cases, exponential in the
size of the presentation of a normal curve by normal coordinates (as in~\cite{yurttas-hall}, where
intersections of closed multiple curves in a punctured disc are counted using the `relaxation' algorithm from~\cite{cumplido}).

Saul Schleimer pointed out to the author that a similar computational efficiency (in worst cases)
to the one of the method proposed here can also be achieved
by means of straight-line programs~\cite{babai}, which provide for another way to efficiently
treat normal curves on a surface (see~\cite{slp}).

The paper is organized as follows. In Section~\ref{present-sec}
we introduce the `natural' geometry on~$\mcg(M,\mathscr P)$
in purely algebraic terms and formulate in these terms the claim about
the efficiency of our approach. Sections~\ref{conventions} and~\ref{pull-sec}
are devoted to preliminaries. The matrix presentation for~$\mcg(M,\mathscr P)$
is constructed in Sections~\ref{coord-sec} and~\ref{matrix-pres-sec}.
In Sections~\ref{trtrsec}--\ref{counting-sec} the algorithmic issues are discussed.
Section~\ref{concluding-sec} proposes a direction of further research motivated
by the matrix presentation of the mapping class groups.

\section{Group presentations and complexity}\label{present-sec}

For two non-negative functions $c_1$ and $c_2$ on a group $G$, we write
$c_1\preceq c_2$, if there exists a constant $C$ such that $c_1(g)<C\cdot c_2(g)$
for all $g\in G$, $g\ne1$. If both $c_1\preceq c_2$ and $c_2\preceq c_1$
hold we say that $c_1$ and $c_2$ are \emph{comparable}.
If only $c_1\preceq c_2$ holds but not $c_2\preceq c_1$,
we write $c_1\prec c_2$.

We say that $c:G\rightarrow\mathbb R_{\geqslant0}$ is a \emph{complexity function} if there
exists a finite alphabet $\mathscr A$ and a language $\mathscr L$
(= a subset of the set of all words) in $\mathscr A$
with an onto mapping $\ev:\mathscr L\rightarrow G$ such that
\begin{enumerate}
\item if $w_1,w_2\in\mathscr L$, then $w_1w_2\in\mathscr L$
and $\ev(w_1w_2)=\ev(w_1)\ev(w_2)$;
\item $c$ is comparable to the following function $f$:
$$f(g)=\inf_{\ev(w)=g}|w|,$$
where by $|w|$ we denote the word length.
\end{enumerate}
A couple $(\mathscr L,\ev)$ satisfying (i) will be referred to as \emph{a $G$-presentation},
and if (ii) also holds then it will be said to be \emph{appropriate} for $c$.

\begin{exam}
An ordinary word length complexity function~$\cwl_{\mathscr A}$,
where $\mathscr A\subset G$ is a finite generating set for~$G$,
is a typical example of a complexity function. An appropriate
$G$-presentation is obtained by letting $\mathscr L$ be the set of all words in $\mathscr A$.
One can see that~$\cwl_{\mathscr A}$ is always a maximal complexity function with respect to~$\preceq$.
\end{exam}

\begin{exam}\label{gln-exam}
The conventional way for encoding integral matrices, by listing their entries written
in a positional numeral system, yields a~$\mathrm{GL}(n,\mathbb Z)$-presentation
appropriate for the complexity function~$c(x)=\log\|x\|$ (any element can also
be represented as a product of other elements encoded in this way).
\end{exam}

A complexity function on~$GL(n,\mathbb Z)$ (viewed as an abstract finitely presented group)
comparable to $c$ in Example~\ref{gln-exam} can also be defined without an explicit reference
to the matrix presentation using the following general construction.

For a finite generating set $\mathscr A$ of a group $G$ we
define \emph{the zipped word length
function} $\zwl_{\mathscr A}$ as follows:
$$\zwl_{\mathscr A}(g)=\min_{\begin{array}cg=a_1^{k_1}\ldots a_m^{k_m},\\
a_1,\ldots,a_m\in\mathscr A,\\
k_1,\ldots,k_m\in\mathbb Z\end{array}}\sum_{i=1}^m\log_2(|k_i|+1).$$
Obviously, this is a complexity function, for which an appropriate $G$-presentation is
obtained by choosing a reasonable encoding for sequences
of the form $((a_1,k_1),\ldots,(a_m,k_m))$, where $a_i\in\mathscr A$, $k_i\in\mathbb Z$,
and interpreting such a sequence as the product $a_1^{k_1}\ldots a_m^{k_m}\in G$. We call it
\emph{the zipped word presentation}.

One can show that the complexity function~$c(x)=\log\|x\|$ on~$\mathrm{GL}(n,\mathbb Z)$
is comparable to $\zwl_{\mathscr A}$ if the generating set~$\mathscr A$ is chosen appropriately.
Namely, it suffices that, for each~$\ell=2,3,\ldots,n$,
the subset~$\mathscr A$ contains an element
whose Jordan normal form has a single Jordan block of size~$\ell$ and~$n-\ell$ blocks of size~$1$
with all eigenvalues equal to~$\pm1$, and there are no elements in~$\mathscr A$ having
eigenvalues other than~$\pm1$.

There is a direct analogy of this statement for the mapping class groups. The `natural'
geometry on~$\mcg(M,\mathscr P)$ can be defined in terms of
the matrix presentation introduced below in Section~\ref{matrix-pres-sec},
and this geometry coincides with the one defined by the zipped word length function
provided that the generating set is chosen appropriately (see Proposition~\ref{complexitiesequivalence}).
As shown in~\cite{shast} this geometry also coincides with the one
coming from the action of the group on the thick part of the corresponding Teichm\"uller space.

\begin{defi}\label{effsol}
For a complexity function $c$ on a group $G$, we call \emph{an efficient solution
of the word problem for $G$ with respect to $c$} an appropriate
$G$-presentation $(\mathscr L,\ev)$ together with
\begin{enumerate}
\item
a mapping $\nf:G\rightarrow\mathscr L$ (the word $\nf(g)$ is thought of
as the normal form of $g$) such that
we have $\ev\circ\nf=\mathrm{id}_G$ and the function
$g\mapsto\left|\nf(g)\right|$ is comparable to $c$,
and
\item
polynomial-time algorithms to decide wether $w\in\mathscr L$
or not and to compute $\nf(\pi(w))$
from $w$ if $w\in\mathscr L$.
\end{enumerate}
\end{defi}

\begin{defi}
Let $a,b$ be elements of a group $G$. We say that $a$ is \emph{a fractional power of} $b$
if $a^k=b^l$ for some $k,l\in\mathbb Z$, $k>0$.
\end{defi}

In particular, any torsion element is a fractional power any other group element.

\begin{theo}\label{maintheo}
Let $M$ be a compact surface, $P_1,\ldots,P_n\in M$ a non-empty collection
of pairwise distinct points such that
the mapping class group $G=\mcg(M,\{P_1,\ldots,P_n\})$ is infinite.
Let~$\mathscr A$ be a finite generating set for $G$ such that
\begin{enumerate}
\item
every element in $\mathscr A$ is a fractional power of a Dehn twist;
\item
every Dehn twist in $G$ is conjugate to a fractional power of an element from~$\mathscr A$.
\end{enumerate}
Then the word problem in $G$ is efficiently solvable with respect to $\zwl_\mathscr A$. Moreover,
the algorithms for this solution can be made quadratic-time.
\end{theo}

There are various generating sets known
satisfying Condition~(i) for the mapping class groups, see~\cite{lick1,lick2,lick3,chil,kork,brfa,szep}.
Condition~(ii) can always be met by adding a few Dehn twists to
the generating set, since up to a homeomorphism there are only finitely
many distinct simple closed curves in $M\setminus\{P_1,\ldots,P_n\}$.

Our settings here are slightly more general than those that one typically considers
(we allow multiple $P_i$'s on a single boundary component and
include orientation reversing homeotopies into the mapping class
group of an orientable surface), but extending the
existing results so as to obtain a generating set
satisfying~(i) and~(ii) is easy. So, Theorem~\ref{maintheo}
applies to any infinite mapping class group of a compact surface with $n\geqslant1$ punctures.

In the particular case when $G$ is the braid group $B_n$ and the generating set~$\mathscr A$
consists of all Garside-like elements $\Delta_{ij}$ (half-twists of strands $i$ through $j$),
Theorem~\ref{maintheo} was established by the present author and Bert Wiest in~\cite{dw}.

Theorem~\ref{maintheo} will be proved in Section~\ref{counting-sec}.

\section{Notation, terminology, and conventions}\label{conventions}

Once and for all until Subsection~\ref{n-dep-subsec} we fix a connected compact surface $M$, orientable or not, which
will be referred to simply as \emph{the surface},
and a non-empty set of \emph{punctures} $\mathscr P=
\{P_1,\ldots,P_n\}\subset M$. If~$M$ is a sphere
we require $n\geqslant4$; if $M$ is a projective plane,
a disk, an annulus, or a M\"obius
band we require $n\geqslant3$;
and if $M$ is a torus or a Klein bottle we require $n\geqslant2$.
The excluded cases will be referred to as \emph{sporadic}
and the remaining ones \emph{nonsporadic}.

By $G$ we will denote \emph{the mapping class group} $\mcg(M,\mathscr P)$, that is,
the quotient of the group $\homeo(M,\mathscr P)$ of self-homeomorphisms of $M$
preserving the subset $\mathscr P$ by the connected
component $\homeo_0(M,\mathscr P)$
containing the identity homeomorphism.

We assume that every boundary component $\gamma$
of $M$ contains at least one of~$P_i$'s. This is not
a loss of generality because otherwise one can contract
$\gamma$ to a point and treat it as a
puncture, which does not affect the mapping class group.

The punctures located at $\partial M$ will be called \emph{boundary punctures}
and all the others \emph{internal punctures}.

By \emph{a proper arc} on $M$ we mean an open simple arc
$\alpha$ in $M\setminus\mathscr P$ approaching
some punctures $P_i$, $P_j$ at the ends such that the
closure $\overline\alpha$ of $\alpha$ does not bound
an \emph{empty} disk, i.e.\ a disk with no puncture inside.
It is allowed, however, that $\overline\alpha$ forms a loop.

By \emph{a simple curve} on $M$ we mean a smooth simple
closed curve in $M\setminus\mathscr P$ that
does not bound an empty disk.

By \emph{a multiple curve} on $M$ we mean a possibly empty union of pairwise
disjoint simple curves and proper arcs on $M$.

Two proper arcs are \emph{parallel} if they coincide or enclose an empty disk.
Two simple curves are \emph{parallel} if they
enclose an empty annulus.

Two curves $\gamma_1$, $\gamma_2$ are said to
be \emph{tight} (with respect to each other)
if they either do not meet or meet transversely,
and there is no empty disk $D\subset M$
bounded by two subarcs $\alpha_1\subset\overline{\gamma_1}$
and $\alpha_2\subset\overline{\gamma_2}$
such that at least one of the common endpoints of $\alpha_1$ and $\alpha_2$ is not
a puncture.

Let $\gamma$ and $\gamma'$ be two multiple curves. We write $\gamma\sim\gamma'$ if they
are isotopic relative to $\mathscr P$.

\begin{defi}\label{intersection-def}
If two multiple curves $\gamma_1$ and $\gamma_2$ are
tight we define their \emph{geometric intersection index}
$\langle\gamma_1,\gamma_2\rangle$ to be the number of intersections
$|\gamma_1\cap\gamma_2|$ less the number of
pairs~$(\alpha_1,\alpha_2)$ of parallel proper arcs such that $\alpha_i\subset\gamma_i$, and
the number of pairs~$(\beta_1,\beta_2)$ of isotopic one-sided simple curves such that $\beta_i\subset\gamma_i$.
For arbitrary multiple curves $\gamma_1,\gamma_2$, the geometric intersection index
$\langle\gamma_1,\gamma_2\rangle$ is defined as $\langle\gamma_1',\gamma_2'\rangle$
with any tight pair $(\gamma_1',\gamma_2')$ of multiple curves such
that~$\gamma_i'\sim\gamma_i$, $i=1,2$. As we will see below
(Proposition~\ref{pulltight1}) this number
is well defined.
\end{defi}

Note that, according to this definition, a pair $(\beta_1,\beta_2)$ of isotopic simple curves
such that~$\beta_i\subset\gamma_i$ does not contribute anything to $\langle\gamma_1,\gamma_2\rangle$.
Indeed, if these curves are two-sided, then they are disjoint in tight position. If they
are one-sided, then they have a single intersection in tight position, but this contribution is cancelled
by subtracting the total number of such pairs.

Note also that, for any proper arc~$\alpha$, we have~$\langle\alpha,\alpha\rangle=-1$.
This will be justified by Proposition~\ref{normalcoord} (see also Section~\ref{concluding-sec}).

By \emph{a triangulation of $M$ with vertices at $\mathscr P$}
we mean a maximal collection of proper arcs $(e_1,\ldots,e_N)$
such that they are pairwise disjoint and nonparallel. The arcs $e_i$ are called \emph{edges}
of the triangulations. We assume additionally that the boundary $\partial M$ is covered by
$\bigcup_{i=1}^N\overline{e_i}$.

In nonsporadic cases, the edges of a triangulation cut the surface
$M$ into \emph{triangles}, which are homeomorphic
images of the interior of
a 2-simplex under a continuous map that sends the
interior of each side of the simplex into
an edge of the triangulation.

It is standard to check that the number of edges of any triangulation of $M$ with vertices
at~$\mathscr P$ is equal to
$$N=-3\chi+3n-m,$$
where $\chi$ is the Euler characteristics of $M$ and $m$ is the number of punctures
at $\partial M$, and the number of triangles equals
$$F=-2\chi+2n-m.$$

In order not to overload the exposition by technical details,
we postpone the discussion on the dependence
of the asymptotic complexity of the proposed algorithms on the complexity of the surface
till Subsection~\ref{n-dep-subsec}, and discuss it there very briefly.

\emph{`Isotopic'} in this paper always means `isotopic relative to $\mathscr P$'.

\section{The pulling tight procedure}\label{pull-sec}
Here we recall some standard facts about curves on a surface, adapted to our settings.
An experienced reader may safely skip this section.

\begin{prop}\label{pulltight1}
Let $\gamma_1$ and $\gamma_2$ be two multiple curves in $M$. Then there exist
multiple curves $\gamma_1'$ and $\gamma_2'$ such that $\gamma_i'\sim\gamma_i$, $i=1,2$,
and $\gamma_1'$ and $\gamma_2'$ are tight.
If $\gamma_1''$ and $\gamma_2''$ are another such pair of multiple curves
then $\langle\gamma_1',\gamma_2'\rangle=\langle\gamma_1'',\gamma_2''\rangle$.
\end{prop}

\begin{proof}
The standard method to produce the desired $\gamma_1'$, $\gamma_2'$ is known
as \emph{the pulling tight procedure}. We start from $\gamma_i'=\gamma_i$, $i=1,2$,
and then modify them. First, we disturb
$\gamma_1'$ and $\gamma_2'$ slightly to make them transverse to each other.

Assume there is \emph{a bigon}, i.e.\ a $2$-disk~$D\subset M$ bounded by arcs
$\overline{\alpha_1}$ and $\overline{\alpha_2}$ with
$\alpha_i\subset\gamma_i'$, $i=1,2$, being non-proper open arcs
such that the interior of $D$
is disjoint from $\gamma_1'$ and $\gamma_2'$.
\begin{figure}[ht]
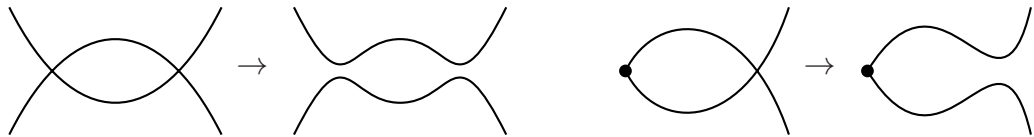

\centerline{\includegraphics[scale=.8]{tight1.eps}\hskip.2cm\raisebox{23pt}{$\rightarrow$}
\hskip.2cm\includegraphics[scale=.8]{tight2.eps}\hskip1.5cm
\includegraphics[scale=.8]{tight3.eps}\hskip.2cm\raisebox{23pt}{$\rightarrow$}
\hskip.2cm\includegraphics[scale=.8]{tight4.eps}}
\caption{Reducing bigons}\label{bigons}
\end{figure}
We replace $\gamma_1'$ and $\gamma_2'$
by $(\gamma_1'\setminus\alpha_1)\cup\alpha_2$
and $(\gamma_2'\setminus\alpha_2)\cup\alpha_1$, respectively,
and then smooth out the obtained curves at the breaking point(s) (see Figure~\ref{bigons}).
This reduces the number of intersections of $\gamma_1'$ and
$\gamma_2'$, so the process terminates after finitely many steps.
Obviously the isotopy class of each curve stays unchanged.

Figure~\ref{diffways} illustrates the fact that the order in which we reduce
bigons does not matter. More precisely, the isotopy
class of the union $\gamma_1'\cup\gamma_2'$
does not depend on that order.
\begin{figure}[ht]
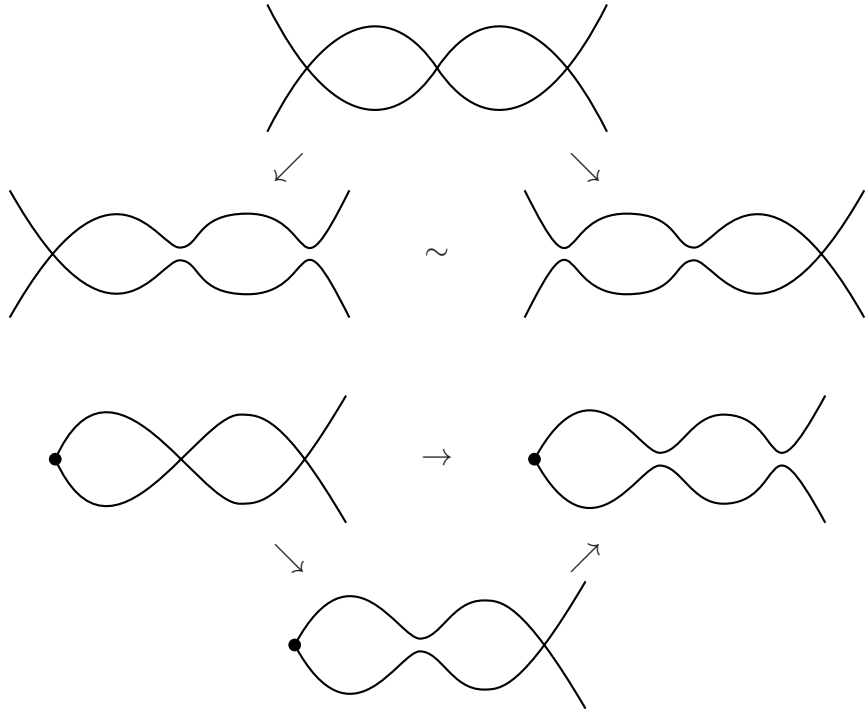

\centerline{\includegraphics[scale=.8]{tight5.eps}}
\smallskip
\centerline{$\swarrow\hskip100pt\searrow$}
\smallskip
\centerline{\includegraphics[scale=.8]{tight7.eps}\hskip1cm\raisebox{22pt}{$\sim$}\hskip1cm
\includegraphics[scale=.8]{tight6.eps}}
\vskip1cm
\centerline{\includegraphics[scale=.8]{tight8.eps}\hskip1cm\raisebox{22pt}{$\rightarrow$}\hskip1cm
\includegraphics[scale=.8]{tight9.eps}}
\smallskip
\centerline{$\searrow\hskip100pt\nearrow$}
\smallskip
\centerline{\includegraphics[scale=.8]{tight10.eps}}
\caption{Different ways of reducing bigons give isotopic pictures}\label{diffways}
\end{figure}

\begin{rem}
However, the isotopy class of the pair $(\gamma_1',\gamma_2')$ can depend
on the arbitrariness in the pulling tight process if
there are connected components $\beta_1\subset\gamma_1$, $\beta_2\subset\gamma_2$
that are isotopic to each other. The issue is illustrated
in Figure~\ref{2bigons}. Whichever bigon we reduce
we get the same \emph{unordered} pair of curves
but which of them will be $\beta_1'$ and which $\beta_2'$
depends on the choice of the bigon(s) being reduced.
\begin{figure}[ht]
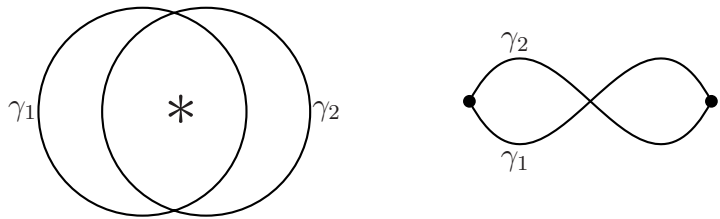

\centerline{\includegraphics[scale=.8]{2bigons1.eps}\put(-56,35){\Huge$*$}\put(-115,39){$\gamma_1$}\put(0,39){$\gamma_2$}
\hskip2cm\raisebox{20pt}{\includegraphics[scale=.8]{2bigons2.eps}\put(-82,0){$\gamma_1$}\put(-82,45){$\gamma_2$}}}
\caption{Pairs of bigons whose reduction gives non-isotopic
\emph{ordered} pairs of curves}\label{2bigons}
\end{figure}
\end{rem}

Let $\beta_1$ and $\beta_2$ be isotopic to $\gamma_1$
and $\gamma_2$, respectively, and
the multiple curves in both pairs are transverse to each other. Let $(\beta_1',\beta_2')$
and $(\gamma_1',\gamma_2')$ be obtained from
the corresponding pairs by pulling them tight.
Then we claim that $\beta_1'\cup\beta_2'$ is isotopic to $\gamma_1'\cup\gamma_2'$.

Indeed by choosing a generic isotopy from $\beta_1$ to $\gamma_1$ and from
$\beta_2$ to $\gamma_2$ we get a finite sequence of
bigon reduction and inverse operations,
that produces $\gamma_1\cup\gamma_2$ from $\beta_1\cup\beta_2$. So, it suffices
to prove the claim for a single bigon reduction. If
$\beta_1\cup\beta_2\mapsto\gamma_1\cup\gamma_2$
is a bigon reduction, then it may be taken for the
first step of pulling $\beta_1$ and $\beta_2$
tight, so the result of pulling tight procedure for $(\beta_1,\beta_2)$ and $(\gamma_1,\gamma_2)$
will be exactly the same.

Applying this to $\beta_i=\gamma_i''$, $i=1,2$,
we get $\gamma_1'\cup\gamma_2'\sim\gamma_1''\cup\gamma_2''$. This implies the second
claim of the proposition.
\end{proof}

The first claim in Proposition~\ref{pulltight1} can be strengthen as follows.

\begin{prop}\label{pulltight2}
Let $\gamma_1$ and $\gamma_2$ be two multiple curves in $M$. Then there exists
a multiple curve~$\gamma_2'$ such that we have $\gamma_2'\sim\gamma_2$
and $\gamma_1,\gamma_2'$ are tight.
\end{prop}

\begin{proof}
One only needs to apply an isotopy that carries $\gamma_1'$ to $\gamma_1$ at the end
of the pulling tight procedure described in the proof of Proposition~\ref{pulltight1}.
One can restore $\gamma_1$ by an isotopy not only
at the very end but also at every step of the procedure.
\end{proof}

\begin{prop}\label{pulltight3}
For any three multiple curves $\gamma_1$, $\gamma_2$, and $\gamma_3$, there
are multiple  curves $\gamma_2'$ , $\gamma_3'$ such
that $\gamma_i'\sim\gamma_i$, $i=2,3$,
and the curves $\gamma_1$, $\gamma_2'$, and $\gamma_3'$
are pairwise tight.
\end{prop}

\begin{proof}
Due to Proposition~\ref{pulltight2} we may assume without loss of generality that
the pairs $(\gamma_1,\gamma_2)$ and $(\gamma_1,\gamma_3)$ are already tight.
We may also assume that there are no triple intersections, i.e.
$\gamma_1\cap\gamma_2\cap\gamma_3=\varnothing$ as we can
achieve this by a small deformation of $\gamma_3$.

Now we apply the pulling tight procedure to $(\gamma_2,\gamma_3)$.
It produces $\gamma_2'$, $\gamma_3'$ that are still tight with respect
to $\gamma_1$. Indeed, $\gamma_2'\cup\gamma_3'$
is obtained from $\gamma_2\cup\gamma_3$ by resolution of
intersections, which occur far from $\gamma_1$. So,
the number of intersection points in $\gamma_1\cap(\gamma_2'\cup\gamma_3')$
is equal to that in $\gamma_1\cap(\gamma_2\cup\gamma_3)$.
If $(\gamma_1,\gamma_2')$ or $(\gamma_1,\gamma_3')$ were not tight,
we could have applied the pulling tight process again and get
$\gamma_2''$ and $\gamma_3''$ such that $\gamma_2''\cup\gamma_3''$
has a smaller number of intersections with~$\gamma_1$ than $\gamma_2\cup\gamma_3$ has,
which contradicts Proposition~\ref{pulltight1}.
\end{proof}

\section{Normal coordinates}\label{coord-sec}
The idea of a normal curve and normal coordinates goes back to~H.\,Kneser~\cite{kneser}
who introduced the concept of a normal surface in a $3$-manifold,
which has had a big impact on low-dimensional topology.
We use a modification the classical notion of a normal curve which allows arcs emanating from
the punctures.

\begin{figure}[ht]
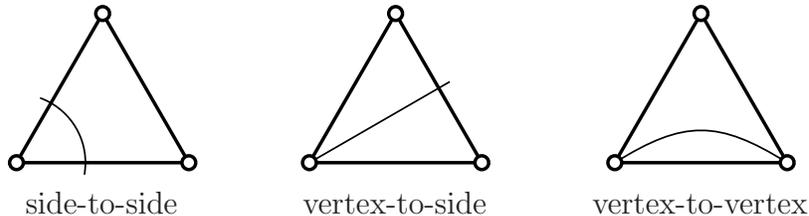

\begin{tabular}{ccc}
\includegraphics[scale=.65]{narc1.eps}&\hbox to 1cm{}%
\includegraphics[scale=.65]{narc2.eps}\hbox to 1cm{}&
\includegraphics[scale=.65]{narc3.eps}\\
side-to-side&vertex-to-side&vertex-to-vertex
\end{tabular}
\caption{Normal arcs}\label{normalarcs}
\end{figure}

Let $T=(e_1,\ldots,e_N)$ be a triangulation of $M$.
A multiple curve $\gamma$ is said to be \emph{a normal curve with respect to $T$} if
$\gamma$ and $\bigcup_{i=1}^Ne_i$ are tight. This is equivalent to
saying that the intersection of $\gamma$ with any triangle $\tau$ of
$T$ consists of arcs each of which either connects points
on different sides of $\tau$ or a vertex to a point on the
opposite side or two vertices (in the latter case such an arc is
parallel to a side of $\tau$); see Figure~\ref{normalarcs}. Such arcs
will be called \emph{normal}. By saying that an arc is normal we will
also assume that its endpoints do not lie in the interior of boundary
edges of $M$ as this never happens to intersections of multiple curves
with triangles. We will call a normal arc
\emph{side-to-side}, \emph{vertex-to-side}, or \emph{vertex-to-vertex}
according the location of its endpoints.

If $\gamma$ is not normal with respect to $T$, then the pulling tight
process for $\gamma$ and $\cup_ie_i$ with the latter staying fixed (see
Proposition~\ref{pulltight2}) will be referred to as \emph{normalization} of $\gamma$
with respect to~$T$.

\begin{defi}
For a multiple curve~$\gamma$, the numbers~$\langle\gamma,e_i\rangle$, $i=1,\ldots,N$,
are called \emph{normal coordinates of~$\gamma$ with respect to~$T$}.
\end{defi}

\begin{prop}\label{normalcoord}
Let $\gamma_1$ and $\gamma_2$ be multiple curves such
that $\langle\gamma_1,e_i\rangle=\langle
\gamma_2,e_i\rangle$ for all $i=1,\ldots, N$. Then $\gamma_1$ and $\gamma_2$ are
isotopic.
\end{prop}

\begin{proof}
Due to Proposition~\ref{pulltight2} we may restrict ourselves to the case when $\gamma_1$
and $\gamma_2$ are normal with respect to $T$.

Let $\gamma_1'$ and $\gamma_2'$ be multiple curves
obtained from $\gamma_1$ and $\gamma_2$,
respectively, by removing all proper arcs parallel to edges of $T$. Then we still have
$\langle\gamma_1',e_i\rangle=\langle\gamma_2',e_i\rangle$
for all $i=1,\ldots,N$, and, moreover,
all these geometric intersection indexes are non-negative.

Since the number of intersections of $\gamma_1'$ with $e_i$ coincides with that
of $\gamma_2'$ for any $i=1,\ldots,N$, we can apply an isotopy that preserves
all $e_i$s and carries $\gamma_1'\cap e_i$ to $\gamma_2'\cap e_i$ for all $i$.
So, we may assume that $\gamma_1'\cap e_i=\gamma_2'\cap e_i$ for $i=1,\ldots,N$.

Now we focus on a single triangle $\tau$ of $T$.

\begin{lemm}\label{lm1}
The intersection of
a normal curve $\gamma$ with $\tau$ can be recovered from $\gamma\cap\partial\tau$ uniquely
up to isotopy relative to $\partial\tau$ provided that $\gamma$ has no components parallel to
the sides of $\tau$.
\end{lemm}

\begin{proof}
Indeed, let $\tau$ be bounded by the edges $e_1$,
$e_2$, $e_3$. Denote $x_i=\langle\gamma,e_i\rangle$.
Each normal arc in $\gamma\cap\tau$ connects either a point at $e_i$, $i\in\{1,2,3\}$, to the
opposite vertex, in which
case we say that it has \emph{type} $(0i)$, or two points on different sides $e_i$ and $e_j$,
$i,j\in\{1,2,3\}$, $i<j$, in which case we attribute it \emph{type} $(ij)$; see Figure~\ref{arctypes}.
\begin{figure}[ht]
\includegraphics[scale=.75]{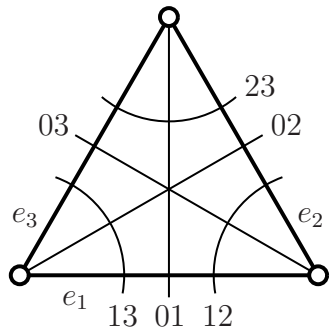}
\put(-68,-10){$01$}
\put(-50,-10){$12$}
\put(-86,-10){$13$}
\put(-24,62){$02$}
\put(-112,62){$03$}
\put(-34,77){$23$}
\put(-103,-2){$e_1$}
\put(-122,29){$e_3$}
\put(-14,29){$e_2$}
\caption{Six types of normal arcs in a single triangle}\label{arctypes}
\end{figure}

One can see that some of these
types are incompatible meaning that normal arcs of those
types can not occur in $\gamma\cap\tau$
simultaneously.
For example, type (01) is incompatible with types (02),
(03), and (23). Therefore, if a normal arc of
type (01) is present in $\gamma\cap\tau$, then the other
arcs may have only types (12) and (13).
One can readily see that, in this case, we have $x_1>x_2+x_3$,
and the number of normal arcs
of type (01), (12), and (13) is equal to $(x_1-x_2-x_3)$, $x_2$, and $x_3$, respectively.

Similarly, if arcs of type (02) or (03) are present, we will have $x_2>x_1+x_3$ or
$x_3>x_1+x_2$, respectively, and, in each case, recover the number of arcs of
each type from $x_1,x_2,x_3$.

If only arcs of types (12), (23), (13) are present, then the triangle inequalities
hold for $x_1$, $x_2$, and $x_3$, and the number of arcs of type (12), (23), and (13)
is equal to $(x_1+x_2-x_3)/2$, $(x_2+x_3-x_1)/2$, $(x_1+x_3-x_2)/2$, respectively.

The sets of all possible triples $(x_1,x_2,x_3)$ obtained in theses four cases do not overlap;
hence, from the knowledge of $x_1$, $x_2$, $x_3$ we can always decide which case occurs.

Clearly, the number of normal arcs of each type defines $\gamma\cap\tau$
up to isotopy relative to $\partial\tau$.
\end{proof}

We resume the proof of Proposition~\ref{normalcoord}.
It follows from Lemma~\ref{lm1} that $\gamma_1'$ and~$\gamma_2'$
are isotopic relative to $\bigcup_{i=1}^Ne_i$.
Thus, we may assume $\gamma_1'=\gamma_2'$ from the beginning.

Now let $k_i$ be $\max(0,-\langle\gamma_1,e_i\rangle)=\max(0,-\langle\gamma_2,e_i\rangle)$.
Each of $\gamma_1$ and $\gamma_2$ is obtained from $\gamma_1'=\gamma_2'$ by adding
$k_i$ proper arcs parallel to $e_i$ for all $i=1,\ldots,N$. Obviously, the result is unique
up to isotopy.
\end{proof}

\begin{rem}\label{setL-rem}
One can see from the proof of Proposition~\ref{normalcoord} that the collections of normal coordinates
of normal curves form a subset~$L$ of~$\mathbb Z^N$ that can be characterized as follows:
$(x_1,x_2,\ldots,x_N)\in\mathbb Z^N\setminus L$ if and only if either there is a boundary
edge~$e_i$ such that~$x_i>0$ or there is a triangle of~$T$
with sides~$e_i$, $e_j$, and~$e_k$ such that the respective coordinates~$x_i$, $x_j$, and~$x_k$
are all positive, satisfy the triangle inequalities, and sum up to an odd number.
\end{rem}

Let $T'=(e_1',\ldots,e_N')$ be another triangulation of $M$ with vertices at $\mathscr P$.
We denote by $\langle T,T'\rangle$ the $N\times N$ matrix whose $(ij)$th entry is
$$\langle T,T'\rangle_{ij}=\langle e_i,e_j'\rangle.$$

\begin{rem}
The determinant of the matrix~$\langle T,T'\rangle$ is always an integral power of two. This fact
plays no role here, but the reader might enjoy trying to prove this.
\end{rem}

Recall from Section~\ref{conventions} that we deal only with nonsporadic cases.

\begin{prop}\label{(i)<=>(ii)}
Let $\varphi$ be a diffeomorphism of $M$. Then the following two
statements are equivalent:
\begin{enumerate}
\item
$\langle T,T'\rangle=\langle T,\varphi(T')\rangle$, where
$\varphi(T')=(\varphi(e'_1),\ldots,\varphi(e'_N))$;
\item
 $\varphi$ is isotopic to the identity.
\end{enumerate}
\end{prop}

\begin{proof}
Implication (ii)$\Rightarrow$(i) follows from Proposition~\ref{pulltight1}.

Suppose now that (i) holds. Then it follows from
Proposition~\ref{normalcoord} that $\varphi(T')$ is isotopic to~$T'$.
So, without loss of generality we may assume $\varphi(T')=T'$.

Let us choose an orientation for each triangle of $T'$. Clearly, if $\varphi$ carries
each triangle to itself and preserves its orientation, then $\varphi$ is isotopic to identity.
Suppose, to the contrary, that this is not the case.

Then either there are two different triangles $\tau$, $\tau'$ of $T'$ such that
$\varphi(\tau)=\tau'$, or $\varphi$ preserves each triangle but flips the orientation.

In the former case, the triangles $\tau$ and $\tau'$ have the same sides,
hence $\overline{\tau\cup\tau'}$ is a closed surface, which is the whole of~$M$.
By gluing up two triangles along all three sides one can obtain only the following
surfaces: a sphere $\mathbb S^2$ with $3$ punctures, a torus $\mathbb T^2$
with a single puncture, a projective plain $\mathbb RP^2$ with two punctures, and
a Klein bottle $\mathbb K^2$ with a single puncture. All these are sporadic cases.

In the latter case, two sides of every triangle $\tau$ of $T'$ are glued together,
hence $\overline\tau$ is either a disk $\mathbb D^2$ with a single puncture inside
and a single puncture at the boundary, or a M\"obius band $\mathbb M^2$
with a single puncture at the boundary. Both are sporadic cases, so there must
be more than one triangle of $T'$.
Two such surfaces glued along the boundary form a closed surface, so the number
of triangles cannot be greater than two. From two triangles we get
either $(M,n)=(\mathbb S^2,3)$ or $(M,n)=(\mathbb RP^2,2)$
or $(M,n)=(\mathbb K^2,1)$, which are also sporadic cases.

So, in all nonsporadic cases we have (i)$\Rightarrow$(ii).
\end{proof}

\section{A matrix presentation of the mapping class groups}\label{matrix-pres-sec}

Whenever $\gamma$ is a multiple curve and $g$ is an element of $G=\mcg(M,\mathscr P)$
we will use the notation $g(\gamma)$ for $\varphi(\gamma)$, where
$\varphi$ is any diffeomorphism representing~$g$. We will do it when
only the isotopy class of $\varphi(\gamma)$ matters.
This will apply also to triangulations in place of curves.

Let us fix a triangulation $T=(e_1,\ldots,e_N)$ of $M$ with
vertices at $\mathscr P$. It follows from Proposition~\ref{(i)<=>(ii)} that
an element $g\in G$ can be recovered uniquely from the matrix
$\langle T,g(T)\rangle$. Thus, by choosing a proper encoding for $N\times N$-matrices
we get a $G$-presentation in which an element $g\in G$ can be presented by
any sequence of matrices $(m_1,\ldots,m_k)$ such that
$m_i=\langle T,g_i(T)\rangle$ with $g_1,\ldots,g_k\in G$,
$g_1\cdot\ldots\cdot g_k=g$.

For any~$g\in G$, we let the intersection matrix~$\langle T,g(T)\rangle$ (encoded in a reasonable way using a finite alphabet)
be the normal from~$\nf(g)$ of~$g$ and define the complexity of~$g$ as
\begin{equation}\label{ct(g)-eq}
c_T(g)=\sum_{i,j=1}^N\log_2(|\langle T,g(T)\rangle_{ij}+\delta_{ij}|+1),
\end{equation}
where $\delta_{ij}$ is the Kroneker delta (which is added just to set the
complexity of the identity element to zero and plays no role otherwise). One can see that~$c_T(g)$ is comparable
to the amount of space needed to encode~$\langle T,g(T)\rangle$. The rest of this section
is devoted to showing that this setup satisfies Condition~(i) of Definition~\ref{effsol}.

The key question about the efficiency of this approach is
how to compute $\langle T,g_1(g_2(T))\rangle$
from $\langle T,g_1(T)\rangle$ and $\langle T,g_2(T)\rangle$ for arbitrary $g_1,g_2\in G$.
We start by observing that this computation has much in common with
the ordinary matrix multiplication.

\begin{prop}\label{product}
The matrix element $\langle T,g_1(g_2(T))\rangle_{ij}$ equals $\langle\gamma,\gamma'\rangle$,
where $\gamma$ and $\gamma'$ are the normal curves whose normal coordinates
with respect to $T$ form the $i$th row of $\langle T,g_1(T)\rangle$
and the $j$th column of $\langle T,g_2(T)\rangle$, respectively.
\end{prop}

\begin{proof}
Let $\gamma=g_1^{-1}(e_i)$ and $\gamma'=g_2(e_j)$. Then we have
$$\langle T,g_1(g_2(T))\rangle_{ij}=
\langle g_1^{-1}(T),g_2(T)\rangle_{ij}=
\langle\gamma,\gamma'\rangle.$$
The $k$th coordinate of $\gamma$ is
$$\langle\gamma,e_k\rangle=\langle g_1^{-1}(e_i),e_k\rangle=
\langle e_i,g_1(e_k)\rangle=\langle T,g_1(T)\rangle_{ik}.$$
The $k$th coordinate of $\gamma'$ is
\begin{align*}\langle\gamma',e_k\rangle&=
\langle e_k,\gamma'\rangle=\langle e_k,g_2(e_j)\rangle=
\langle T,g_2(T)\rangle_{kj}.\end{align*}
The claim follows.
\end{proof}

For $1\leqslant i,j\leqslant N$ let $\mu_{ij}$ be equal to the number of triangles of $T$
adjacent to both $e_i$ and $e_j$ if $i\ne j$,
and $1$ otherwise.

\begin{prop}\label{intersectionbound}
For any two curves $\gamma$, $\gamma'$ we have
$$|\langle\gamma,\gamma'\rangle|\leqslant\sum_{i,j=1}^N
|\langle\gamma,e_i\rangle|\cdot\mu_{ij}\cdot|\langle\gamma',e_j\rangle|.$$
\end{prop}

\begin{proof}
Due to Proposition~\ref{pulltight3} we may assume that
$\gamma$ and $\gamma'$ are tight and each of them is normal
with respect to $T$. We may also assume that $\gamma\cap\gamma'$
is disjoint from the edges of $T$.

Denote by $X_i$
the intersection set $\gamma\cap e_i$, if
$\langle\gamma,e_i\rangle\geqslant0$, and
the set of proper arcs in $\gamma$ parallel to $e_i$ otherwise.
In both cases we have $|X_i|=|\langle\gamma,e_i\rangle|$.
We define $X_i'$ similarly, with $\gamma'$ in place of $\gamma$.

Denote by $Y$ the set of transverse intersections in $\gamma\cap\gamma'$
joined with the set of all pairs $(\alpha,\alpha')$ of parallel proper arcs
with $\alpha\subset\gamma$, $\alpha'\subset\gamma'$.
We clearly have $|\langle\gamma,\gamma'\rangle|\leqslant|Y|$.

Now define maps $f,f'$ from $Y$ to $\bigl(\cup_iX_i\bigr)$
and $\bigl(\cup_iX_i'\bigr)$, respectively, as follows.
Let $P\in Y$ be an intersection point of $\gamma$ and $\gamma'$. Let $\tau$
be the triangle of $T$ in which this intersection occurs,
and $\alpha\subset\gamma$, $\alpha'\subset\gamma'$
be the normal arcs that contain $P$. If $\alpha$ is a proper arc
parallel to an edge $e_i$ we put $f(P)=\alpha$.
Otherwise, $\alpha$ must have an endpoint $Q$ at
some edge of $T$. In this case we put $f(P)=Q$.
We define the map $f'$ similarly, by replacing $\alpha$
with $\alpha'$.

Now let $P=(\alpha,\alpha')\in Y$ be a pair of parallel proper arcs.
If they are parallel to some $e_i$ we put
$f(P)=\alpha$ and $f'(P)=\alpha'$. Otherwise,
$\alpha$ and $\alpha'$ must intersect some edge $e_i$.
Then we choose $Q\in\alpha\cap e_i$ and $Q'\in\alpha'\cap e_i$ and put $f(P)=Q$, $f'(P)=Q'$.

It is now easy to check that due to normality and tightness of $\gamma$, $\gamma'$
the number of preimages of any $(Q,Q')\in\bigl(\cup_iX_i\bigr)\times\bigl(\cup_iX_i'\bigr)$
under the map $f\times f'$ does not exceed $\mu_{ij}$ if $Q\in X_i$ and~$Q'\in X_j'$.
Therefore, we have
\begin{align*}|\langle\gamma,\gamma'\rangle|\leqslant|Y|\leqslant\sum_{i,j=1}^N|X_i|\cdot\mu_{ij}\cdot|X_j'|&=
\sum_{i,j=1}^N
|\langle\gamma,e_i\rangle|\cdot\mu_{ij}\cdot|\langle\gamma',e_j\rangle|.\qedhere\end{align*}
\end{proof}

\begin{prop}\label{cccc}
There exists a constant $C$ depending on $M$ and $\mathscr P$ such that
\begin{equation}\label{normnorm}
c_T(g_1g_2\ldots g_k)\leqslant C(c_T(g_1)+c_T(g_2)+\ldots+c_T(g_k)).\end{equation}
for any $k\in\mathbb N$, $g_1,g_2,\ldots, g_k\in G$.
\end{prop}

\begin{proof}
For a matrix $A$, we denote by $\|A\|_{\mathrm E}$ the standard Euclidean norm of $A$:
\begin{equation}\label{eu-norm-eq}\|A\|_{\mathrm E}=\sqrt{\sum_{i,j}A_{ij}^2}.
\end{equation}

Since the numbers of summands in~\eqref{ct(g)-eq} and~\eqref{eu-norm-eq}
are fixed, the function $E:G\rightarrow\mathbb R$ defined by
$$E(g)=\log_2\|\langle T,g(T)\rangle\|_{\mathrm E}$$
is comparable to $c_T$. Therefore, it suffices to prove~\eqref{normnorm}
for $E$ in place of $c_T$. This is done by using Propositions~\ref{product}
and~\ref{intersectionbound}, which imply
$$E(g_1g_2)\leqslant E(g_1)+E(g_2)+\mu,$$
where $\mu=\log_2\|(\mu_{ij})\|_{\mathrm E}$. The rest of the proof is easy.
\end{proof}

Thus, we are done with showing that the matrix presentation introduced in
this section satisfies Condition~(i) of Definition~\ref{effsol}.
The key question now is how to compute $\langle\gamma,\gamma'\rangle$ efficiently for
two normal curves given by their normal coordinates.

\section{Train tracks}\label{trtrsec}
Train tracks, which have been introduced by W.\,Thurston, are widely used
for studying homeomorphisms of surfaces and related problems~\cite{bh,split,mosher-expansions,ph,thurston}.
Here by a train track track we mean what is known as a train track with terminals~\cite{mosher-expansions}.

Whenever we deal with a finite graph (i.e.~a $1$-dimensional CW-complex)
$\theta$ embedded in~$M$ we assume that all edges
of $\theta$ are smooth images of a closed interval, and that
an open contractible neighborhood $U_v$ is chosen around each vertex $v$ of $\theta$
so that $U_v$ and $U_{v'}$ do not overlap for any two different
vertices $v$ and $v'$, and the intersection $U_v\cap\theta$
is contractible for all $v$. The closure $t$ of any
connected component of $(U_v\cap\theta)\setminus\{v\}$ is
called \emph{a tail} of the edge whose closure contains~$t$.

Loops and multiple edges with the same endpoints are allowed for graphs.

The edges of graphs that we consider are not allowed to pass through a puncture, but
a puncture may be a vertex of a graph.

A connected component $A$ of $U_v\setminus\theta$ is called \emph{a cusp} in two cases:
\begin{enumerate}
\item[(a)]
the boundary $\partial A$ contains two tails whose tangent rays at $v$ coincide;
\item[(b)]
the vertex $v$ coincides with a puncture, i.e.\ $v\in\mathscr P$.
\end{enumerate}
In the latter case the cusp is called \emph{special}, and otherwise \emph{ordinary}.

\begin{defi}\label{ttdef}
By a \emph{train track} we mean an embedded 1-dimensional CW-complex $\theta\subset M$
consisting of two disjoint parts $\theta_1,\theta_2$ such that
\begin{enumerate}
\item
$\theta_1$ is a union of pairwise disjoint smooth simple closed curves disjoint from $\mathscr P$;
\item
$\theta_2$ is a graph whose edges have interiors disjoint from~$\mathscr P$;
\item
every vertex $v$ of $\theta$ such that $v\notin\mathscr P$ is
a 3-valent \emph{switch}, which means the following.
There are exactly three tails attached to $v$, and they can
be numbered $t_1,t_2,t_3$ so that
$t_1$ forms a smooth arc together with any of $t_2$ and
$t_3$ (thus, $t_2$ and $t_3$ give rise to a cusp).
The tail $t_1$ will be referred to as \emph{outgoing}, and $t_2$, $t_3$ \emph{ingoing};
\item\label{c4}
no connected component of $M\setminus\theta$ is an empty
disk with exactly two ordinary cusps and no special cusp;
\item\label{c5}
no connected component of $M\setminus\theta$ is an empty disk with less than two cusps.
\end{enumerate}
Connected components of $\theta_1$ and edges of
$\theta_2$ will be referred to as \emph{branches} of $\theta$.
Branches that are not attached to at least one switch are called \emph{free}. In particular,
all branches contained in~$\theta_1$ are such.

If both tails of an edge of $\theta_2$ are outgoing or one is outgoing and the other
is attached to a puncture, then the edge is called
\emph{a wide branch} of $\theta$.
\end{defi}

So, our train tracks may have vertices at punctures, and those vertices are \emph{not}
switches. For instance, the closure $\overline\gamma$
of any multiple curve $\gamma$ is a train track in our sense.

Let $(\theta,w)$ be a pair in which $\theta$ is a train track, and $w$ is an assignment to every
branch a non-negative integer, which is referred to as \emph{the width} of
the branch, such that, for every switch, the sum of the widths
of the ingoing tails equals to the width of the outgoing
one. We will call such a pair \emph{a measured train track}.

\emph{The complexity} $|(\theta,w)|$ of a measured train track $(\theta,w)$ is defined as
$$|(\theta,w)|=\sum_{\alpha}\bigl(1+\log_2(w(\alpha)+1)\bigr),$$
where the sum is taken over all branches of $\theta$. One can see
from this formula that the complexity decreases whenever the width
of a branch decreases or a branch with zero width is removed.

\begin{figure}[ht]
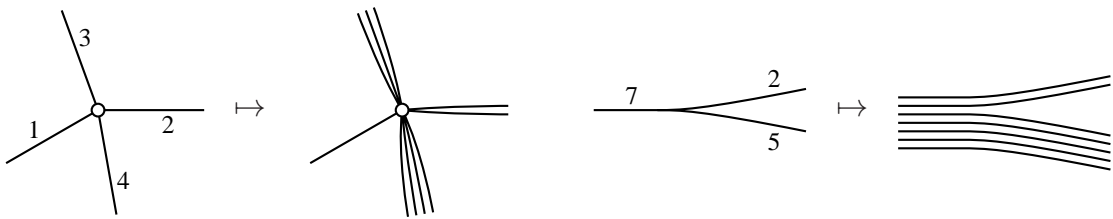

\includegraphics[scale=.8]{thicktrack1.eps}\hskip.4cm\raisebox{38pt}{$\mapsto$}\hskip.4cm\includegraphics[scale=.8]{thicktrack2.eps}
\hskip1cm
\includegraphics[scale=.8]{thicktrack3.eps}\hskip.4cm\raisebox{38pt}{$\mapsto$}\hskip.4cm\includegraphics[scale=.8]{thicktrack4.eps}
\caption{Turning a measured train track into a multiple curve}\label{track2curve}
\end{figure}

Every measured train track $(\theta,w)$ encodes a multiple curve as follows.
Each branch $\alpha$ of~$\theta$ is replaces by as many as~$w(\alpha)$
`parallel' copies $\alpha_1,\ldots,\alpha_{w(\alpha)}$ of $\alpha$.
If $\alpha$ is attached to a puncture $P_i$, then the corresponding
arcs~$\alpha_j$ approach $P_i$ at the corresponding end. At every switch,
the parallel copies of ingoing tails are attached to
that of the outgoing one so as to get a non-selfintersecting curve; see Figure~\ref{track2curve}.
(If $\alpha$ is a branch of $\theta$ having the form of a one-sided closed simple curve, then `$w(\alpha)$ parallel copies of $\alpha$'
should be understood `locally'. Precisely this means `$[w(\alpha)/2]$ parallel copies of the boundary of a small tubular
neighborhood of $\alpha$ and, if $w(\alpha)$ is odd, $\alpha$ itself'. Here $[x]$ stands for the integral part of $x$.)

A curve $\gamma$ obtained in this way from $\theta$ for some
choice of branch widths is said to be \emph{carried} by the train track~$\theta$.

\begin{figure}[ht]
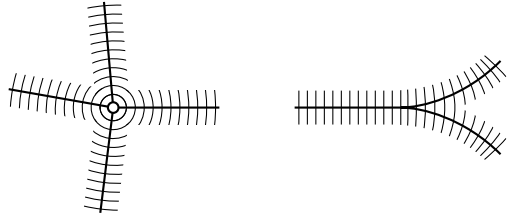

\centerline{\includegraphics[scale=.8]{foliation2.eps}\hskip1cm\includegraphics[scale=.8]{foliation1.eps}}
\caption{Foliation $\mathscr F_\theta$ near punctures and switches}\label{foli}
\end{figure}

More formally, the correspondence between measured train tracks
and curves can be described as follows.
For every train track $\theta$, we fix a singular foliation~$\mathscr F_\theta$ on $M$ such that:
\begin{enumerate}
\item
every branch of $\theta$ is transverse to $\mathscr F_\theta$
everywhere except at the punctures;
\item
$\mathscr F_\theta$ has only isolated singularities;
\item
$\mathscr F_\theta$ has a center-like singularity at every puncture. All other
singularities are outside of $\theta$ (see Figure~\ref{foli});
\item
every connected component of $M\setminus\theta$ contains a singularity of $\mathscr F_\theta$.
\end{enumerate}
In order to construct such a foliation one first
defines it in a small neighborhood of $\theta$ so as
to enforce (i) and (iii), then in a small disk in every connected component of $M\setminus\theta$
so as to enforce (iv), and then continue to the whole surface generically.

\begin{defi}
A union $\gamma$ of pairwise disjoint proper arcs and simple curves
is said to be \emph{carried by a train track} $\theta$ if
$\gamma$ is transverse to $\mathscr F_\theta$ and
there exists a homotopy $f:\gamma\times[0,1]\rightarrow M$ such that
\begin{enumerate}
\item
there are no singularities of $\mathscr F_\theta$ in $f(\gamma\times[0,1])$;
\item
for all $x\in\gamma$ we have
$f(x,0)=x$, $f(x,1)\in\theta$;
\item
all the leaves of the foliation on $\gamma\times[0,1]$ induced by $f$ from $\mathscr F_\theta$
have the form $x\times[0,1]$, $x\in\gamma$.
\end{enumerate}

The map $\pi_\gamma:\gamma\rightarrow\theta$ defined by $\pi_\gamma(x)=f(x,1)$ is called
\emph{the projection} of $\gamma$ to $\theta$. Due to Condition~(i) in this definition and
Condition~(iv) in the definition of $\mathscr F_\theta$ one can see that $\pi_\gamma$ does
not depend on a particular choice of the homotopy $f$.

The measured train track $(\theta,w)$ that \emph{encodes} $\gamma$ is defined by
letting $w(\alpha)$, where $\alpha$ is a branch of~$\theta$,
be the number of points in $\pi_\gamma^{-1}(y)$
with $y$ a point from the interior of $\alpha$.
\end{defi}

\begin{figure}[ht]
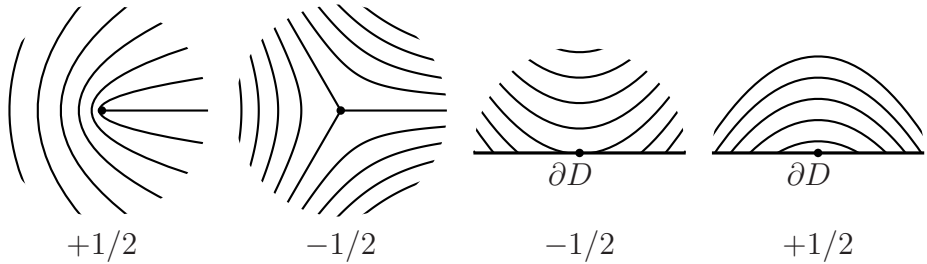

\begin{tabular}{cccc}
\includegraphics[scale=.8]{sing1.eps}&\includegraphics[scale=.8]{sing2.eps}&
\includegraphics[scale=.8]{sing3.eps}\put(-52,12){$\partial D$}&
\includegraphics[scale=.8]{sing4.eps}\put(-52,12){$\partial D$}\\
$+1/2$&$-1/2$&$-1/2$&$+1/2$
\end{tabular}
\caption{Generic singularities of a foliation and their topological indexes}\label{folsing}
\end{figure}

We will use the well known relation between the Euler characteristics of a compact
surface~$D$ and singularities of a generic foliation $\mathscr F$ on $D$.
Namely, the Euler characteristics $\chi(D)$ is equal to the sum of topological indexes
of all singularities of $\mathscr F$ provided that the points of $\partial D$ in which the
leaves of $\mathscr F$ are not transverse to $\partial D$
are also regarded as singularities. The simplest singularities
and their topological indexes are shown in Figure~\ref{folsing}.
Thus, Conditions~\ref{c4} and \ref{c5} in Definition~\ref{ttdef} simply
mean that the sum of indexes of all singularities of $\mathscr F_\theta$ inside
any connected component $D$ of $M\setminus\theta$ is non-positive unless
there is a puncture inside $D$,
and the sum is strictly negative unless there is a puncture
inside $D$ or at the boundary $\partial D$.

This implies, in particular, the following.

\begin{prop}
Any curve $\gamma$ encoded by a measured train track satisfies the conventions
that we introduced in Section~\ref{conventions}.
Namely, if $D\subset M$ is a disk that is bounded by
the closure of a connected component of $\gamma$,
then $D$ is not empty \emph(i.e.\ contains a puncture\emph).
\end{prop}

\begin{proof}
Indeed, if the connected component in question is a closed
curve, then it is transverse to $\mathscr F_\theta$. If it is an
arc whose closure forms a loop, then (after an appropriate smoothing)
it will contribute just $1/2$ to the sum of the singularity indexes
whereas we have $\chi(D)=1$. So, in both
cases the total contribution of singularities from
the interior of $D$ must be positive.

Since the boundary $\partial D$ can be homotoped to its
projection $\pi_\gamma(\partial D)\subset\theta$ through a family of curves that remain
transverse to $\mathscr F_\theta$ (except at one point
in the case when $\partial D$ is the closure of a proper arc)
there is a family of connected components $D_1,\ldots,D_k$
of $M\setminus\theta$ such that the set of singularities inside $D$ coincides
with that inside $D_1\cup\ldots\cup D_k$. Since the sum
of the topological indexes of singularities inside $D$
is positive, some $D_i$ contains a puncture, and so does $D$.\end{proof}

\section{Universal train tracks}\label{univtr-sec}
With every triangulation $T=(e_1,\ldots,e_N)$ we associate a train track $\theta_T$
having the following property:
for any multiple curve $\gamma$, there is another multiple curve
$\gamma'$ isotopic to $\gamma$ such that~$\theta_T$ carries $\gamma'$.
For this reason we call this train track \emph{universal}.
It is not uniquely defined but the arbitrariness
in its definition will not matter.

We construct $\theta_T$ in three steps.

\begin{figure}[ht]
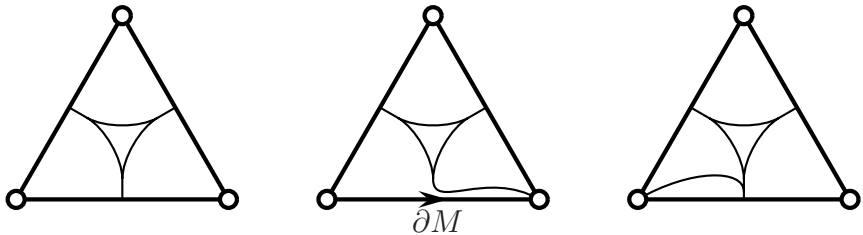

\centerline{\includegraphics[scale=.8]{univ1.eps}\hskip1cm
\includegraphics[scale=.8]{univ2.eps}\put(-52,-6){$\partial M$}\hskip1cm
\includegraphics[scale=.8]{univ3.eps}}
\caption{Constructing the train track $\theta_T$}\label{univtrack}
\end{figure}

\smallskip\noindent\emph{Step 1.} Put three switches in each triangle of $T$
and mark a single point in each edge of $T$. Outgoing tails are connected by arcs
to the marked points, and ingoing ones are paired so as to make
three-cusped disk in each triangle (see Figure~\ref{univtrack} on the left).

\smallskip\noindent\emph{Step 2.}
Orient connected components of $\partial M$ arbitrarily. Then we detach
the edges of the graph under construction
from the marked points at $\partial M$ and pull them in
the direction defined by the orientation of the corresponding
edge of $T$ toward the nearest puncture at $\partial M$ (see Figure~\ref{univtrack} in the center).

\smallskip\noindent\emph{Step 3.}
Let $\theta$ be the graph constructed so far. For every internal puncture $P$, the
connected component of $M\setminus\theta$ containing $P$ is a disk with
smooth boundary, and this disk contains no other puncture. We put an additional
switch at its boundary and connect it by a new branch with $P$.
We put
the new switches at wide branches of $\theta$ and position the new branches
as shown in Figure~\ref{univtrack} on the right.
Namely, each new branch must be contained entirely in a single triangle, and
its smooth extension through the new switch should point to the nearest
edge of the triangle, i.e.\ away of the 3-cusped disk located inside the triangle.

\begin{figure}[ht]
\centerline{\includegraphics{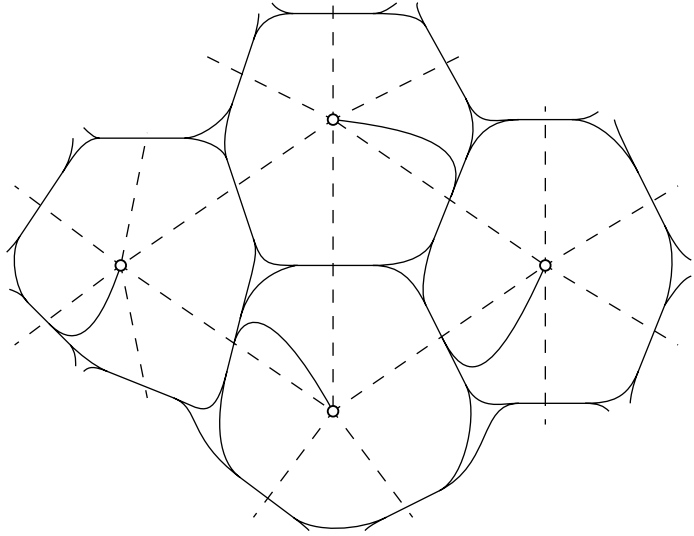}}
\caption{A universal train track $\theta_T$. Dotted lines are the edges of $T$}\label{triang2track}
\end{figure}

The result may look as shown in Figure~\ref{triang2track}.

\begin{prop}\label{minimalencoding}
\begin{enumerate}
\item
For any multiple curve $\gamma$ there is an isotopic multiple curve $\gamma'$ such
that $\theta_T$ carries $\gamma'$.
\item
Among multiple curves isotopic to $\gamma$ and carried by $\theta_T$ there is
a multiple curve $\gamma_{\min}$ that is minimal in the following sense:
if $(\theta_T,w_{\min})$ encodes $\gamma_{\min}$ and $(\theta_T,w)$ encodes
any other multiple curve isotopic to $\gamma$, then
$w_{\min}(\alpha)\leqslant w(\alpha)$ for any branch branch $\alpha\subset\theta_T$.
Clearly, such width assignment $w_{\min}$ is unique.
\item
There is a linear time algorithm to produce $w_{\min}$ from normal coordinates of $\gamma$, and we
have $|(\theta_T,w_{\min})|\leqslant C|\gamma|_T$, with $C$ not depending on $\gamma$,
where by $|\gamma|_T$ we denote the following complexity measure:
\begin{equation}\label{gamma-norm-eq}
|\gamma|_T=\sum_i\log_2(|\langle\gamma,e_i\rangle|+1).
\end{equation}
\end{enumerate}
\end{prop}

\begin{proof}
By construction, for every puncture $P$, we have a single branch of $\theta_T$ approaching~$P$.
Denote this branch by $\alpha_P$, and the triangle of $T$ containing $\alpha_P$
by $\tau_P$. If $P$ is an internal puncture, then the other end of $\alpha_P$ approaches a switch
from an ingoing side,
and the latter will be used to choose an orientation of $M$ at $P$, by which we
mean a sign designation to either rotation direction. Namely, if the cusp
at the switch occurs on the left when one travels along~$\alpha_P$ from $P$ to the switch,
then the counterclockwise direction will be positive and clockwise negative,
and vice versa if the cusp occurs on the right; see Figure~\ref{orient}.
\begin{figure}[ht]
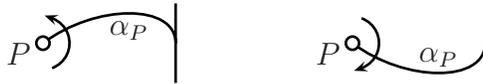

\centerline{\includegraphics{or.eps}\put(-69,7){$P$}\put(-30,18){$\alpha_P$}\hskip2cm
\includegraphics{or2.eps}\put(-69,7){$P$}\put(-30,8){$\alpha_P$}}
\caption{The positive rotation direction at an internal puncture $P$}\label{orient}
\end{figure}

We also choose an orientation of the surface at every boundary vertex $P$ so that a tangent vector
to the boundary $\partial M$ having positive direction will
point inward~$M$ after a small rotation in the positive direction around $P$; see Figure~\ref{orient-boundary}.
\begin{figure}[ht]
\centerline{\includegraphics[scale=.8]{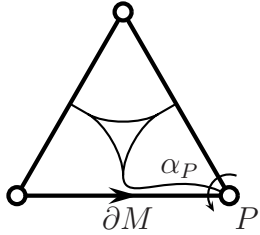}\put(-52,-6){$\partial M$}\put(-2,-6){$P$}\put(-30,15){$\alpha_P$}}
\caption{The positive rotation direction at a boundary vertex $P$}\label{orient-boundary}
\end{figure}

An arc in a triangle $\tau$ of $T$ will be called \emph{almost normal}
if it connects a vertex of $\tau$ with an interior point of an adjacent side; see Figure~\ref{almostnorm}.
\begin{figure}[ht]
\centerline{\includegraphics[scale=.8]{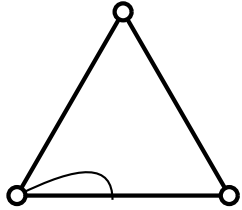}}
\caption{An almost normal arc}\label{almostnorm}
\end{figure}

Two arcs in a triangle $\tau$ of $T$ are called \emph{similar} if they are ambient isotopic in $\tau$
relative to the vertices of $\tau$. If an arc is similar to a smooth arc contained in $\theta_T$
we say that it is \emph{supported} by $\theta_T$.

Our train track $\theta_T$ is designed so that any side-to-side normal arc
is supported by $\theta_T$ (and in a unique way).
For a vertex-to-side or vertex-to-vertex normal
arc this is typically not true. Exceptions occur in triangles having one or two edges at the boundary,
see Figure~\ref{orient-boundary} and Figure~\ref{two-boundary-edges}.
\begin{figure}[ht]
\centerline{\includegraphics[scale=.8]{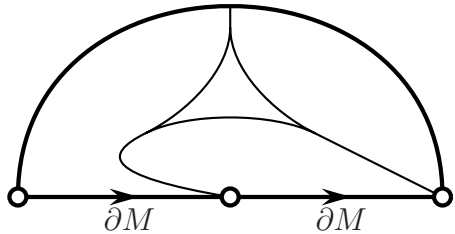}\put(-52,-6){$\partial M$}\put(-132,-6){$\partial M$}}
\caption{The train track $\theta_T$ in a triangle that has two edges at the boundary}\label{two-boundary-edges}
\end{figure}
The intersection of $\theta_T$ with such a triangle supports one
vertex-to-side arc and, in the case of two boundary edges, one vertex-to-vertex arc.

The idea behind the construction of $\gamma'$ is to normalize the original curve with
respect to~$T$ and then deform all normal arcs that are not supported by $\theta_T$ so
as to obtain a composition of normal and almost normal arcs that are supported.
After that we can push the obtained curve toward $\theta_T$ so that
all normal and almost normal arcs become close to the corresponding arcs in $\theta_T$.

In order to see how it works we start from the opposite side, i.e.\ from a
multiple curve $\gamma'$ that \emph{is} carried by $\theta_T$. Let $w$
be the corresponding width assignment to branches of $\theta_T$.

If $P$ is not a boundary puncture with just one triangle adjacent to it
(consult Figure~\ref{two-boundary-edges}), then there is a unique,
up to similarity,
almost normal arc attached to $P$ that is supported by~$\theta_T$.
It is obtained by a smooth extension of $\alpha_P$ along $\theta_T$
up to the boundary of the triangle. Denote this almost normal arc by
$\widetilde\alpha_P$.

If $w(\alpha_P)>0$, then $\gamma'$ contains an arc similar
to $\widetilde\alpha_P$, hence, it is not normal with
respect to~$T$ as $\widetilde\alpha_P$ cuts a bigon off $\tau_P$.
Now see what happens if we run the normalization
procedure for $\gamma'$.

The following assertions remain true during the normalization  process:
\begin{enumerate}
\item
at every normalization step the multiple curve $\gamma'$ is composed of normal and almost normal arcs;
\item
every bigon reduction results in rotating the tail of an almost
normal arc around the corresponding puncture in the negative direction.
\end{enumerate}

Indeed, it is easy to see that a bigon whose boundary is disjoint from punctures
cannot appear in the pulling tight process unless it was present at the beginning.
Figure~\ref{normalisation} demonstrates a single bigon reduction
for all possible types of arcs extending the almost normal arc
being reduces.
\begin{figure}[ht]
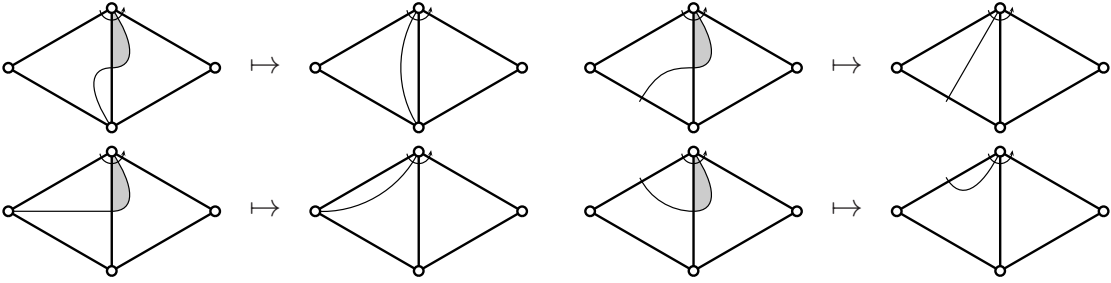

\begin{tabular}{ccccccc}
\includegraphics[scale=.45]{norm1.eps}&\raisebox{23pt}{$\mapsto$}&\includegraphics[scale=.45]{norm2.eps}&&
\includegraphics[scale=.45]{norm3.eps}&\raisebox{23pt}{$\mapsto$}&\includegraphics[scale=.45]{norm4.eps}\\
\includegraphics[scale=.45]{norm5.eps}&\raisebox{23pt}{$\mapsto$}&\includegraphics[scale=.45]{norm6.eps}&&
\includegraphics[scale=.45]{norm7.eps}&\raisebox{23pt}{$\mapsto$}&\includegraphics[scale=.45]{norm8.eps}
\end{tabular}
\caption{Reduction of a bigon cut off by an almost normal arc}\label{normalisation}
\end{figure}
In the first three cases, a normal arc is produced. In the last case, a new
almost normal arc appears, and it is `oriented' in the same way as the
original one meaning that a small rotation in the positive direction
around the puncture pushes it off the corresponding bigon.

Figure~\ref{normalisation} shows all possible ways in which a normal arc that is not supported
by $\theta_T$ may appear. So, it is clear how to invert this procedure.

Namely, we start from a normal curve $\gamma'$ isotopic to $\gamma$. Then
we keep repeating the following step until~$\gamma'$ is carried by~$\theta_T$:
if~$\gamma'$ contains an unsupported almost normal arc, we apply an isotopy to~$\gamma'$
that modifies such an arc by a transformation inverse to one of those
shown in Fig~\ref{normalisation}. It is not hard to see that the procedure will eventually stop.

There is an arbitrariness in the process affecting the result in the following two ways. First,
at some stages of the process there may be
more than one unsupported arc to modify and more than
one way to modify the chosen unsupported arc (the latter case
may occur for a normal vertex-to-vertex arc). 
Second, if $\gamma'$ contains proper arcs isotopic to edges of $T$, then
their initial position is not unique. Both issues are illustrated in Figure~\ref{denorm-fig}, where the universal train
track is shown in grey dashed line.
\begin{figure}[ht]
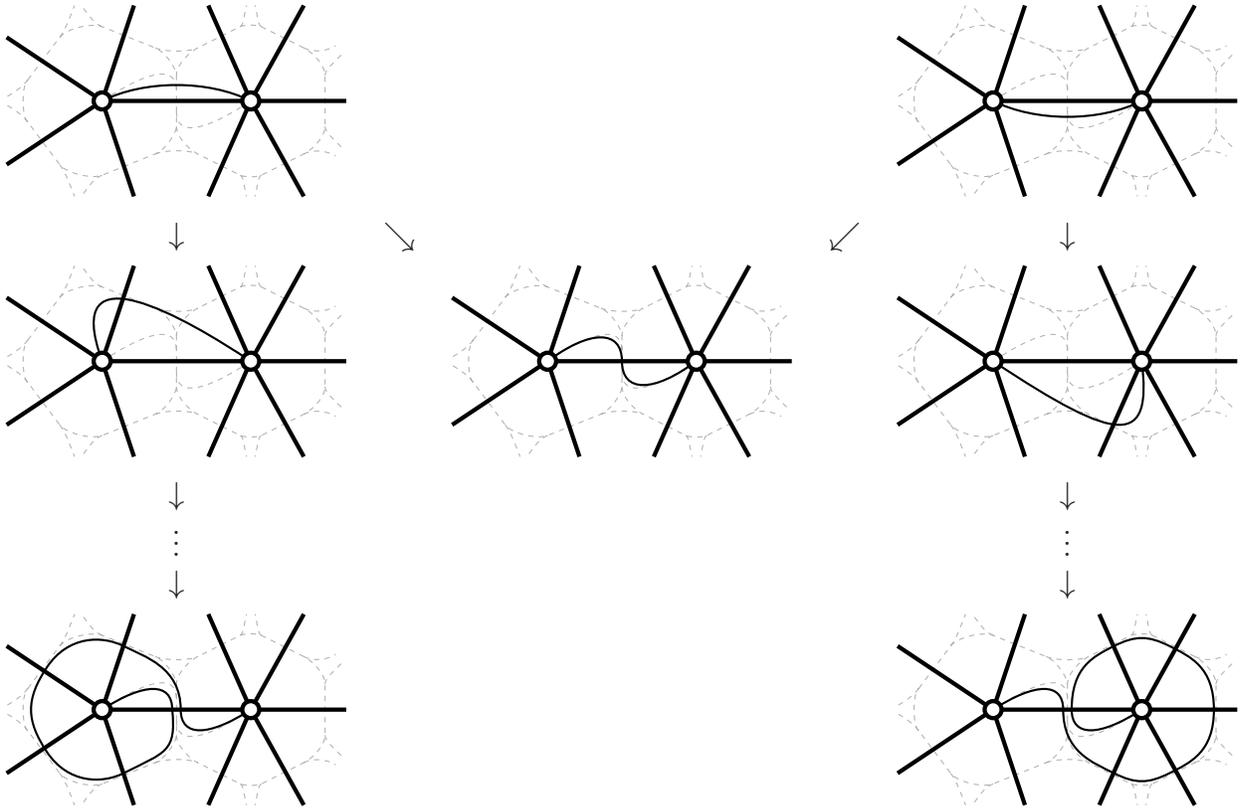

$$\begin{array}{ccccc}
\includegraphics[scale=.8]{denorm1.eps}&&&&\includegraphics[scale=.8]{denorm2.eps}\\
\downarrow&\searrow&&\swarrow&\downarrow\\
\includegraphics[scale=.8]{denorm4.eps}&&\includegraphics[scale=.8]{denorm3.eps}&&\includegraphics[scale=.8]{denorm6.eps}\\
\downarrow&&&&\downarrow\\
\vdots&&&&\vdots\\
\downarrow&&&&\downarrow\\
\includegraphics[scale=.8]{denorm5.eps}&&&&\includegraphics[scale=.8]{denorm7.eps}
\end{array}$$
\caption{Non-uniqueness of a multiple curve isotopic to a given one and carried by the universal train track}\label{denorm-fig}
\end{figure}
\begin{figure}
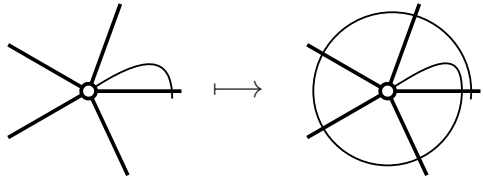

\centerline{\includegraphics[scale=.7]{circle1.eps}\quad\raisebox{33pt}{$\longmapsto$}\quad\includegraphics[scale=.7]{circle2.eps}}
\caption{Creation of an additional spiral turn}\label{circle}
\end{figure}

However, the possible results are not very much different from each other. Namely, each one can be obtained from any other
by creating and/or removing additional spiral turns like the one shown in Figure~\ref{circle}.

Clearly, each spiral turn contributes positively into widths of the branches
of $\theta_T$ that make a full turn around a puncture, so, in order to minimize
the widths we must avoid the spiral turns.
The branches $\alpha_P$ were constructed so that the spiral turns around
different punctures do not overlap. So, there is always a unique way (up to
isotopy preserving the triangulation) to remove them,
which gives the sought-for~$\gamma_{\min}$.

Computing the width assignment $w_{\min}$ corresponding to $\gamma_{\min}$ is
now very simple. There are only finitely many different types of normal arcs.
For each of them we implement the procedure described above and
find an isotopic arc decomposed in the optimal way into normal and almost normal arcs
supported by $\theta_T$. In this way the contribution of each normal
arc type into the width assignment is computed and recorded. This is
done only once, before any multiple curve is given.
Note that a single normal arc of the normalized
form of the original curve contributes at most two to the width of
any branch of $\theta_T$.

Then, given the normal coordinates of a multiple curve $\gamma$
one computes the number of normal arcs of each type (as described in the
proof of Lemma~\ref{lm1}) and sums up their contributions.
The running time estimation and that for the complexity of the result are straightforward.
\end{proof}

\begin{exam}Figure~\ref{5edges} illustrates how the curve $\gamma_{\min}$ and the corresponding
width assignment (the non-zero widths) look like for $\gamma$
the union of the five edges connecting the four punctures in Figure~\ref{triang2track},
where the choice of $\theta_T$ is shown in grey dashed line.
\begin{figure}[ht]
\centerline{\includegraphics{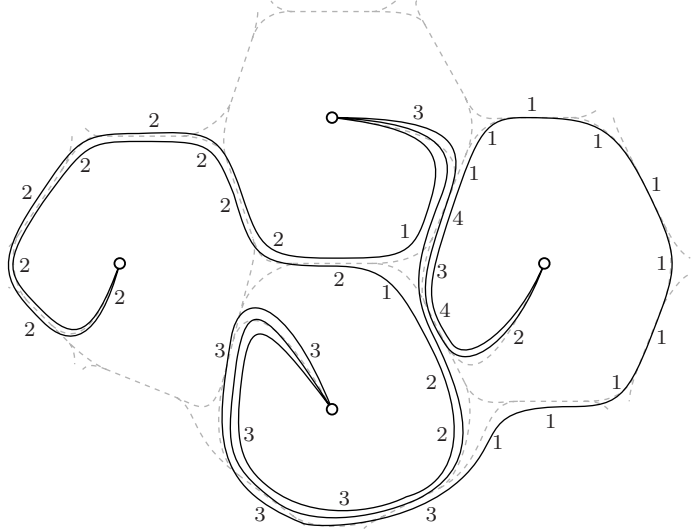}\put(-233,85){\tiny2}
\put(-267,73){\tiny2}\put(-269,97){\tiny2}\put(-268,125){\tiny2}\put(-246,135){\tiny2}
\put(-220,152){\tiny2}\put(-202,137){\tiny2}\put(-193,120){\tiny2}\put(-173,107){\tiny2}
\put(-150,92){\tiny2}\put(-125,110){\tiny1}\put(-105,115){\tiny4}\put(-120,155){\tiny3}
\put(-111,95){\tiny3}\put(-132,87){\tiny1}\put(-110,80){\tiny4}\put(-115,53){\tiny2}
\put(-82,70){\tiny2}\put(-99,132){\tiny1}\put(-92,145){\tiny1}\put(-77,158){\tiny1}
\put(-52,144){\tiny1}\put(-30,128){\tiny1}\put(-28,98){\tiny1}\put(-28,70){\tiny1}
\put(-45,53){\tiny1}\put(-70,38){\tiny1}\put(-90,30){\tiny1}\put(-111,33){\tiny2}
\put(-115,3){\tiny3}\put(-148,9){\tiny3}\put(-180,3){\tiny3}\put(-184,33){\tiny3}
\put(-195,65){\tiny3}\put(-159,65){\tiny3}}
\caption{Five edges of the triangulation $T$ put in a minimal position carried
by $\theta_T$}\label{5edges}
\end{figure}
\end{exam}

\section{Simplifying train tracks}\label{simplifying}
Simplification procedure introduced in this section is one of the many
similar ones that mimic the accelerated Euclidean algorithm.
The general principle for constructing such algorithms in
low-dimensional topology settings was learnt by the author
from the work of I.\,Agol, J.\,Hass, and W.\,Thurston~\cite{aht}.

Here we describe transformations $(\theta,w)\mapsto(\theta',w')$
of measured train tracks
such that the multiple curves encoded by $(\theta,w)$ and $(\theta',w')$
are isotopic.
To every such transformation we assign two numbers
that are called \emph{the gain} and \emph{the cost}
of the transformation. Vaguely speaking, the former indicates how much $(\theta',w')$
is simpler than $(\theta,w)$, and the latter measures `the algorithmic complexity'
of the operation.

Recall that by complexity $|(\theta,w)|$ of a measured train track $(\theta,w)$ we mean the sum
\begin{equation}\label{|thetaw|-eq}|(\theta,w)|=\sum_{\alpha}\bigl(1+\log_2(w(\alpha)+1)\bigr),
\end{equation}
taken over all branches of $\theta$.
It is comparable to the amount of space needed to encode $(\theta,w)$.
However, for technical reasons,
we will need a slightly more subtle measure of complexity.

Denote by $A(\theta)$ the set of non-free branches of $\theta$. Define
\begin{equation}\label{complexity0}
|(\theta,w)|_0=|A(\theta)|+\sum_{\alpha\in A(\theta)}\log_2(w(\alpha)+1).
\end{equation}
This is obtained from~\eqref{|thetaw|-eq} by dropping the contribution of free branches.

Whatever a transformation $(\theta,w)\mapsto(\theta',w')$
is \emph{the gain} of this transformation is defined as
the difference $|(\theta,w)|_0-|(\theta',w')|_0$. If we have a sequence
$$(\theta_0,w_0)\mapsto(\theta_1,w_1)\mapsto\ldots\mapsto(\theta_k,w_k)$$
of transformations, then \emph{the total gain} of the sequence
is set to $|(\theta_0,w_0)|_0-|(\theta_k,w_k)|_0$.

Now we introduce transformations $(\theta,w)\mapsto(\theta',w')$ of our interest.
They will be referred to as \emph{simplification moves} and
include removing trivial branches, splittings (ordinary and multiple),
and slidings defined below.

\smallskip\noindent\emph{Removing trivial branches}.
The train track $\theta'$ is obtained from $\theta$ by removing
all non-free branches $\alpha$ such that $w(\alpha)=0$.
If~$\theta$ contains a switch
to which exactly one free branch of~$\theta$ is attached, then the other two branches
approaching this switch become parts of a single branch of~$\theta'$.
The width $w'(\alpha')$ of any branch $\alpha'$ of~$\theta'$ is
set to $w(\alpha)$ with any branch $\alpha$
of $\theta$ such that $\alpha\subset\alpha'$ (clearly the choice of $\alpha$ does not matter).

We set \emph{the cost} of
this operation to be equal to the number of non-free branches
$\alpha\subset\theta$ such that $w(\alpha)=0$.

\smallskip\noindent\emph{Ordinary splitting}.
Recall that a branch $\alpha$ of $\theta$ is called \emph{wide} in the following two cases:
\begin{enumerate}
\item[(1)]
both tails of $\alpha$ are outgoing for some switches;
\item[(2)]
one tail of $\alpha$ is outgoing, and the other approaches a puncture.
\end{enumerate}

\begin{figure}[ht]
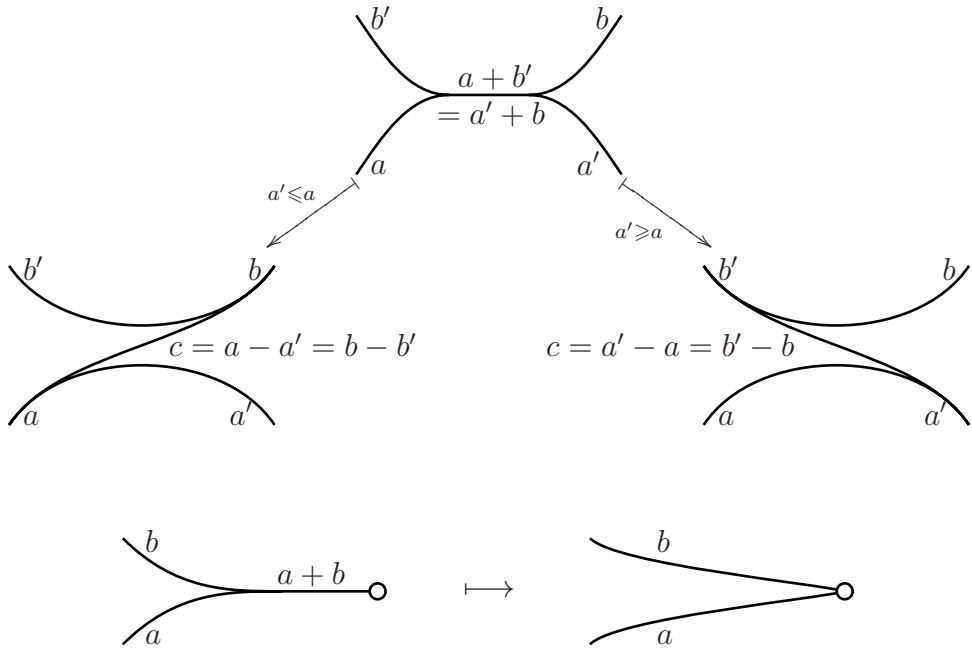

\centerline{\xymatrix{%
&
\raisebox{-30pt}{\includegraphics{split1.eps}\put(-95,0){$a$}\put(-95,55){$b'$}\put(-17,0){$a'$}\put(-10,55){$b$}\put(-62,33){$a+b'$}\put(-70,19){$=a'+b$}}
\ar@{|->}[dr]_{a'\geqslant a}\ar@{|->}[dl]_{a'\leqslant a}&\\
\raisebox{-30pt}{\includegraphics{split2.eps}\put(-95,0){$a$}\put(-95,55){$b'$}\put(-17,0){$a'$}\put(-10,55){$b$}
\put(-40,26){$c=a-a'=b-b'$}}
&&\raisebox{-30pt}{\includegraphics{split2a.eps}\put(-95,0){$a$}\put(-95,55){$b'$}\put(-17,0){$a'$}\put(-10,55){$b$}
\put(-160,26){$c=a'-a=b'-b$}}
}}
\vskip1cm
\centerline{\includegraphics{split3.eps}\put(-92,10){$a$}%
\put(-92,45){$b$}\put(-42,33){$a+b$}
\hskip1cm\raisebox{28pt}{$\longmapsto$}\hskip1cm
\includegraphics{split4.eps}\put(-75,10){$a$}\put(-75,45){$b$}}
\caption{Ordinary splittings}\label{ordsplit}
\end{figure}

\emph{An ordinary splitting} $(\theta,w)\mapsto(\theta',w')$ on
a wide branch $\alpha$ is a modification of the measured train track
$(\theta,w)$ that occurs in a small neighborhood of $\alpha$ and
has the form shown in Figure~\ref{ordsplit}, where
widths of the involved branches are also indicated.
Widths of the other branches are preserved.

To every ordinary splitting we assign \emph{cost}~$1$.

\smallskip\noindent\emph{Multiple splitting}.
Suppose that the train track $\theta$ contains two branches $\beta$ and $\gamma$,
say, whose union is a two-sided simple closed curve. There must be two tails outside of
$\beta\cup\gamma$ that approach switches at $\beta\cup\gamma$.
We additionally suppose that they do it from different sides
of $\beta\cup\gamma$. Finally, we suppose $w(\beta)<w(\gamma)\leqslant2w(\beta)$.

\begin{figure}[ht]
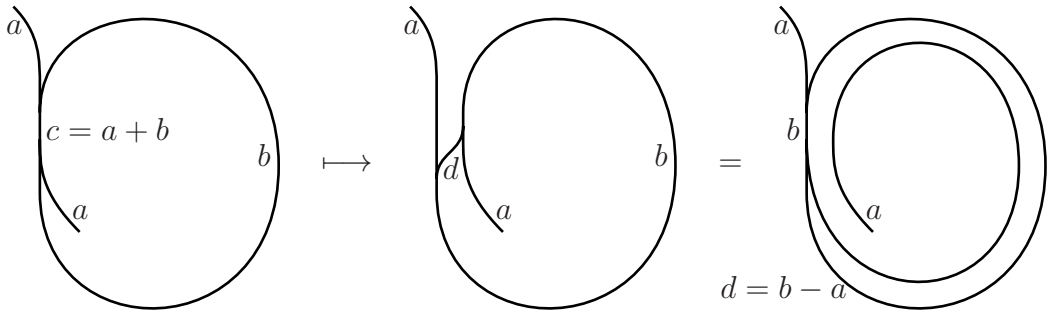

\vskip1cm
\centerline{\includegraphics{spiral1.eps}\put(-10,60){$b$}\put(-80,40){$a$}
\put(-105,110){$a$}\put(-90,70){$c=a+b$}\hskip.5cm\raisebox{58pt}{$\longmapsto$}\hskip.5cm
\includegraphics{spiral2.eps}\put(-10,60){$b$}\put(-70,40){$a$}\put(-90,55){$d$}
\put(-105,110){$a$}\hskip.5cm\raisebox{58pt}{$=$}\hskip.5cm
\includegraphics{spiral3.eps}\put(-100,70){$b$}\put(-70,40){$a$}
\put(-105,110){$a$}\put(-125,10){$d=b-a$}}
\caption{If $b/a\geqslant k\in\mathbb N$, then we can apply $k$
splittings at once, which will be a multiple splitting}\label{multsplit}
\end{figure}

Then $\gamma$ must be a wide branch, and we have a
situation shown in Figure~\ref{multsplit} on the left,
where the widths $b=w(\beta)$, $c=w(\gamma)$, $a=c-b$
are indicated near the respective branches.
By assumption, we have $c\leqslant2b$, hence $b\geqslant a$.
After a splitting on the branch $\gamma$
we get a measured train track $(\theta',w')$ that is obtained from $(\theta,w)$ by a Dehn twist
along $\beta\cup\gamma$ and making the branches $\beta$, $\gamma$ narrower by $a$,
see Figure~\ref{multsplit}. So, if $b\geqslant ka$, $k\in\mathbb N$,
we can apply $k$ successive splittings to this portion of $\theta$,
which result in the application of the $k$th power of a Dehn twist
along $\alpha\cup\beta$ to $\theta$ and
making the branches $\alpha$ and $\beta$ narrower by $ka$.

Such application of $k$ successive splittings will be treated as a single
operation called \emph{a $k$-times multiple splitting on the circle $\beta\cup\gamma$}. Its \emph{cost}
is set to $\log_2(k+1)$.

\smallskip\noindent\emph{Sliding}. Let $\alpha$ be a branch
of $\theta$ having one ingoing and one outgoing tail.
\emph{A sliding along~$\alpha$} is a modification of $(\theta,w)$
that occurs in a small neighborhood of $\alpha$
\begin{figure}[ht]
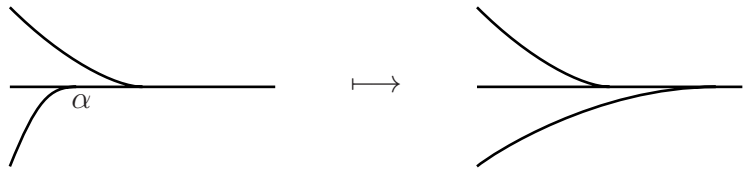

\centerline{\includegraphics{slide1.eps}\put(-77,22){$\alpha$}
\hskip1cm\raisebox{28pt}{$\longmapsto$}\hskip1cm\includegraphics{slide2.eps}}
\caption{A sliding}\label{slide}
\end{figure}
as shown in Figure~\ref{slide}. The widths of all branches outside
of this neighborhood are preserved.
The \emph{cost} is set to $1$.
Note that unlike other simplification moves a sliding may
have a negative gain, so, sometimes
it does not justify the name `simplification move'.

If $(\theta_0,w_0),(\theta_1,w_1),\ldots,(\theta_k,w_k)$ is
a sequence of measured train tracks in which
every transition $(\theta_i,w_i)\mapsto(\theta_{i+1},w_{i+1})$ is a simplification move,
then the sum of their costs is called \emph{the total cost} of the sequence.

\begin{prop}\label{costbound}
Let $(\theta,w)$ be a measured train track. Then there
exists a sequence of simplification moves
starting from $(\theta,w)$ and ending with a measured train track
without switches, such that the total cost of the
sequence does not exceed $3\cdot|(\theta,w)|$.

There is an algorithm that produces such a sequence in $O(|(\theta,w)|^2)$ operations
on a RAM machine.
\end{prop}

\begin{proof}
We prove the first statement with $|(\theta,w)|$ replaced by $|(\theta,w)|_0$,
which is stronger as we always have $|(\theta,w)|\geqslant|(\theta,w)|_0$.
We proceed by induction in $[|(\theta,w)|_0]$, where $[\ ]$ stands for the integral part.
The equality $[|(\theta,w)|_0]=0$ means that $\theta$ has no switches, and we are done.

For the induction step we just need to find a sequence
$(\theta,w)=(\theta_0,w_0)\mapsto(\theta_1,w_1)\mapsto\ldots\mapsto(\theta_l,w_l)$
of simplification moves such that its total gain $g$ and total
cost $p$ satisfy the following inequalities:
$$g\geqslant1,\quad g\geqslant p/3.$$

If $(\theta,w)$ has trivial branches, we remove them,
which gives $g\geqslant p\geqslant 1$. In the sequel we assume that all non-free branches
have positive widths.

If there is a wide branch that is attached to a puncture (see the lower part
of Figure~\ref{ordsplit}) we apply an ordinary splitting on it,
which gives $g>p=1$. In the sequel we assume that there is no such branch.

Among all wide branches of $(\theta,w)$ choose a widest
one $\alpha$, say, i.e.\ having the largest width.

By the assumption we have just made, both ends of $\alpha$ are switches. An ordinary splitting on~$\alpha$
may then have arbitrarily small gain, so, our strategy will depend on the structure
of~$(\theta,w)$ around $\alpha$. We consider below a bunch of cases that are summarized
in Figure~\ref{cases}, where $\alpha$ is the branch that has width $a+b+c$.
The sign `$\bigotimes$' in the pictures denotes an orientation flip.

\begin{figure}[ht]
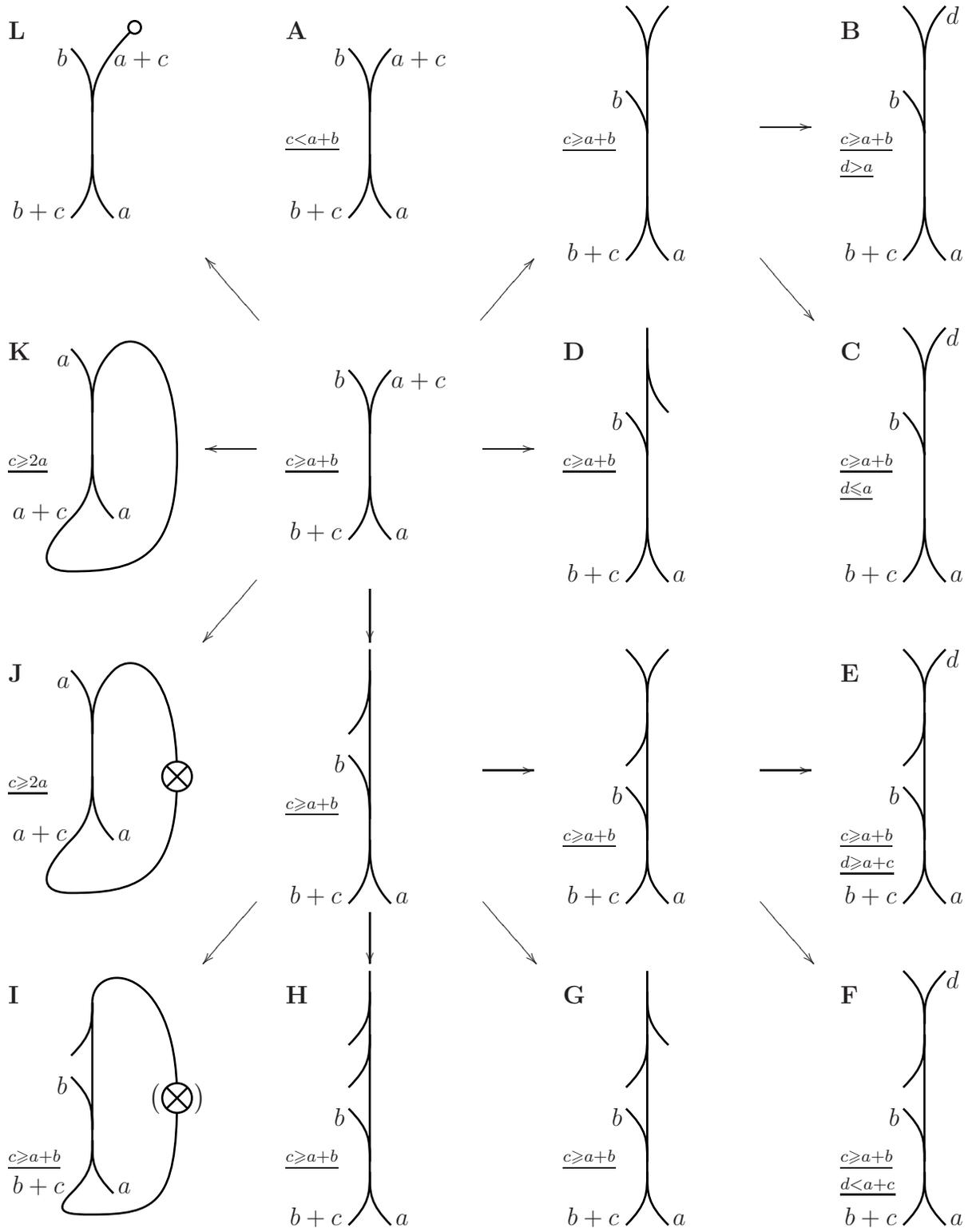

\xymatrix{
\raisebox{-60pt}{\includegraphics{case3.eps}\put(-88,20){$b+c$}%
\put(-38,20){$a$}\put(-67,92){$b$}\put(-40,92){$a+c$}
\put(-90,105){\bf L}}%
&
\raisebox{-60pt}{\includegraphics{case0.eps}\put(-88,20){$b+c$}%
\put(-38,20){$a$}\put(-67,92){$b$}\put(-40,92){$a+c$}
\put(-90,105){\bf A}\put(-90,55){\underline{$\scriptstyle c<a+b$}}}
&
\raisebox{-60pt}{\includegraphics{case21.eps}\put(-88,0){$b+c$}\put(-38,0){$a$}\put(-67,72){$b$}
\put(-90,55){\underline{$\scriptstyle c\geqslant a+b$}}}\ar[r]\ar[dr]
&
\raisebox{-60pt}{\includegraphics{case21.eps}\put(-88,0){$b+c$}\put(-38,0){$a$}\put(-67,72){$b$}
\put(-40,112){$d$}
\put(-90,105){\bf B}
\put(-90,55){\underline{$\scriptstyle c\geqslant a+b$}}
\put(-90,42){\underline{$\scriptstyle d>a$}}}
\\
\raisebox{-60pt}{\includegraphics{case25.eps}\put(-88,30){$a+c$}%
\put(-38,30){$a$}\put(-67,102){$a$}
\put(-90,105){\bf K}
\put(-90,55){\underline{$\scriptstyle c\geqslant2a$}}}
&
\raisebox{-60pt}{\includegraphics{case0.eps}\put(-88,20){$b+c$}%
\put(-38,20){$a$}\put(-67,92){$b$}\put(-40,92){$a+c$}
\put(-90,55){\underline{$\scriptstyle c\geqslant a+b$}}}\ar[ur]\ar[r]\ar[l]\ar[dl]\ar[d]\ar[ul]
&
\raisebox{-60pt}{\includegraphics{case22.eps}\put(-88,0){$b+c$}\put(-38,0){$a$}\put(-67,72){$b$}
\put(-90,105){\bf D}
\put(-90,55){\underline{$\scriptstyle c\geqslant a+b$}}}
&
\raisebox{-60pt}{\includegraphics{case21.eps}\put(-88,0){$b+c$}\put(-38,0){$a$}\put(-67,72){$b$}
\put(-40,112){$d$}
\put(-90,55){\underline{$\scriptstyle c\geqslant a+b$}}
\put(-90,105){\bf C}
\put(-90,42){\underline{$\scriptstyle d\leqslant a$}}}
\\
\raisebox{-60pt}{\includegraphics{case24.eps}\put(-88,30){$a+c$}\put(-38,30){$a$}\put(-67,102){$a$}
\put(-90,105){\bf J}
\put(-90,55){\underline{$\scriptstyle c\geqslant 2a$}}}
&
\raisebox{-60pt}{\includegraphics{case23.eps}\put(-88,0){$b+c$}\put(-38,0){$a$}\put(-67,62){$b$}
\put(-90,45){\underline{$\scriptstyle c\geqslant a+b$}}}\ar[r]\ar[d]\ar[dr]\ar[dl]
&
\raisebox{-60pt}{\includegraphics{case233.eps}\put(-88,0){$b+c$}\put(-38,0){$a$}\put(-67,47){$b$}
\put(-90,30){\underline{$\scriptstyle c\geqslant a+b$}}}\ar[r]\ar[dr]
&
\raisebox{-60pt}{\includegraphics{case233.eps}\put(-88,0){$b+c$}\put(-38,0){$a$}\put(-67,47){$b$}
\put(-40,112){$d$}
\put(-90,105){\bf E}
\put(-90,17){\underline{$\scriptstyle d\geqslant a+c$}}
\put(-90,30){\underline{$\scriptstyle c\geqslant a+b$}}}
\\
\raisebox{-60pt}{\includegraphics{case234.eps}\put(-88,15){$b+c$}\put(-38,15){$a$}\put(-67,62){$b$}
\put(-90,105){\bf I}
\put(-23.5,57){$\bigl($}\put(-2,57){$\bigr)$}
\put(-90,30){\underline{$\scriptstyle c\geqslant a+b$}}}
&
\raisebox{-60pt}{\includegraphics{case231.eps}\put(-88,0){$b+c$}\put(-38,0){$a$}\put(-67,47){$b$}
\put(-90,105){\bf H}
\put(-90,30){\underline{$\scriptstyle c\geqslant a+b$}}}
&
\raisebox{-60pt}{\includegraphics{case232.eps}\put(-88,0){$b+c$}\put(-38,0){$a$}\put(-67,47){$b$}
\put(-90,105){\bf G}
\put(-90,30){\underline{$\scriptstyle c\geqslant a+b$}}}
&
\raisebox{-60pt}{\includegraphics{case233.eps}\put(-88,0){$b+c$}\put(-38,0){$a$}\put(-67,47){$b$}
\put(-40,112){$d$}
\put(-90,105){\bf F}
\put(-90,17){\underline{$\scriptstyle d<a+c$}}
\put(-90,30){\underline{$\scriptstyle c\geqslant a+b$}}}
}
\caption{The chart of simplification cases}\label{cases}
\end{figure}

\noindent{\bf Case A.}
We have $a+b\geqslant c+1$, hence, for an ordinary splitting on~$\alpha$, we have
$$2^g=\frac{a+b+c+1}{c+1}\geqslant2.$$
Thus, $g\geqslant1=p$.

\smallskip

\noindent{\bf Case B.}
Let $e=d-a>0$, $f=c-e>0$; consult Figure~\ref{caseb}. After one sliding and one ordinary splitting
the two branches of width
\begin{figure}[ht]
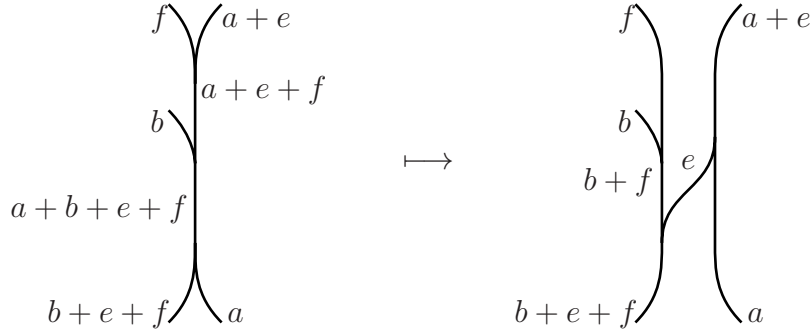

\centerline{\includegraphics{case21.eps}\put(-106,0){$b+e+f$}%
\put(-38,0){$a$}\put(-67,72){$b$}\put(-67,112){$f$}\put(-40,112){$a+e$}
\put(-120,40){$a+b+e+f$}\put(-48,85){$a+e+f$}
\hskip1cm\raisebox{58pt}{$\longmapsto$}\hskip1cm
\includegraphics{case21a.eps}\put(-126,0){$b+e+f$}\put(-38,0){$a$}\put(-87,72){$b$}
\put(-87,112){$f$}\put(-40,112){$a+e$}
\put(-100,50){$b+f$}\put(-63,58){$e$}}
\caption{Simplification in Case B}\label{caseb}
\end{figure}
$a+b+e+f$ and $a+e+f$ are replaced by those of width $b+f$ and $e$.
We have
$$\begin{aligned}
2^g&=
\frac{(a+b+e+f+1)(a+e+f+1)}{(b+f+1)(e+1)}\geqslant
\frac{(a+b+e+f+1)(a+e+f+1)}{(a+b+f)(e+1)}\\
&\geqslant
\frac{4(a+b+e+f+1)(a+e+f+1)}{(a+b+e+f+1)^2}=\frac{4(a+e+f+1)}{(a+b+e+f+1)}
>\frac{2(a+e+f+1)}{e+f+1}>2
\end{aligned}
$$
as $e+f=c\geqslant a+b$ in this case. Thus, $g>1$, $p=2$.

\smallskip
\noindent{\bf Case C.}
Let $e=a-d$; see Figure~\ref{casec}. We apply two ordinary splittings, which give
\begin{figure}[ht]
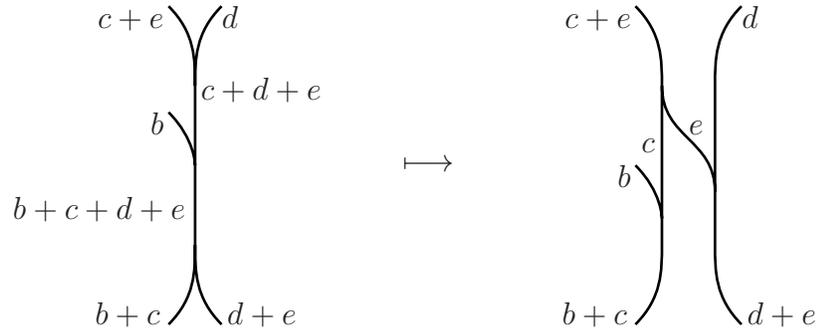

\centerline{\includegraphics{case21.eps}\put(-88,0){$b+c$}\put(-38,0){$d+e$}\put(-67,72){$b$}
\put(-40,112){$d$}\put(-119,40){$b+c+d+e$}\put(-48,85){$c+d+e$}\put(-87,112){$c+e$}
\hskip1cm\raisebox{58pt}{$\longmapsto$}\hskip1cm
\includegraphics{case21b.eps}\put(-108,0){$b+c$}\put(-38,0){$d+e$}\put(-87,52){$b$}
\put(-40,112){$d$}\put(-107,112){$c+e$}\put(-78,65){$c$}\put(-60,72){$e$}}
\caption{Simplification in Case C}\label{casec}
\end{figure}
$$2^g=\frac{(b+c+d+e+1)(c+d+e+1)}{(c+1)(e+1)}>2$$
as $c\geqslant a+b=b+d+e>e+1$. Thus, we have $g>1$, $p=2$.

\smallskip
\noindent{\bf Case D.}
We apply an ordinary splitting and a sliding as shown in Figure~\ref{cased}.
Since the branch $\alpha$ is the widest one, we have $d\leqslant b$.
Together with $c\geqslant a+b$ this gives:
\begin{figure}[ht]
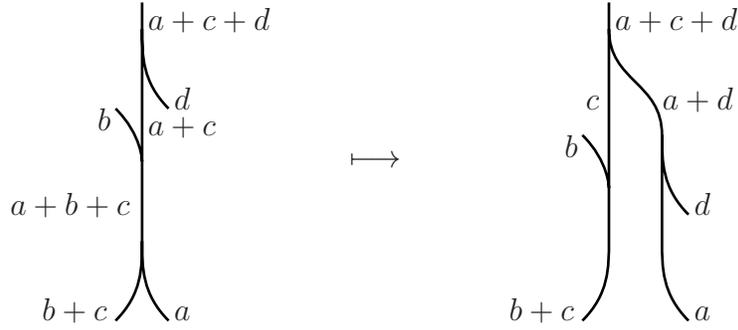

\centerline{\includegraphics{case22.eps}\put(-88,0){$b+c$}%
\put(-38,0){$a$}\put(-67,72){$b$}\put(-38,80){$d$}
\put(-100,40){$a+b+c$}\put(-48,70){$a+c$}\put(-48,110){$a+c+d$}
\hskip1cm\raisebox{58pt}{$\longmapsto$}\hskip1cm
\includegraphics{case22a.eps}\put(-108,0){$b+c$}\put(-38,0){$a$}\put(-87,62){$b$}\put(-38,40){$d$}
\put(-79,80){$c$}\put(-50,80){$a+d$}\put(-68,110){$a+c+d$}}
\caption{Simplification in Case D}\label{cased}
\end{figure}
$$\begin{aligned}
2^g&=
\frac{(a+b+c+1)(a+c+1)}{(a+d+1)(c+1)}\geqslant
\frac{(a+b+c+1)(a+c+1)}{(a+b+1)(c+1)}>
\frac{(a+b+c+2)(a+c)}{(a+b+1)(c+1)}\\
&\geqslant\frac{4(a+c)}{a+b+c+2}\geqslant\frac{2(a+c)}{c+1}\geqslant2.
\end{aligned}$$
Hence, in this case, $g>1$, $p=2$.

\smallskip
\noindent{\bf Case E.}
Consult Figure~\ref{casee} for notation. We have $c\geqslant a+b$. Since $\alpha$ is
a widest branch we also have $e+f\leqslant b$.
\begin{figure}[ht]
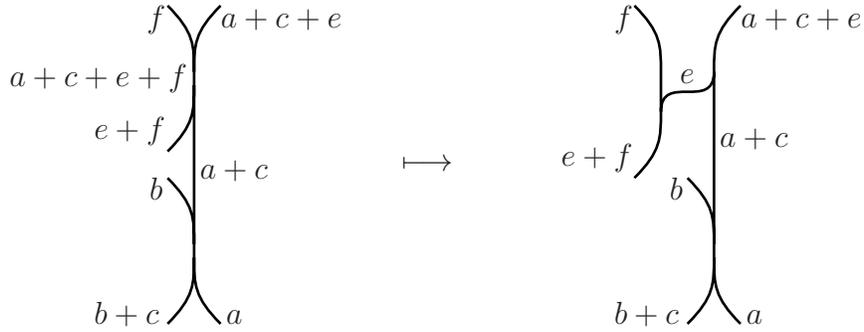

\centerline{\includegraphics{case233.eps}\put(-88,0){$b+c$}\put(-38,0){$a$}\put(-67,47){$b$}
\put(-40,112){$a+c+e$}\put(-68,112){$f$}%
\put(-88,70){$e+f$}\put(-48,55){$a+c$}\put(-120,90){$a+c+e+f$}
\hskip1cm\raisebox{58pt}{$\longmapsto$}\hskip1cm
\includegraphics{case233a.eps}\put(-88,0){$b+c$}\put(-38,0){$a$}\put(-67,47){$b$}
\put(-40,112){$a+c+e$}\put(-88,112){$f$}\put(-108,60){$e+f$}\put(-48,67){$a+c$}\put(-63,91){$e$}}
\caption{Simplification in Case E}\label{casee}
\end{figure}
Thus, we have $e<b<c<a+c+f$, which implies that splitting on the upper wide branch
in Figure~\ref{casee} will do the job (the situation is identical to Case A).
We will have $g\geqslant1=p$.

\smallskip
\noindent{\bf Case F.} After one splitting we come to the situation
\begin{figure}[ht]
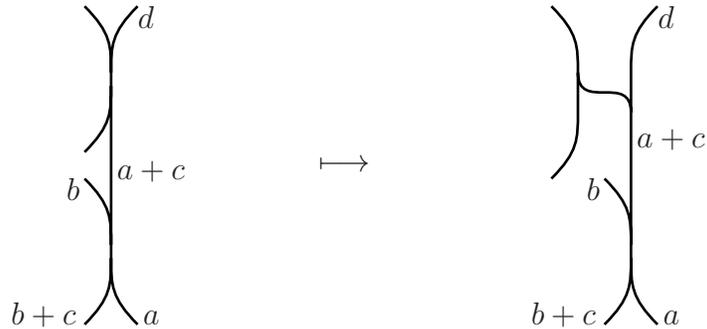

\centerline{\includegraphics{case233.eps}\put(-88,0){$b+c$}\put(-38,0){$a$}\put(-67,47){$b$}
\put(-40,112){$d$}\put(-48,55){$a+c$}
\hskip1cm\raisebox{58pt}{$\longmapsto$}\hskip1cm
\includegraphics{case233b.eps}\put(-88,0){$b+c$}\put(-38,0){$a$}\put(-67,47){$b$}
\put(-40,112){$d$}\put(-48,67){$a+c$}}
\caption{Simplification in Case F}\label{casef}
\end{figure}
of Case B or Case C; see Figure~\ref{casef}. Thus, after applying two
more simplification moves we have $g>1$, $p=3$.

\smallskip
\noindent{\bf Case G.}
Consult Figure~\ref{caseg} for notation.
Since there are no branches wider than $a+b+c$, we must have $d+e\leqslant b$.
\begin{figure}
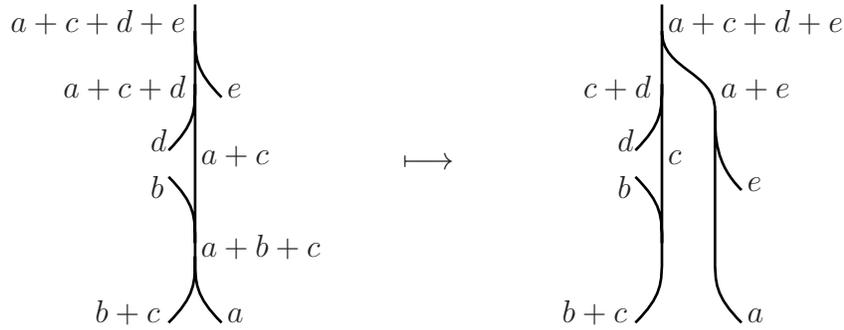

\centerline{\includegraphics{case232.eps}\put(-88,0){$b+c$}\put(-38,0){$a$}\put(-67,47){$b$}
\put(-48,25){$a+b+c$}\put(-48,60){$a+c$}\put(-67,65){$d$}\put(-100,85){$a+c+d$}
\put(-38,85){$e$}\put(-120,110){$a+c+d+e$}
\hskip1cm\raisebox{58pt}{$\longmapsto$}\hskip1cm
\includegraphics{case232a.eps}\put(-108,0){$b+c$}\put(-38,0){$a$}\put(-87,47){$b$}
\put(-48,85){$a+e$}\put(-68,60){$c$}\put(-87,65){$d$}\put(-100,85){$c+d$}
\put(-38,50){$e$}\put(-68,110){$a+c+d+e$}}
\caption{Simplification in Case G}\label{caseg}
\end{figure}
An ordinary splitting followed by two slidings gives:
$$\begin{aligned}
2^g&=
\frac{(a+b+c+1)(a+c+1)(a+c+d+1)}{(c+1)(a+e+1)(c+d+1)}\geqslant
\frac{(a+b+c+1)(a+c+1)(a+c+d+1)}{(c+1)(a+b)(c+d+1)}\\
&\geqslant
\frac{4(a+c+1)(a+c+d+1)}{(a+b+c+1)(c+d+1)}>
\frac{2(a+c+1)(a+c+d+1)}{(c+1)(c+d+1)}>2.
\end{aligned}$$
Thus, we have $g>1$, $p=3$.

\smallskip
\noindent{\bf Case H.}
Since there are no branches wider than $a+b+c$, we have
$d+e\leqslant b$ (see Figure~\ref{caseh}). After an ordinary splitting and a sliding shown in Figure~\ref{caseh} we have
\begin{figure}[ht]
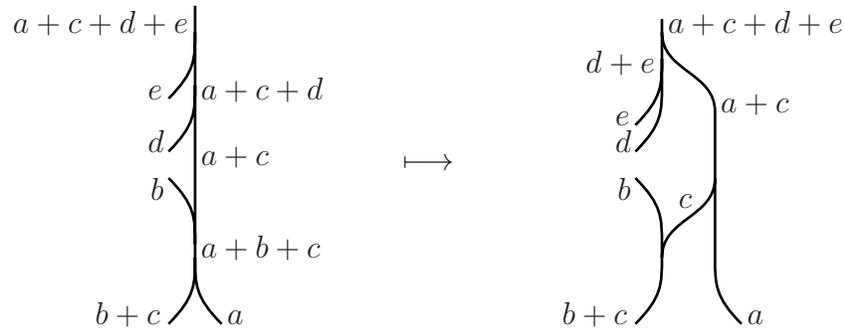

\centerline{\includegraphics{case231.eps}\put(-88,0){$b+c$}
\put(-38,0){$a$}\put(-67,47){$b$}\put(-68,65){$d$}
\put(-68,85){$e$}\put(-48,60){$a+c$}\put(-48,85){$a+c+d$}
\put(-119,110){$a+c+d+e$}\put(-48,25){$a+b+c$}
\hskip1cm\raisebox{58pt}{$\longmapsto$}\hskip1cm
\includegraphics{case231a.eps}\put(-108,0){$b+c$}
\put(-38,0){$a$}\put(-87,47){$b$}\put(-88,65){$d$}
\put(-88,75){$e$}\put(-64,44){$c$}\put(-48,80){$a+c$}
\put(-68,110){$a+c+d+e$}\put(-99,95){$d+e$}}
\caption{Simplification in Case H}\label{caseh}
\end{figure}
$$2^g=\frac{(a+b+c+1)(a+c+d+1)}{(c+1)(d+e+1)}\geqslant
\frac{(a+b+c+1)(a+c+d+1)}{(b+1)(c+1)}>2
$$
as $c\geqslant a+b$. So, $g>1$, $p=2$.

\smallskip
\noindent{\bf Case I.} The fragment of $(\theta,w)$
contains the configuration symmetric to that
covered by cases~B and~C; see Figure~\ref{casei}.
\begin{figure}[ht]
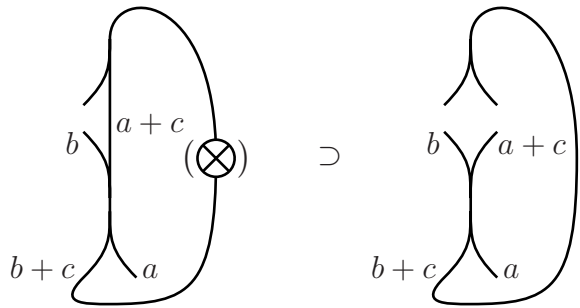

\centerline{\includegraphics{case234.eps}\put(-88,15){$b+c$}
\put(-38,15){$a$}\put(-67,62){$b$}\put(-48,70){$a+c$}
\put(-23.5,57){$\bigl($}\put(-2,57){$\bigr)$}
\hskip1cm\raisebox{58pt}{$\supset$}\hskip1cm
\includegraphics{case234a.eps}\put(-118,15){$b+c$}
\put(-68,15){$a$}\put(-97,62){$b$}\put(-70,62){$a+c$}}
\caption{Case I reduces to Cases B and C}\label{casei}
\end{figure}

\smallskip
\noindent{\bf Case J.} After two splittings we can remove
a trivial branch, which gives $g>3$, $p=3$ (see Figure~\ref{casej}).
\begin{figure}[ht]
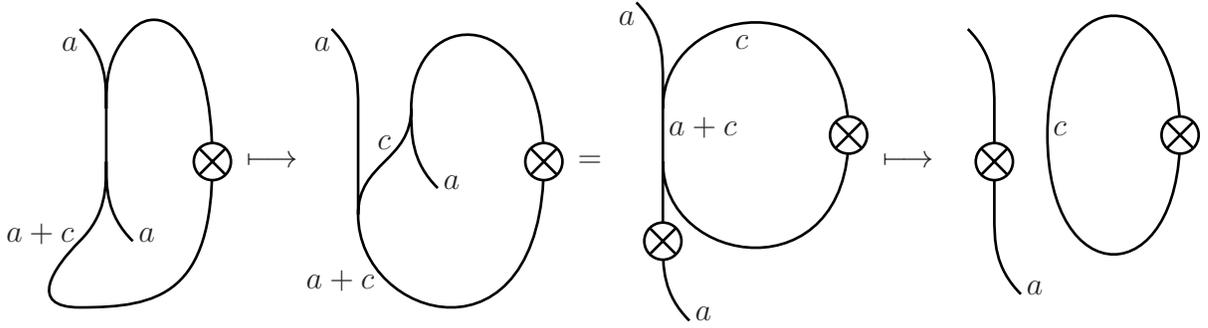

\centerline{\includegraphics{case24.eps}\put(-88,30){$a+c$}
\put(-38,30){$a$}\put(-67,102){$a$}
\hskip0.1cm\raisebox{58pt}{$\longmapsto$}\hskip0.1cm
\includegraphics{case24a.eps}\put(-100,12){$a+c$}
\put(-48,50){$a$}\put(-97,102){$a$}\put(-73,65){$c$}
\hskip0.1cm\raisebox{58pt}{$=$}\hskip0.1cm
\includegraphics{case24b.eps}\put(-78,70){$a+c$}
\put(-68,0){$a$}\put(-97,112){$a$}\put(-53,103){$c$}
\hskip0.1cm\raisebox{58pt}{$\longmapsto$}\hskip0.1cm
\includegraphics{case24c.eps}\put(-68,10){$a$}\put(-58,70){$c$}}
\caption{Simplification in Case J}\label{casej}
\end{figure}

\smallskip
\noindent{\bf Case K.}
Let $k=[c/a]+1$, $d=c-(k-1)a$. We have $0\leqslant d<a$, $k\geqslant3$.
We apply a $k$-times multiple splitting (see Figure~\ref{casek}),
which gives
\begin{figure}[ht]
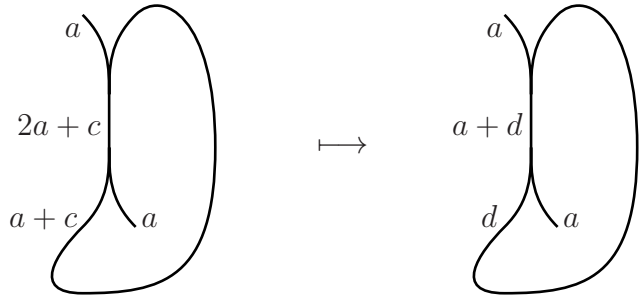

\centerline{\includegraphics{case25.eps}\put(-88,30){$a+c$}
\put(-38,30){$a$}\put(-67,102){$a$}\put(-85,65){$2a+c$}
\hskip1cm\raisebox{58pt}{$\longmapsto$}\hskip0.4cm
\includegraphics{case25.eps}\put(-69,30){$d$}
\put(-38,30){$a$}\put(-67,102){$a$}\put(-80,65){$a+d$}}
\caption{Simplification in Case K}\label{casek}
\end{figure}
$$2^g=\frac{(ka+d+1)((k+1)a+d+1)}{(d+1)(a+d+1)}
\geqslant\frac{(ka+d+1)((k+1)a+d+1)}{2(d+1)a}>\frac{(k+1)^2}2>k+1.$$
Thus, $g>\log_2(k+1)=p$.

\smallskip
\noindent{\bf Case L.}
We apply two splittings (see Figure~\ref{casel}), which gives
$$2^g=\frac{(a+b+c+1)(a+c+1)}{c+1}>4,$$
so, $g>2$, $p=2$.
\begin{figure}[ht]
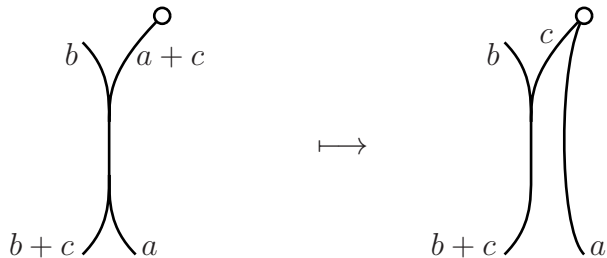

\centerline{\includegraphics{case3.eps}\put(-88,20){$b+c$}%
\put(-38,20){$a$}\put(-67,92){$b$}\put(-40,92){$a+c$}%
\hskip1cm\raisebox{58pt}{$\longmapsto$}\hskip0.4cm%
\includegraphics{case3a.eps}\put(-88,20){$b+c$}%
\put(-28,20){$a$}\put(-67,92){$b$}\put(-47,100){$c$}}
\caption{Simplification in Case L}\label{casel}
\end{figure}

\smallskip

This completes the proof of the first claim of the proposition.

The procedure above gives explicitly an algorithm to find the desired simplification
sequence. The only thing we need is to estimate the number of operations.
Throughout the procedure we operate with integers whose absolute value
is bounded by $2^{|(\theta,w)|}$.

In every case except Case K we need to perform a bounded number of additions (subtractions)
of such numbers. Thus, the amount of work at every step where we don't have Case~K
is $O(|(\theta,w)|)$. In Case K we perform additionally a single division, which consumes
time $O\bigl(|(\theta,w)|\cdot\log_2(k+1)\bigr)$, where $k$ is the multiplicity of the splitting.

Thus, in all cases the time consumed at each step of the algorithm
is bounded by the total cost of the step multiplied by $|(\theta,w)|$,
which implies the second claim of the Proposition.
\end{proof}

\begin{prop}\label{complexitiesequivalence}
Let $\mathscr A\subset G$ be a generating set as in Theorem~\ref{maintheo}, and let~$T$ be
a triangulation of~$M$. Then the zipped word length function
$\zwl_{\mathscr A}$ is comparable to the matrix complexity function $c_T$.
\end{prop}

\begin{proof}
We start from proving that $c_T\preceq\zwl_{\mathscr A}$.

Due to Proposition~\ref{cccc} it
suffices to show that $c_T(g^k)$ grows with $k$ not faster than $\log|k|$ when~$g$ is a Dehn twist.

Let $g$ be a Dehn twist along a simple closed curve $\gamma$. For any multiple
curve $\alpha$ such that~$\gamma$ and $\alpha$ are tight,
the image $g^k(\alpha)$ can be obtained from the union of $\alpha$ and
$|k|\cdot\langle\gamma,\alpha\rangle$ parallel copies of $\gamma$ by
resolving intersections of those copies with $\alpha$. Therefore,
for any edge~$e_i$ of the triangulation~$T$ we have
$$\langle g^k(\alpha),e_i\rangle\leqslant|k|\cdot\langle\gamma,\alpha\rangle\cdot\langle
\gamma,e_i\rangle+
|\langle e_i,\alpha\rangle|,$$
which implies
$$c_T(g^k)\leqslant C\cdot\log_2(|k|+1)$$
for some constant $C$.

Now we will show that $\zwl_{\mathscr A}\preceq c_T$.

Denote by~$K$ the number of branches in a universal train track (which is clearly
independent on the choice of the latter).
There are only finitely many, up to a self-homeomorphism of $(M,\mathscr P)$, train tracks in $M$ with
at most~$K$ branches. So, we
can fix a finite subset $X$ of train tracks such that:
\begin{enumerate}
\item
every train track in $X$ carries a triangulation;
\item
for any train track $\theta$ that carries a triangulation and has not more than~$K$ branches, there is an element $g\in G$ such that $g(\theta)\in X$.
\end{enumerate}

If a train track $\theta$ carries a triangulation, then the set of $g\in G$ such that $g(\theta)\sim\theta$
is also finite. Therefore, there is a finite subset $H$ of $G$ such that for any
simplification move $(\theta_1,w_1)\mapsto(\theta_2,w_2)$ with $\theta_1,\theta_2$ carrying a
triangulation and having at most~$K$ branches, and any $g_1,g_2\in G$ such that $g_1(\theta_1),g_2(\theta_2)\in X$,
the following holds:
\begin{enumerate}
\item
if $(\theta_1,w_1)\mapsto(\theta_2,w_2)$ is not a multiple splitting, then
$g_1g_2^{-1}\in H$;
\item
if $(\theta_1,w_1)\mapsto(\theta_2,w_2)$ is a $k$-times multiple splitting,
then there is a Dehn twist $d\in H$ such that $g_1g_2^{-1}=ad^k$ for
some $a\in H$.
\end{enumerate}
Thus, in both cases $\zwl_{\mathscr A}(g_1g_2^{-1})$ is bounded from above by
$C\cdot p$, where $C$ is a constant and~$p$ is the cost of the move $(\theta_1,w_1)\mapsto(\theta_2,w_2)$.
In the multiple splitting case this is due to the hypothesis that every Dehn twist
is conjugate to a fractional power of an element from $\mathscr A$.

We may assume without loss of generality that $\cup_{i=1}^Ne_i\in X$ and $\theta_T\in X$, where~$\{e_i\}_{i=1,\ldots,N}$
is the set of all edges of~$T$.

Now let $g\in G$ be any element different from $1$, and let~$w_0$ be a width assignment to
the branches of $\theta_T$ such that $(\theta_T,w_0)$ encodes $g(\cup_{i=1}^Ne_i)$ in the minimal way.
Pick a sequence of simplification moves
$$(\theta_T=\theta_0,w_0)\mapsto(\theta_1,w_1)\mapsto\ldots\mapsto
(\theta_r=g(\cup_{i=1}^Ne_i),w_r)$$
with total cost not larger than $3|(\theta_T,w_0)|$.
Such a sequence exists according to Proposition~\ref{costbound}.
Since the number of branches never grows under a simplification move,
all train tracks in this sequence have at most~$K$ branches.

Now for every $i=0,\ldots,r$ chose $g_i\in G$ so that $g_i(\theta_i)\in X$.
Specifically for $i=0$ and $r$ we put $g_0=1$ and $g_r=g^{-1}$.
We will have
$$\zwl_{\mathscr A}(g)=\zwl_{\mathscr A}\bigl((g_0g_1^{-1})(g_1g_2^{-1})\ldots(g_{r-1}g_r^{-1})\bigr)\leqslant
\sum_{i=1}^r\zwl_{\mathscr A}(g_{i-1}g_i^{-1})\leqslant3C|(\theta_T,w_0)|.$$
An application of Proposition~\ref{minimalencoding} completes the proof.
\end{proof}

\section{Counting intersections}\label{counting-sec}

Here we present the main technical result of the paper and prove Theorem~\ref{maintheo}, which asserts the existence of an efficient solution
of the word problem with respect to $\zwl_{\mathscr A}$.
Before starting the actual proof we mention briefly a strategy that we are \emph{not} going to follow,
but which yields another proof of the theorem.

For any fractional power $a$ of a fixed Dehn twist, one can construct an algorithm that produces
the normal coordinates of $a^k(\gamma)$ from the normal coordinates of a multiple curve~$\gamma$ and an
integer~$k$ in time $O(|\gamma|_T\cdot\log_2k)$, where~$|\gamma|_T$ is defined by~\eqref{gamma-norm-eq}.

Doing so for all generators from $\mathscr A$ yields \emph{a translation algorithm} from
the zipped word presentation to the matrix presentation. Given a zipped word representing an element
$g\in G$ it computes $\langle T,g(T)\rangle$ in time $O\bigl(\zwl_{\mathscr A}(g)^2\bigr)$
(if implemented properly).

The procedure from the proof of Proposition~\ref{complexitiesequivalence} used
to establish $\zwl_{\mathscr A}\preceq c_T$ can be turned into an actual algorithm that
performs the inverse translation, from the matrix presentation to the zipped word presentation,
and also consumes $O\bigl(\zwl_{\mathscr A}(g)^2\bigr)$ amount of time. The output
of the algorithm is a zipped word representing $g$ and depending only on $g$ but not
on the original presentation. Thus, this output can be taken for the normal form of $g$.

The strategy that we do follow is not to translate back and forth,
and use only the matrix presentation. The key ingredient missing so far
is the following statement.

\begin{prop}\label{compute-index-prop}
There exists an algorithm that, given the normal coordinates
of two multiple curves $\gamma_1$ and $\gamma_2$, computes
$\langle\gamma_1,\gamma_2\rangle$ in time $O\bigl(|\gamma_1|_T\cdot|\gamma_2|_T\bigr)$
on a RAM machine.
\end{prop}

\begin{proof}
We will use a modification of the procedure from the proof of Proposition~\ref{costbound}. This time we are going
to simplify two train tracks simultaneously to an extent that allows to detect all intersections
between~$\gamma_1$ and $\gamma_2$. We subdivide the proof into several subsections.

\subsection{General strategy and notation}
At every step of the algorithm, the multiple curves~$\gamma_1$ and $\gamma_2$ are encoded
by measured train tracks denoted $(\theta_1,w_1)$ and $(\theta_2,w_2)$, respectively, which are being modified
during the process. What data representing~$\theta_1$ and~$\theta_2$
is actually kept in computer's memory is described in Subsection~\ref{data-subsec}.

Here is the skeleton of the algorithm.

\begin{itemize}
\item The algorithm receives as input the vectors of normal coordinates of~$\gamma_1$ and~$\gamma_2$.

\item We start from $\theta_1=\theta_2=\theta_T$ and compute $w_1$ and $w_2$ so as to obtain the minimized representation
of the isotopy classes of $\gamma_1$ and $\gamma_2$ by $\theta_T$ as described in the proof of Proposition~\ref{minimalencoding}.

\item We run the simplification process for~$(\theta_1,w_1)$ and~$(\theta_2,w_2)$ as described in Subsections~\ref{moves-subsec} and~\ref{rules-subsec}.
For branches $\alpha$ and $\beta$ of $\theta_1$ and $\theta_2$, respectively, the number of
their transverse intersections is counted during the simplification process. By abusing
notation slightly we denote this number by~$\langle\alpha,\beta\rangle$. If $\alpha=\beta$ is
the closure of a proper arc we set $\langle\alpha,\beta\rangle=-1$.

Initially we set~$\langle\alpha,\beta\rangle=0$ for all branches~$\alpha$ and~$\beta$ of~$\theta_1$
and~$\theta_2$, respectively, as there are no transverse intersections
of~$\theta_1$ and~$\theta_2$ and no proper arcs in any of them.
These numbers are updated during the simplification process
whenever a new intersection (or coincidence of proper arcs)
is detected or any of~$\theta_1$ and~$\theta_2$ is modified.

\item
When the simplification finishes, we compute
\begin{equation}\label{bigsum}
\langle\gamma_1,\gamma_2\rangle=\sum_{\alpha,\beta}\langle\alpha,\beta\rangle\,w_1(\alpha)w_2(\beta),
\end{equation}
where the sum is taken over all branches $\alpha$ of $\theta_1$ and $\beta$ of $\theta_2$, and this is the output.
\end{itemize}

At every step of the simplification process, the train tracks $\theta_1$ and $\theta_2$
partially coincide and have a finite number of transverse intersection points. The latter
may not occur at switches of~$\theta_1$ and~$\theta_2$. We denote by $\theta_\cap$
the set of all non-isolated points of $\theta_1\cap\theta_2$, and by $\theta_\pitchfork$ the set
of all isolated ones. At every step, the intersection $\theta_1\cap\theta_2$ is homeomorphic
to a simplicial complex of dimension $\leqslant1$, with $\theta_\cap$
being the $1$-dimensional part of $\theta_1\cap\theta_2$
and $\theta_\pitchfork$ the $0$-dimensional one.

We think of $\theta_1\cup\theta_2$ as `a train track with self-intersections' and use
notation $\theta_\cup$ for the `abstract train track' of which $\theta_1\cup\theta_2$
is the image under an immersion $\theta_\cup\rightarrow M$. The formal meaning of $\theta_\cup$
will not be needed, but it will be handy to define branches and switches of $\theta_\cup$.

\begin{defi}
By \emph{a branch of $\theta_\cup$} we mean any of the following:
\begin{enumerate}
\item
a connected component of~$\alpha\cap\beta$ different from a single point, where $\alpha$
and $\beta$ are branches of $\theta_1$ and $\theta_2$, respectively;
\item
the closure of a connected component of $\alpha\setminus\theta_\cap$
where $\alpha$ is a branch of $\theta_1$ or $\theta_2$.
\end{enumerate}

By \emph{a switch of $\theta_\cup$} we mean a point $p\in(\theta_1\cup\theta_2)\setminus(\theta_\pitchfork
\cup\mathscr P)$ such
that the intersection of $\theta_1\cup\theta_2$ with any small
neighborhood of $p$ is not an arc.
\end{defi}

We allow only $3$-valent switches of $\theta_\cup$, which means that
exactly three branches of $\theta_\cup$ (counted with multiplicity) join at every
switch. There are, however, four ways how a switch of $\theta_\cup$ can arise.

A switch of $\theta_\cup$ can be a switch of both $\theta_1$ and $\theta_2$.
For a small enough neighborhood $U$ of such a switch we have $U\cap\theta_1=U\cap\theta_2$.

A switch of $\theta_\cup$ can be a switch of $\theta_1$ but not of $\theta_2$.
For a small enough neighborhood $U$ of such a switch we have $\theta_2\cap U\subset\theta_1\cap U$.
In particular, the intersection $\theta_2\cap U$ may be empty.

Similarly, a switch of $\theta_\cup$ can be a switch of $\theta_2$ bot not of $\theta_1$.

Finally, a switch of $\theta_\cup$ can be neither a switch of $\theta_1$ nor a switch of $\theta_2$.
The intersection of a small enough neighborhood $U$ of such a point with either $\theta_1$
or $\theta_2$ is an arc. We call such a switch \emph{a divergence point}.
Among three branches of $\theta_\cup$ joining at a divergence point, exactly one is contained in $\theta_\cap$,
one in $\theta_1$ but not in $\theta_2$, and one in $\theta_2$ but not in $\theta_1$ (see Figure~\ref{diverg}).
\begin{figure}[ht]
\centerline{\includegraphics{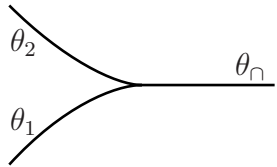}\put(-100,13){$\theta_1$}\put(-100,43){$\theta_2$}\put(-15,34){$\theta_\cap$}}
\caption{A divergence point}\label{diverg}
\end{figure}

By \emph{the width} $w(\alpha)$ of a branch $\alpha$ of $\theta_\cup$ we call the pair $(u_1,u_2)$ in which
$u_i$ is equal to~$w_i(\beta)$ if a branch $\beta$ of $\theta_i$ contains $\alpha$, and $0$
if $\alpha$ is not contained in $\theta_i$.

The terms `ingoing', `outgoing', `free', and `wide' have
the same meaning for branches of $\theta_\cup$ as for branches of an ordinary train track.

\subsection{Simplification moves}\label{moves-subsec}
Here we introduce certain transformations of the pair $\bigl((\theta_1,w_1),(\theta_2,w_2)\bigr)$
under which each of the two measured train tracks is modified either by an isotopy or by a simplification move introduced
in Section~\ref{simplifying}.

\smallskip\noindent\emph{Removing trivial branches.}
For each of the measured train tracks $(\theta_1,w_1)$ and $(\theta_2,w_2)$
this move consists, as before, in removing trivial branches. On the level of $\theta_\cup$
this means that all non-free branches of $\theta_\cup$ width $(0,0)$ are removed, and those
whose width has the form $(0,u)$ (respectively, $(u,0)$) with $u>0$
are thought of as being contained in $\theta_2$ but not in $\theta_1$ (respectively, in $\theta_1$ but not in $\theta_2$).

This move has the highest priority. So, in the sequel, whenever the width of a branch has the form $(0,u)$
of $(u,0)$ we assume that it is no longer contained in $\theta_\cap$ and any non-free branch of $\theta_\cup$
of width $(0,0)$ has to be erased.

\smallskip\noindent\emph{Ordinary splitting of $\theta_\cup$.}
If $\alpha$ is a wide branch of $\theta_\cup$ we can perform a splitting on $\alpha$.
This will be done only if~$w_1(\alpha)>0$ and~$w_2(\alpha)>0$.
If one of the endpoints of $\alpha$
is a puncture, then it is done exactly as in the case of a single train track, see the bottom of Figure~\ref{ordsplit}.
For each of the measured train tracks $(\theta_i,w_i)$, $i=1,2$, this will result in an ordinary splitting or just an isotopy
depending on whether or not $\alpha$ is a wide branch of~$\theta_i$.

If both endpoints of $\alpha$ are switches of $\theta_\cup$, a splitting on $\alpha$ will mean the modification of $\theta_\cup$, $w_1$, $w_2$
shown in Figure~\ref{doublesplit}, where the indicated widths of the branches are related as follows:
$$\begin{aligned}
c_1&=\max(0,a_1+b_1-a_1'-b_1'),&
c_2&=\max(0,a_2+b_2-a_2'-b_2'),\\
d_1&=\max(0,a_1'+b_1'-a_1-b_1),&
d_2&=\max(0,a_2'+b_2'-a_2-b_2).
\end{aligned}$$
and it is understood
that branches of width $(0,0)$, if any, must be erased.
\begin{figure}[ht]
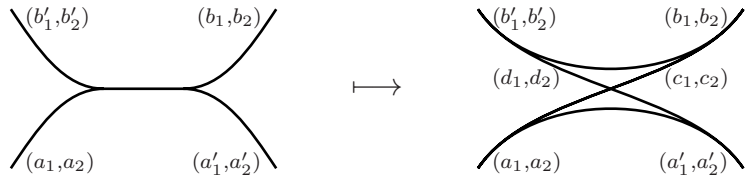

\centerline{\includegraphics{split1.eps}\put(-95,0){$\scriptstyle(a_1,a_2)$}\put(-95,55){$\scriptstyle(b_1',b_2')$}%
\put(-32,0){$\scriptstyle(a_1',a_2')$}\put(-30,55){$\scriptstyle(b_1,b_2)$}
\hskip1cm\raisebox{28pt}{$\longmapsto$}\hskip1cm
\includegraphics{split5.eps}\put(-95,0){$\scriptstyle(a_1,a_2)$}\put(-95,55){$\scriptstyle(b_1',b_2')$}%
\put(-32,0){$\scriptstyle(a_1',a_2')$}\put(-30,55){$\scriptstyle(b_1,b_2)$}%
\put(-30,32){$\scriptstyle(c_1,c_2)$}\put(-95,32){$\scriptstyle(d_1,d_2)$}}
\caption{Splitting two train tracks simultaneously}\label{doublesplit}
\end{figure}

For each of the measured train tracks $(\theta_i,w_i)$, $i=1,2$, this operation may be an isotopy, an ordinary splitting, or an ordinary splitting
followed by removing a trivial branch. All combinations can occur. A new point of $\theta_\pitchfork$ may
or may not be introduced depending of the widths of the branches adjacent to the endpoints of $\alpha$.
If it is introduced we say that $\theta_1$ and $\theta_2$ \emph{disagree} on the branch being splitted.

\smallskip\noindent\emph{Multiple splitting of $\theta_\cup$.}
We refer again to Figure~\ref{multsplit}. Now all widths are not just integers but integral vectors from $\mathbb Z^2_{\geqslant0}$:
$a=(a_1,a_2)$, etc. Multiple splitting of $\theta_\cup$
will be used only if $b_i\geqslant a_i>0$, $i=1,2$. Again, it is equivalent
to applying $k$ times an ordinary splitting on the wide branch in the fragment,
where for $k$ we take the largest integer satisfying $b_1\geqslant ka_1$, $b_2\geqslant ka_2$. For each of
the measured train tracks $(\theta_i,w_i)$, $i=1,2$, this will result in a $k$-multiple
splitting in the previously defined sense.

\smallskip\noindent\emph{Separation of circles.}
Suppose $\theta_\cap$ has a connected component $\sigma$ that is a two-sided
simple curve containing exactly four switches
of $\theta_\cup$ on $\sigma$ two of which are switches of $\theta_1$
and the other two of $\theta_2$. Suppose also that branches of $\theta_i$
approach $\sigma$ from both sides, $i=1,2$. Thus, the parts of $\gamma_1$ and $\gamma_2$
located in a small neighborhood of $\sigma$ have the form of `spirals'.

Finally, suppose that the `spirals' of $\gamma_1$ and $\gamma_2$
are twisted in opposite ways. Formally this means the following.
Let $\alpha$ be a smooth arc in a small neighborhood of $\sigma$ such that:
\begin{enumerate}
\item $\alpha$ is contained in $\theta_1\cup\theta_2$;
\item one of the endpoints of $\alpha$ is in $\theta_1\setminus\theta_2$ and
the other in $\theta_2\setminus\theta_1$;
\item
the intersection $\alpha\cap\theta_\cap$ is an arc contained in $\sigma$.
\end{enumerate}
Then the endpoints of $\alpha$ are on the same side of $\sigma$, see the left picture on Figure~\ref{separating}.

Then we can deform $\theta_1$ and $\theta_2$ so as to obtain disjoint simple
closed curves $\sigma_1\subset\theta_1$, $\sigma_2\subset\theta_2$ close to $\sigma$
and such that $\sigma_1$ has a single intersection point with $\theta_2$ and so does
$\sigma_2$ with $\theta_1$; see Figure~\ref{separating}.
\begin{figure}[ht]
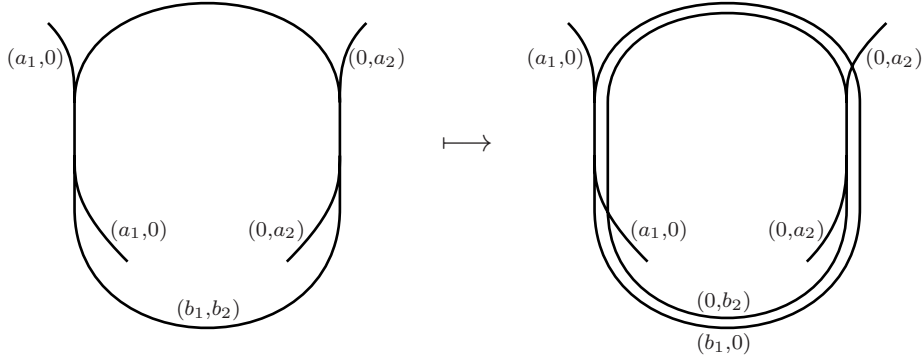

\centerline{\includegraphics{spiral4.eps}\put(-72,10){$\scriptstyle(b_1,b_2)$}\put(-97,40){$\scriptstyle(a_1,0)$}%
\put(-136,105){$\scriptstyle(a_1,0)$}%
\put(-45,40){$\scriptstyle(0,a_2)$}\put(-7,105){$\scriptstyle(0,a_2)$}%
\hskip1cm\raisebox{72pt}{$\longmapsto$}\hskip1cm\includegraphics{spiral5.eps}\put(-72,13){$\scriptstyle(0,b_2)$}%
\put(-72,-3){$\scriptstyle(b_1,0)$}\put(-97,40){$\scriptstyle(a_1,0)$}%
\put(-136,105){$\scriptstyle(a_1,0)$}%
\put(-45,40){$\scriptstyle(0,a_2)$}\put(-8,105){$\scriptstyle(0,a_2)$}%
}
\caption{Separation of circles}\label{separating}
\end{figure}

Such modification of $\theta_1$, $\theta_2$ will be referred to as \emph{a separation of circles}.

We also will use this name in the situation when a connected component $\sigma$
of $\theta_\cap$ has the form of a two-sided simple curve that is a free branch of $\theta_1$
and contains two switches of $\theta_2$ such that there are branches of $\theta_2$
approaching $\sigma$ from both sides. Separation of circles in this case
works as before, Figure~\ref{separating} will illustrate this if the branches of width $(a_1,0)$ are erased.

\subsection{Simplification rules}\label{rules-subsec}
The simplification procedure starts with two
measured train tracks $(\theta_1,w_1)$, $(\theta_2,w_2)$ such that $\theta_1=\theta_2=\theta_\cap=\theta_T$
and modifies them so as to end up with the situation in which $\theta_\cap$ consists of common free branches of
$\theta_1$ and $\theta_2$.

The simplification rules are not symmetric with respect to $\theta_1$ and $\theta_2$. The process
starts from checking which of the multiple curves $\gamma_1$ and $\gamma_2$ is simpler.
If $|\gamma_1|_T>|\gamma_2|_T$, then their roles are exchanged. So, we suppose in the sequel
that $|\gamma_1|_T\leqslant|\gamma_2|_T$. (More honestly, we compute an approximate
value of~$|\gamma_i|_T$ by using the~$[\log_2]$ function instead of~$\log_2$, which is much quicker. This produces a bounded
error, which can be ignored in this context.)

The simplification runs as follows.
\begin{description}
\item[Step 1] Remove all trivial branches of $\theta_1$, $\theta_2$.
\item[Step 2] If $\theta_\cap\setminus\mathscr P$ has a contractible connected component
with at least one switch of $\theta_\cup$ in it then:
\begin{enumerate}
\item
do a splitting on a wide branch of $\theta_\cup$ contained
in (the closure of) this component;
\item
repeat this step.
\end{enumerate}
\item[Step 3] If there is a switch of $\theta_1$ contained in $\theta_\cap$ such that
its outgoing branch $\alpha$ of $\theta_1$ contains a divergence point (equivalently,
is not covered by $\theta_\cap$) then:
\begin{enumerate}
\item
run the cleanup process described below, on $\alpha$;
\item
return to Step 2;
\end{enumerate}
\item[Step 4] If there is a branch of $\theta_1$ that is also a wide branch of $\theta_\cup$
on which $\theta_1$ and $\theta_2$ disagree then:
\begin{enumerate}
\item
split this branch;
\item
return to Step 2;
\end{enumerate}
\item[Step 5] Run the simplification procedure from the proof of Proposition~\ref{costbound} for $(\theta_1,w_1)$
with the following modifications:
\begin{enumerate}
\item
at every step, the widest branch $\alpha$ should be chosen only among wide branches of~$\theta_1$
contained in $\theta_\cap$;
\item
do not perform any splitting or sliding on a branch of~$\theta_1$ not contained in~$\theta_\cap$.
In any of the Cases~B--L shown in Figure~\ref{cases}, if the branch of width~$a+c$ (in~$(\theta_1,w_1)$)
is not covered by~$\theta_2$ just perform a splitting of the branch of
width~$a+b+c$;
\item should an ordinary splitting or sliding be performed on a branch $\alpha$ of $\theta_1$, first
remove all switches of $\theta_\cup$ from $\alpha$ by running the cleanup process on $\alpha$ and
then perform a splitting or sliding on $\alpha$ with $\theta_\cup$ so as to have the desired modification of $\theta_1$;
\item should a multiple splitting be performed on a circle $\sigma$ of $\theta_1$, first
run the cleanup process for $\sigma$ and then do either a multiple splitting on $\sigma$
with $\theta_\cup$ or a separation of circles, whichever is applicable. If none of these
can be applied perform an ordinary splitting of $\theta_\cup$ on the wide branch
contained in $\sigma$;
\item
whenever during this process the number of branches of $\theta_1$ contained in~$\theta_\cap$
decreases interrupt the process and return to Step~2;
\item
after each round of simplification return to Step~3.
\end{enumerate}
\item[Step 6] If $\theta_\cap$ has a connected component $\sigma$ that is a free branch of $\theta_1$
having the form of a simple curve and containing a switch of $\theta_2$, then:
\begin{enumerate}
\item
run the cleanup procedure for this component;
\item
if two switches of $\theta_2$ remain on $\sigma$ do a separation of circles on $\sigma$;
\item
repeat this step.
\end{enumerate}
\end{description}

Now we describe the cleanup procedure. The general principles are as follows:
\begin{enumerate}
\item
we apply simplification moves to $\theta_\cup$ so that $\theta_1$ does not change (or changes by isotopy);
\item
we remove switches of $\theta_\cup$ from a branch or a circle consisting of two branches of $\theta_1$ so that either the desired
simplification move of $\theta_1$ becomes extendable to a simplification move of $\theta_\cup$ or
it becomes possible to apply a separation of circles.
\end{enumerate}

The cleanup procedure appears in three different versions.

\smallskip\noindent\emph{Cleanup of a single branch of $\theta_1$ having the form of an arc.}
Let $\alpha$ be a branch of $\theta_1$ not forming a simple closed curve.
We suppose that there are some switches
of $\theta_\cup$ in the interior of $\alpha$ and we want to get rid of them.
A tail of a branch of $\theta_\cup$ not contained in $\theta_1$ is attached
to every such switch. We call these tails \emph{shavings}.

The branch $\alpha$ is locally two-sided, so we can choose one side to be top and the other to be bottom.
We can also orient $\alpha$ and think of this orientation as being from left to right. Having fixed this orientations
we can sort shavings and the corresponding
switches of $\theta_\cup$ contained in~$\alpha$ into four types: bottom-left, bottom-right, top-left, and top-right
according to the direction from which the corresponding shaving approaches the switch.

At the first stage of the cleanup we move left (top and bottom) shavings to the right and right shavings to the left of $\alpha$ as well
as reduce the number of shavings of each type to at most one. This is done by performing
splittings on wide branches of $\theta_\cup$ contained in the interior of~$\alpha$ (see Figure~\ref{cleanup1} a,b),
slidings on branches of $\theta_\cup$ connecting switches of the same type (Figure~\ref{cleanup1} c), and,
if neither of these is possible but still there are two shavings of the same type, slidings on branches of~$\theta_\cup$
contained in the interior of $\alpha$ followed by another sliding reducing the number of shavings (Figure~\ref{cleanup1} d).
\begin{figure}[ht]
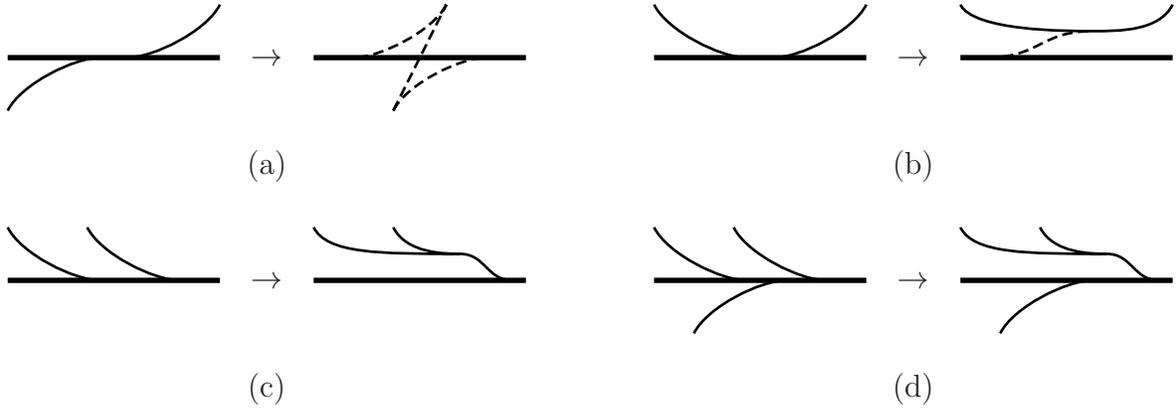

$$\begin{array}{ccc}
\includegraphics{cleanup1.eps}\quad\raisebox{27pt}{$\rightarrow$}\quad\includegraphics{cleanup2.eps}&\hskip1cm&\includegraphics{cleanup3.eps}\quad\raisebox{27pt}{$\rightarrow$}\quad\includegraphics{cleanup4.eps}\\
\text{(a)}&&\text{(b)}
\end{array}$$
$$\begin{array}{ccc}
\includegraphics{cleanup5.eps}\quad\raisebox{27pt}{$\rightarrow$}\quad\includegraphics{cleanup6.eps}&\hskip1cm&\includegraphics{cleanup7.eps}\quad\raisebox{27pt}{$\rightarrow$}\quad\includegraphics{cleanup8.eps}\\
\text{(c)}&&\text{(d)}
\end{array}$$
\caption{Cleanup procedure, first stage. Shown in bold is the branch $\alpha$, the others are branches of $\theta_\cup$
contained only in $\theta_2$. Dashed lines show branches  that are not necessarily all present after the move}\label{cleanup1}
\end{figure}

At the second stage of the cleanup we do splittings (Figure~\ref{cleanup2} a, b) and/or slidings (Figure~\ref{cleanup2} c, d),
whichever are applicable, on branches of $\theta_\cup$
contained in $\alpha$ and sharing an endpoint with~$\alpha$.
Each operation removes one switch of $\theta_\cup$ from $\alpha$.
\begin{figure}[ht]
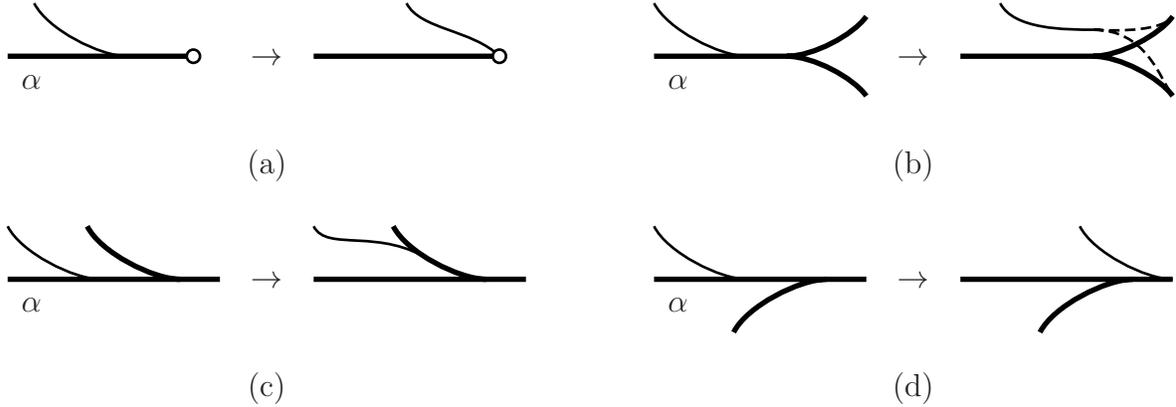

$$\begin{array}{ccc}
\includegraphics{cleanup9.eps}\put(-75,18){$\alpha$}\quad\raisebox{27pt}{$\rightarrow$}\quad\includegraphics{cleanup10.eps}&\hskip1cm&
\includegraphics{cleanup11.eps}\put(-75,18){$\alpha$}\quad\raisebox{27pt}{$\rightarrow$}\quad\includegraphics{cleanup12.eps}\\
\text{(a)}&&\text{(b)}
\end{array}$$
$$\begin{array}{ccc}
\includegraphics{cleanup13.eps}\put(-75,18){$\alpha$}\quad\raisebox{27pt}{$\rightarrow$}\quad\includegraphics{cleanup14.eps}&\hskip1cm&
\includegraphics{cleanup15.eps}\put(-75,18){$\alpha$}\quad\raisebox{27pt}{$\rightarrow$}\quad\includegraphics{cleanup16.eps}\\
\text{(c)}&&\text{(d)}
\end{array}$$
\caption{Cleanup procedure, the second stage. Shown in bold are branches of~$\theta_1$, dashed lines show branches of $\theta_\cup$
contained in $\theta_2$ one of which may not be present after the move}\label{cleanup2}
\end{figure}

\smallskip\noindent\emph{Cleanup of a circle consisting of two branches of $\theta_1$.} Let $\sigma$ be a circle consisting
of two branches of $\theta_1$ such that a multiple splitting of $\theta_1$ can be performed on $\sigma$.
We denote by $\alpha$ the wide branch of $\theta_1$ contained in $\sigma$, and the other branch by $\beta$.
We do the cleanup for $\alpha$ as described above, then for $\beta$, and again for $\alpha$. One can see
that after the second cleanup (for $\beta$) at most two shavings may remain on $\sigma$, and
after the third one they either escape from $\sigma$ or shift to $\beta$, in which case
we get the situation shown in Figure~\ref{separating} on the left.

\smallskip\noindent\emph{Cleanup of a free branch of $\theta_1$ having the form of a simple closed curve.}
Let $\sigma$ be a circular free branch of $\theta_1$ with some switches of $\theta_\cup$ on it.
Choose a point $p\in\sigma$ disjoint from those switches. Run the first stage of
cleanup for $\sigma\setminus\{p\}$ as if it is an ordinary branch of $\theta_1$. At most
four switches of $\theta_\cup$ will remain on $\sigma$. Then, if necessary,
do splittings on wide branches of $\theta_\cup$ contained in $\sigma$ until
at most one shaving remains on each side of $\sigma$. At most three splittings
are needed for that.

\smallskip
In Subsection~\ref{complexity-subsec} we will show that the simplification process defined above eventually stops,
(and estimate the asymptotic complexity of the algorithm). For the moment, we take it for granted.

\begin{lemm}\label{simplify-finished-lem}
When the simplification is completed, $\theta_1$ and~$\theta_2$ satisfy the following conditions\emph:
\begin{enumerate}
\item
each branch of $\theta_1$
is transverse to all branches of $\theta_2$ with an exception that some free branches of $\theta_1$
may coincide with free branches of~$\theta_2$\emph;
\item
the switches of both $\theta_1$ and $\theta_2$ are disjoint from $\theta_1\cap\theta_2$\emph.
\end{enumerate}
\end{lemm}

\begin{proof}
The simplification rules are designed so that after any elementary operation
that can result in creating a new contractible connected component
of~$\theta_\cap\setminus\mathscr P$ (these are operations at Step~3, Step~4, and those at Step~5
that satisfy Condition~(v)) we return to Step~2 and destroy such components or make
them free of switches of~$\theta_\cup$.

At each of Steps~3, 4, and~5 we either change something in~$\theta_\cup$ and return to one of the previous steps or do nothing.
So, proceeding with Step~6 means that, just before that, we passed Steps 3, 4, and~5 with no action.

Let~$\Omega$ be a connected component of~$\theta_\cap\setminus\mathscr P$ at the moment when we proceed to Step~6,
and let~$\alpha$ be the widest branch of~$\theta_1$ having a non-empty intersection with~$\Omega$.
Then~$\alpha$ cannot be a wide branch of~$\theta_1$ contained in~$\overline\Omega$, since otherwise
we would do something non-trivial at Step~4 or~5. It cannot be an ingoing branch of~$\theta_1$ for a switch
contained in~$\theta_\cap$, since otherwise a wider branch of~$\theta_1$ than~$\alpha$ would have
a non-empty intersection with~$\Omega$. It cannot also be an outgoing branch of~$\theta_1$ for a switch contained in~$\theta_\cap$,
since otherwise we would either do something non-trivial at Step~3 or one of the previously ruled out cases occurs. Finally,
it cannot happen that~$\Omega$ is a portion of~$\alpha$ not containing the endpoints,
since whenever such connected components of~$\theta_\cap\setminus\mathscr P$ appear they are destroyed by running Step~2.

We are left with the following two options.
\begin{enumerate}
\item
$\Omega$ is a proper arc and~$\alpha=\overline\Omega$ is a free branch of~$\theta_1$. In this case, $\alpha$
is also a free branch of~$\theta_2$. Indeed,
due to rule~(v) at Step~5, whenever a free branch of~$\theta_1$ having the form of the closure of a proper arc
emerges as a result of a splitting,
it is either disjoint from~$\theta_\cap\setminus\mathscr P$ or we return to~Step~2 and make
this branch free of switches of~$\theta_\cup$.
\item
$\Omega=\alpha$ is a free branch of~$\theta_1$ having the form of a simple closed curve.
If some switches of~$\theta_2$ remain on~$\Omega$ they are removed at Step~6.
\end{enumerate}

This completes the proof of the first claim of the proposition. The second claim follows
from the definition of the simplification moves: transverse intersections of~$\theta_1$ and~$\theta_2$
are never created at switches.
\end{proof}

\subsection{Data representation}\label{data-subsec}
The efficiency of implementation of the simplification procedure described above
depends heavily on the way in which the combinatorial data is represented.
The topological description of the algorithm might suggest that
at every step we should somehow keep track of how $\theta_1$ and $\theta_2$
are embedded in $M$, but this is \emph{not} the case. We only need `a local description'
of $\theta_\cup$ to implement the algorithm, which means the following.

We create a family of data objects, one per each branch and each switch of $\theta_\cup$, and
each puncture. For every switch and puncture we fix (arbitrarily) a surface orientation in its
small neighborhood. Each object keeps references to related objects (e.g., each puncture keeps references
to branches of $\theta_\cup$ adjacent to it) together
with orientation information: for each puncture and each switch it is a cyclic order of
the attached tails, and for each branch there is a boolean saying weather the orientations at
the endpoints agree along the branch.

Additionally, for every branch, we keep information about its width,
and, for every pair of branches~$\alpha$ and~$\beta$ of~$\theta_1$ and~$\theta_2$, respectively,
the number $\langle\alpha,\beta\rangle$ of their transverse intersections which is set to~$-1$
if~$\alpha=\beta$ is the closure of a proper arc).
All these data are updated at every simplification step in the way reflecting the change of~$(\theta_1,w_1)$, $(\theta_2,w_2)$, and~$\theta_\cup$.

One can see that the result of each simplification move can be computed
in terms of this data without any reference to the actual immersion $\theta_\cup\rightarrow M$,
and the number of operations needed for implementing a single move is bounded
by a constant. By an operation here we mean a creation or removal of an object, an assignment, or
an arithmetic operation.

\subsection{Why the output is $\langle\gamma_1,\gamma_2\rangle$?}
\indent
The aim of this subsection is to show that equality~\eqref{bigsum}
holds when the simplification described in the Subsections~\ref{moves-subsec} and~\ref{rules-subsec} halts. The proof is based on the following fact.

\begin{lemm}\label{tight-and-carried-lem}
The multiple curves~$\gamma_1$ and~$\gamma_2$ can be deformed by isotopies, and
the simplification of the train tracks~$\theta_1$, $\theta_2$ can be implemented so as to satisfy the following
two conditions:
\begin{enumerate}
\item
$\gamma_1$ and~$\gamma_2$ are in tight position;
\item
$\gamma_1$ and~$\gamma_2$ are carried by~$\theta_1$ and~$\theta_2$, respectively, at any stage
of the simplification process in a consistent way. The latter means that the foliations~$\mathscr F_{\theta_1}$
and~$\mathscr F_{\theta_2}$ \emph(see Section~\ref{trtrsec} for the notation\emph), to which~$\gamma_1$ and~$\gamma_2$ are
transverse, coincide and remain unchanged during the simplification.
Initially, $\gamma_1$ and~$\gamma_2$ are carried by~$\theta_T$ in the minimal way \emph(see Proposition~\ref{minimalencoding}\emph).
\end{enumerate}
\end{lemm}

\begin{proof}
We use the construction of the foliation~$\mathscr F_\theta$ from Section~\ref{trtrsec}. Let~$\mathscr F=\mathscr F_{\theta_T}$.
All modifications of the train tracks~$\theta_1$ and~$\theta_2$ made during the simplification process
can be performed so that all branches of~$\theta_1$ and~$\theta_2$
remain transverse to~$\mathscr F$ and close to~$\theta_T$, which means
that~$\mathscr F$ can be taken for~$\mathscr F_{\theta_1}$ and~$\mathscr F_{\theta_2}$ at any stage of the process.
Moreover, if the curves~$\gamma_1,\gamma_2$ have been chosen initially to be carried by~$\theta_T$, then they need not be changed
during the simplification to be carried by~$\theta_1$
and~$\theta_2$, respectively. We assume that they are carried by~$\theta_T$ in the minimal way.

Let~$U$ be an open neighborhood of~$\theta_T\setminus\mathscr P$
such that~$\mathscr F$ has no singularities in~$U$
and every leaf of~$\mathscr F|_U$ is an arc intersecting~$\theta_T\setminus\mathscr P$.
Clearly, we may assume that~$\gamma_1,\gamma_2$ are contained in~$U$, and so are~$\theta_1\setminus\mathscr P$
and~$\theta_2\setminus\mathscr P$ during the whole simplification process.

Now consider the pulling tight process for~$\gamma_1$, $\gamma_2$. 
Let~$D$ be a bigon of~$\gamma_1$ and~$\gamma_2$
which is going to be reduced, and let~$p$ be one of its corners.
\begin{figure}[ht]
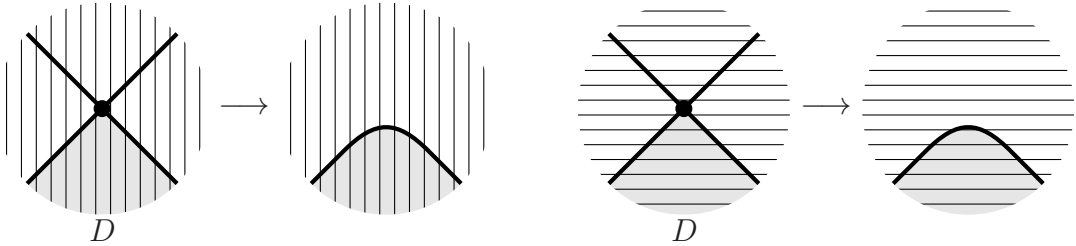

\includegraphics[scale=.8]{corner3.eps}\put(-45,-10){$D$}
\ \raisebox{38pt}{$\longrightarrow$}\
\includegraphics[scale=.8]{corner4.eps}
\hskip1cm
\includegraphics[scale=.8]{corner1.eps}\put(-45,-10){$D$}
\ \raisebox{38pt}{$\longrightarrow$}\
\includegraphics[scale=.8]{corner2.eps}
\caption{Two possible positions of a corner of~$D$ with respect to~$\mathscr F$}\label{smoothing-fig}
\end{figure}
Smoothing out the boundary~$\partial D$ near~$p$ by a small
perturbation (see Figure~\ref{smoothing-fig}) either makes it transverse to~$\mathscr F$ or creates an isolated tangency point of topological index~$1/2$
(like the one on the last picture of Figure~\ref{folsing}). We claim that the former case is impossible. Indeed,
the sum of the indexes of the singularities of~$\mathscr F$ located in the interior of~$D$ is non-positive,
since every connected component of~$M\setminus\overline U$ having a non-empty intersection with~$D$ is contained entirely in~$D$. So, the only way how we can have~$\chi(D)=1$ is that both corners contribute~$1/2$.

This means that, after the reduction of~$D$, the new curves~$\gamma_1$, $\gamma_2$ remain transverse to~$\mathscr F$,
and hence, are still carried by~$\theta_T$. Clearly,
the sum~$w_1+w_2$ of the width functions corresponding to~$\gamma_1,\gamma_2$ does not change.
Since both have been chosen initially in the minimal way, this implies that neither of them changes
under a bigon reduction.
\end{proof}

In what follows we assume that~$\gamma_1$ and~$\gamma_2$ satisfy the conditions from Lemma~\ref{tight-and-carried-lem} and keep the notation~$U$ from its proof.

A kind of tightness holds also for~$\theta_1$ and~$\theta_2$.
Namely, we have the following.

\begin{lemm}\label{theta-tight-lem}
At no stage of the simplification
procedure, there is a $2$-disc~$D\subset\overline U$ such that~$D\cap(\theta_1\cup\theta_2)=\partial D$
and~$\partial D$ consists of two smooth arcs.
\end{lemm}

\begin{proof}
 Indeed, this is true initially, since~$\theta_1$ and~$\theta_2$ coincide with~$\theta_T$,
and for no connected component~$V$ of~$U\setminus\theta_T$ we have~$\partial V\subset\theta_T$.
During the subsequent simplification, new connected components of~$U\setminus(\theta_1\cup\theta_2)$
are created (consult Figures~\ref{doublesplit} and~\ref{separating}), but the boundary of each of them
has three breaking points. Also, some of these connected components may join during the simplification forming
either an open $2$-disc with at least four breaking points at the boundary or a non-simply-connected domain.
\end{proof}

Denote by~$\pi_1$ and~$\pi_2$ the projections~$\pi_{\gamma_1}:\gamma_1\rightarrow\theta_1$
and~$\pi_{\gamma_2}:\gamma_2\rightarrow\theta_2$, respectively, defined for the state of~$\theta_1,\theta_2$
after the simplification has finished. We extend~$\pi_1$ and~$\pi_2$ to the closures~$\overline\gamma_1$,
$\overline\gamma_2$ by continuity, that is, put~$\pi_i(P_j)=P_j$ whenever~$P_j\in\overline\gamma_i$, $i=1,2$.

For any~$p\in\gamma_i$, $i=1,2$, define~$\delta_i(p)$ to be a (unique)
closed subarc, possibly degenerate to a point, of the leaf of~$\mathscr F|_U$
passing through~$p$ such that~$\partial\delta_i(p)=\{p,\pi_i(p)\}$.
Denote by~$\Gamma$ be the following subset of~$\gamma_1\times\gamma_2$:
$$\Gamma=\{(p_1,p_2)\in\gamma_1\times\gamma_2:\delta_1(p_1)\cap\delta_2(p_2)\ne\varnothing\},$$
and by~$\mathscr X$ the set of connected components of~$\Gamma$. Define a map~$s:\mathscr X\rightarrow\{-1,0,1\}$
as follows:
$$s(\beta)=\left\{\begin{aligned}
-1,&\text{ if $\beta$ is homeomorphic to an open interval};\\
1,&\text{ if $\beta$ is homeomorphic to a closed interval and $(p,p)\in\beta$ for some~$p\in\gamma_1\cap\gamma_2$};\\
0,&\text{ otherwise}.
\end{aligned}\right.$$

First, we claim that the following equality holds:
\begin{equation}\label{<>=s(b)}
\langle\gamma_1,\gamma_2\rangle=\sum_{\beta\in\mathscr X}s(\beta).
\end{equation}
Indeed, let~$\beta$ be a connected component of~$\Gamma$. Denote
the projections of~$\beta$ to~$\gamma_1$ and~$\gamma_2$ by~$\beta_1$ and~$\beta_2$,
respectively.

Suppose that~$\beta$ is an open arc. Since~$\Gamma$ is a closed subset of~$\gamma_1\times\gamma_2$,
this means that~$\beta_1$ and~$\beta_2$ are proper arcs. By construction, their projections~$\pi_1(\beta_1)$
and~$\pi_2(\beta_2)$ are isotopic, which implies, by Lemma~\ref{theta-tight-lem},
that~$\pi_1(\beta_1)=\pi_2(\beta_2)\subset\theta_\cap$.

We also have~$\pi_1(\beta_1')=\pi_2(\beta_2')\subset\theta_\cap$
for any parallel proper arcs~$\beta_1'\subset\gamma_1$ and~$\beta_2'\subset\gamma_2$.
Therefore, the number of connected components~$\beta$ of~$\Gamma$ with~$s(\beta)=-1$
is exactly the number of such pairs~$(\beta_1',\beta_2')$.

Now suppose that~$\beta$ is a connected component of~$\Gamma$ such
that~$(p,p)\in\beta$ for some~$p\in\gamma_1\cap\gamma_2$.
Due to the tightness of~$\gamma_1$ and~$\gamma_2$ there is no~$p'\in\gamma_1\cap\gamma_2$
distinct from~$p$ with~$(p',p')\in\beta$, and~$\beta$ is homeomorphic either to a closed
interval or to a circle. In the latter case, the closed curves~$\beta_1$, $\beta_2$ are isotopic to one another
and intersect once, which implies that they are one-sided. In the former case,
the connected components of~$\gamma_1$, $\gamma_2$ containing the arcs~$\beta_1$, $\beta_2$, respectively,
are not isotopic. Thus, the number of connected components~$\beta$ of~$\Gamma$ with~$s(\beta)=1$
is equal to the number of points in~$\gamma_1\cap\gamma_2$ less the number of pairs of isotopic
one-sided curves~$\beta_1'\subset\gamma_1$, $\beta_2'\subset\gamma_2$.

Thus, equality~\eqref{<>=s(b)} is settled.

Now we show that
\begin{equation}\label{s(b)=w_1w_2}
\sum_{\beta\in\mathscr X}s(\beta)=\sum_{\alpha_1,\alpha_2}\langle\alpha_1,\alpha_2\rangle\,w_1(\alpha_1)w_2(\alpha_2),
\end{equation}
where the sum is taken over all branches $\alpha_1$ of $\theta_1$ and $\alpha_2$ of $\theta_2$.
For a connected component~$\beta$ of~$\Gamma$, we again denote by~$\beta_1$ and~$\beta_2$
the projections of~$\beta$ to~$\gamma_1$, $\gamma_2$, respectively.

As we have seen above, the equality~$s(\beta)=-1$ means that~$\pi_1(\beta_1)=\pi_2(\beta_2)$ is a common
proper arc contained in~$\theta_\cap$. The number of such~$\beta$ is thus equal to
$$\sum_{\alpha}w_1(\alpha)w_2(\alpha),$$
where the sum is taken over all common free branches of~$\theta_1$ and~$\theta_2$.

If~$s(\beta)=1$, then~$\beta_1$ and~$\beta_2$ intersect once and do not belong to isotopic
one-sided closed components of~$\gamma_1$, $\gamma_2$, respectively. This implies that~$\pi_1(\beta_1)$ and~$\pi_2(\beta_2)$
also intersect once (more intersections would contradict Lemma~\ref{simplify-finished-lem} or
Lemma~\ref{theta-tight-lem}).
Therefore, the number of connected components~$\beta\subset\Gamma$ with~$s(\beta)=1$
is equal to the number of triples~$(p,p',q)\in\gamma_1\times\gamma_2\times\theta_\pitchfork$
such that~$\pi_1(p)=\pi_2(p')=q$, which, in turn, is equal to
$$\sum_{\alpha_1,\alpha_2}\langle\alpha_1,\alpha_2\rangle\,w_1(\alpha_1)w_2(\alpha_2),$$
where the sum is taken over all pairs of branches~$\alpha_1\subset\theta_1$, $\alpha_2\subset\theta_2$
intersecting transversely.

Thus, equality~\eqref{s(b)=w_1w_2} is also settled. Together with~\eqref{<>=s(b)}, it implies~\eqref{bigsum}.

\subsection{The asymptotic complexity of the algorithm}\label{complexity-subsec}
Recall that we assume the surface~$M$ and the set of punctures~$\mathscr P$ to be fixed once and for all.
So, in what follows, `bounded' means `bounded from above by a constant depending on~$M$ and~$\mathscr P$
but not on anything else'.

We need to show that the number of elementary arithmetic operations needed to compute~$\langle\gamma_1,\gamma_2\rangle$
has growth~$O(|\gamma_1|_T\cdot|\gamma_2|_T)$. The way to prove this is essentially the same
as that of Proposition~\ref{costbound}, so, we stop only on the differences.

We may assume that the number of switches of~$\theta_\cup$ remains bounded during the whole
simplification process. Indeed, the number of switches of~$\theta_1$ and~$\theta_2$
never increases, so, we should worry only about the number of divergence points.

When no branch of~$\theta_\cap$ connecting two divergence points is present, the number
of divergence points is clearly not larger than~$n+3q$, where~$q$
is the total number of switches of~$\theta_1$ and~$\theta_2$.

Suppose that, at some stage of the simplification process, there is a branch~$\alpha$ of~$\theta_\cap$
having the form of an arc whose both ends are divergence points. Such a branch
is untouched during the simplification until a moment when a splitting is performed
on it. As a result of this splitting, the branch~$\alpha$ and the two divergence points disappear.
We may reorder the simplification moves so that such splittings are performed
immediately as they become possible. Any other simplification move
may create at most two such branches of~$\theta_\cap$, which contribute at most four
divergence points in excess of the previous estimate. Thus,
we may assume that the number of divergence points never exceeds~$n+4+3q$.

This implies, in particular, that every cleanup procedure takes a bounded number of elementary operations.
Also bounded is the number of simplification moves performed in a row under Step~2.

We think of the procedure under Steps~1 and~5 as \emph{the regular simplification}, which
is \emph{interrupted} several times to perform simplification moves under Steps~2--4. Let~$C_1,C_2$ be
the number of branches and the number of switches of~$\theta_1$, respectively, contained in~$\theta_\cap$.
These numbers may only decrease during the simplification.

Interruptions of the regular simplification occur either immediately after one of~$C_1,C_2$ decreases
or just before such an event. Therefore, the total number of interruptions is bounded,
and so is the total cost of all simplification moves performed at the interruptions of the regular simplification.

Now we reconsider the regular simplification. To make the gain/cost ratio of every simplification step
bounded away from zero we put~$(\theta,w)=(\theta_1,w_1)$ and
redefine the gain by substituting~$A(\theta)$ in~\eqref{complexity0}
with the set of non-free branches of~$\theta_1$ \emph{contained in~$\theta_2$}.

In any of the Cases~B--L shown in Figure~\ref{cases}, if the branch of width~$a+c$ (in~$(\theta_1,w_1)$)
is not covered by~$\theta_2$, a splitting of the branch of
width~$a+b+c$ will reduce~$C_1$, and thus will have gain at least one.

Another feature of the simultaneous simplification of two train tracks is that
at every step we recompute~$w_2$. So, every simplification step
that does not involve a multiple splitting consumes time~$O(|\gamma_2|_T)$.

A special care is needed only in Case~K if a multiple splitting is involved.
This occurs when the two branches forming a circle
are contained in~$\theta_2$. Let their widths in~$(\theta_2,w_2)$ be~$a'+c'$ and~$2a'+c'$.

First, we compute~$k=[c/a]+1$, as before. Then we check whether or not~$ka'\leqslant c'$.
If the inequality holds true, we make a $k$-times multiple splitting,
which has the same gain as before, that is at least~$\log_2(k+1)$,
and consumes time~$O(\log_2(k+1)\cdot|\gamma_2|_T)$.

If we find out that~$ka'>c'$ we compute~$k'=[c'/a']+1$ and perform a $k'$-times multiple splitting,
which produces a wide branch~$\alpha$ at which~$w_1,w_2$ are equal to~$(c-(k'-2)a)$ and~$(c'-(k'-2)a')$,
respectively, such that~$\theta_1,\theta_2$ disagree on~$\alpha$. Subsequent
ordinary splitting on~$\alpha$ (which is done under Step~4) will produce a branch of~$\theta_1$ of width~$(c-k'a)$
which is no longer covered by~$\theta_2$.
The total gain will be at least
$$\log_2\frac{(2a+c+1)(a+c+1)}{(c-(k'-1)a+1)}\geqslant\log_2(2a+c+1)>\log_2(k+1),$$
whereas the computational time is again~$O(\log_2(k+1)\cdot|\gamma_2|_T)$.

When the regular simplification is finished, it is only bounded number of circle separations performed at Step~6,
which consumes bounded computational time.

Thus, the whole simplification procedure consumes time~$O(|\gamma_1|_T\cdot|\gamma_2|_T)$.
The same estimate works for the final computation of~$\langle\gamma_1,\gamma_2\rangle$
using formula~\eqref{bigsum}. This completes the proof of Proposition~\ref{compute-index-prop}.
\end{proof}

Proposition~\ref{compute-index-prop} gives an estimate for worst cases, but
for `typical' cases the running time of the algorithm might be better than~$O(|\gamma_1|_T\cdot|\gamma_2|_T)$.
This is due to the fact that, for random~$\gamma_1$, $\gamma_2$
having large complexity, it is likely that the simplification procedures
for~$(\theta_1,w_1)$ and~$(\theta_2,w_2)$ diverge well before~$(\theta_1,w_1)$ and~$(\theta_2,w_2)$
get much simpler, which means
that the simplification process described in the proof of Proposition~\ref{compute-index-prop}
will actually have much smaller cost than~$|\gamma_1|_T$.
This means, in turn, that the main contribution to the running
time will come from the computation of the expression in the right hand side of~\eqref{bigsum},
that is, from multiplying large numbers~$w_1(\alpha)$ and~$w_2(\beta)$ for all appropriate pairs
of branches~$(\alpha,\beta)$.

There exist faster methods for multiplication of natural numbers than the grade-school
algorithm. The first such method was proposed by A.\,Karatsuba~\cite{kar1,kar2}. Thus, it is plausible
that the average running time for computing~$\langle\gamma_1,\gamma_2\rangle$ can also be improved.

Propositions~\ref{product} and~\ref{compute-index-prop} immediately imply the following.

\begin{coro}\label{main-coro}
There is an algorithm for computing~$\langle T,(g_1g_2)(T)\rangle$ from~$\langle T,g_1(T)\rangle$ and~$\langle T,g_2(T)\rangle$
in time~$O\bigl(c_T(g_1)\cdot c_T(g_2)\bigr)$ on a RAM machine.
\end{coro}

\subsection{Proof of Theorem~\ref{maintheo}}
Due to Proposition~\ref{complexitiesequivalence} we can choose a triangulation~$T$
of~$M$ with vertices at~$\mathscr P$ and substitute~$c_T$ for~$\zwl_{\mathscr A}$
in the formulation of the theorem. So, we use the matrix presentation for
elements of~$G$.

Computing the normal form of an element amounts to
computing~$\langle T,(g_1\circ g_2\circ\ldots\circ g_k)(T))\rangle$
from~$\langle T,g_i(T)\rangle$, $i=1,\ldots,k$,
which can be done in quadratic time due to Propositions~\ref{product}
and~\ref{compute-index-prop}.

It remains to describe a procedure for checking whether a given word
represents a sequence of matrices of the form~$\langle T,g(T)\rangle$, $g\in G$.
First, it is a simple syntax check whether the given word
represents a sequence of integral matrices of the appropriate size.
Then each of the matrices should be tested for having the form~$\langle T,g(T)\rangle$.

Proposition~\ref{compute-index-prop}
allows to check in quadratic time whether any two distinct
columns of an $N\times N$-matrix represent normal curves having zero geometric
intersection index. Further, if this is true, take the linear combination of the columns
with coefficients~$1,2,\ldots,N$ and simplify
the obtained multiple curve by using
Propositions~\ref{minimalencoding} and~\ref{costbound}.
After the simplification, check whether the result
is, up to a homeomorphism, the triangulation~$T$
with the~$i$th edge repeated~$i$ times, $i=1,\ldots,N$.

\subsection{Asymptotic dependence on~$N$}\label{n-dep-subsec}
Throughout the paper, we have been assuming the punctured surface~$(M,\mathscr P)$ to
be fixed, and have been ignoring the question on the dependence of the complexity of the proposed algorithms on
the complexity of the surface. Now we will have a quick look on this dependence.
For the measure of complexity of~$(M,\mathscr P)$ we use the number~$N$
of edges of any triangulation of~$M$ with vertices at~$\mathscr P$ (which is clearly comparable
to the number of faces).

We make the following, quite realistic, assumption on the computational model:
we use a RAM machine with standard memory unit size large enough
to store integers comparable to the size of the input. So, numbers
like~$N$, $|\gamma|_T$, or~$c_T(g)$ are thought of as fitting
the standard memory unit size, whereas the numbers like~$\langle\gamma,e_i\rangle$
may occupy many standard memory units.

With these settings, Proposition~\ref{compute-index-prop} and Corollary~\ref{main-coro} can be strengthen as follows.

\begin{theo}\label{n-dep-theo}
There is an algorithm that accepts as input a triangulation~$T$ of a punctured
surface~$(M,\mathscr P)$ and two vectors of the normal coordinates of
multiple curves~$\gamma_1$, $\gamma_2$ on~$M$ \emph(or two matrices~$\langle T,g_1(T)\rangle$, $\langle T,g_2(T)\rangle$,
where~$g_1,g_2\in\mcg(M,\mathscr P)$\emph) and computes~$\langle\gamma_1,\gamma_2\rangle$
\emph(respectively, $\langle T,g_1g_2(T)\rangle$\emph) in time~$O(N^3\cdot|\gamma_1|_T\cdot|\gamma_2|_T)$
\emph(respectively, $O(N^5\cdot c_T(g_1)\cdot c_T(g_2))$\emph),
where~$N$ is the number of edges of~$T$.
\end{theo}

\begin{proof}[Sketch of proof]
We  revisit the proofs of all previously made statements about proposed algorithms. First, consider
the problem of computing~$\langle\gamma_1,\gamma_2\rangle$.

With a reasonable way to encode~$T$ at the input, the amount of work to produce
a universal train track~$\theta_T$ and to encode~$\gamma_1$ and~$\gamma_2$ by assigning weights to
the branches of~$\theta_T$ as described in Section~8 takes time~$O(N\cdot(|\gamma_1|_T+|\gamma_2|_T))$,
which is considerably faster than the required~$O(N^3\cdot|\gamma_1|_T\cdot |\gamma_2|_T)$.
Deciding whether we need to exchange~$\gamma_1$ and~$\gamma_2$ to have~$\gamma_1$
not (much) more complicated than~$\gamma_2$ is also pretty quick.

In the transition from a normal multiple curve~$\gamma$ to an isotopic multiple curve carried by~$\theta_T$,
as described in the proof of Proposition~\ref{minimalencoding} a normal arc in~$\gamma$
converts into~$O(N)$ almost normal arcs supported by~$\theta_T$, since there are only~$O(N)$
branches of~$\theta_T$. This means, that, for the obtained with assignment~$w$,
we have in worst case~$|(\theta_T,w)|=O(N\cdot |\gamma|_T)$.

So, the total gain that we should collect during the simplification process described
in the proof of~Proposition~\ref{costbound} with modifications described in Subsections~\ref{moves-subsec} and~\ref{rules-subsec}
does not exceed~$O(N\cdot|\gamma_1|_T)$.

Every round of regular simplification of~$\theta_\cup$ may include clean up procedures, each of which
requires~$O(N)$ simplification moves. For every round, we also need to find the appropriate
place in~$\theta_\cap$ where the simplification should occur, and this also consumes~$O(N)$ arithmetic operations.
Therefore, the estimates for the cost/gain ratio from the proof of~Proposition~\ref{costbound}
should be multiplied by~$O(N)$. In total, we get the estimate~$O(N^2\cdot|\gamma_1|_T\cdot|\gamma_2|_T)$
for the time consumed by the regular part of the simplification.

The regular simplification is interrupted~$O(N)$ times to destroy contractible connected components
of~$\theta_\cap\setminus\mathscr P$ containing switches of~$\theta_\cup$. To discover
that the given connected component is contractible it suffices to perform~$O(N)$ arithmetic operations
on `small' numbers. One can show that a contractible connected components of~$\theta_\cap\setminus\mathscr P$
having~$k$ switches in it is destroyed in~$O(k^2)$ splittings. So, the total amount of splittings
involved at Step~2 is bounded by~$O(N^3)$, and the time lost on the interruptions is estimated
as~$O(N^3\cdot|\gamma_2|_T)$.

At Step~6 we need to perform only~$O(N)$ splittings and slidings, which consume time~$O(N\cdot|\gamma_2|_T)$.

Finally, we use formula~\eqref{bigsum} to compute the output, which requires time~$O(N^2\cdot|\gamma_1|_T\cdot|\gamma_2|_T)$.
We see that none of the asymptotic terms we have encountered exceeds~$O(N^3\cdot|\gamma_1|_T\cdot|\gamma_2|_T)$.

The estimate for the running time to compute~$\langle T,g_1g_2(T)\rangle$ from~$\langle T,g_1(T)\rangle$
and~$\langle T,g_2(T)\rangle$ now follows from that for computing~$\langle\gamma_1,\gamma_2\rangle$,
since we simply need to compute~$N^2$ geometric intersection indexes.
\end{proof}

Our estimations in the proof of Theorem~\ref{n-dep-theo} are very rough,
so it is likely that the power of~$N$ can be lowered.

To conclude this section, we note that the question on the dependence of the running time on the complexity
of the surface makes sense for algorithms whose existence is established in Theorem~\ref{maintheo}
only when a concrete family of generating sets for \emph{all }mapping class groups has been fixed. The author is unaware if such a family
can be chosen so as to make this dependence polynomial.

\section{A concluding remark}\label{concluding-sec}
Suppose for time being that a triangulation~$T$ with vertices at~$\mathscr P$ is fixed on~$M$ such that
the first~$N_0$ edges of~$T$ are non-boundary edges, and the last~$N-N_0$ are boundary ones. Using Proposition~\ref{normalcoord}
let us identify the set of all multiple curves on~$(M,\mathscr P)$ not containing a proper arc
parallel to a boundary edge with a subset~$L\subset\mathbb Z^{N_0}$,
and treat the geometric intersection index (see Definition~\ref{intersection-def}) as a function~$\langle\,,\rangle:L\times L\rightarrow\mathbb Z$.
It is not hard to show that this function can be extended, in a unique way, to a continuous function~$\mathbb R^{N_0}\times\mathbb R^{N_0}\rightarrow\mathbb R$
satisfying the following homogeneity condition: $\langle\lambda x,y\rangle=\langle x,\lambda y\rangle=\lambda\langle x,y\rangle$ for all~$x,y\in\mathbb R^{N_0}$,
$\lambda\in\mathbb R_{\geqslant0}$.

The whole space~$\mathbb R^{N_0}$ in this construction admits a natural topological interpretation as the space of measured
laminations in the sense of W.\,Thurston~\cite{thurston}, which can be defined without a reference to a concrete triangulation~$T$.
We denote this space by~$\mathscr M(M,\mathscr P)$. Topologically it is just an open $N_0$-dimensional ball.
Different triangulations just give rise to different global coordinate systems on it.

Now every triangulation with vertices at~$\mathscr P$ (and with non-boundary edges numbered first) gives rise also to a point of the following space:
$$\mathscr Z(M,\mathscr P)=\{(e_1,e_2,\ldots,e_{N_0})\in\mathscr M(M,\mathscr P)^{N_0}:
\langle e_i,e_j\rangle=-\delta_{ij}\ \forall i,j=1,\ldots,N_0\}.$$
The following question sounds intriguing to the author: is~$\mathscr Z(M,\mathscr P)$ a topological manifold of dimension~$N_0(N_0-1)/2$?
If the answer is positive, what kind of manifold is it?

This is interesting because the natural action of~$\mcg(M,\mathscr P)$ on~$\mathscr Z(M,\mathscr P)$ is likely to be free and properly discontinuous, so knowing
the structure of this space may be useful for studying the group~$\mcg(M,\mathscr P)$.

The author was able to answer the questions above only in the following three simplest cases.
\begin{enumerate}
\item
$M=\mathbb D^2$, $|\mathscr P|=5$, $\mathscr P\subset\partial\mathbb D^2$. In this case, $\mathscr Z(M,\mathscr P)\cong\mathbb S^1\times\mathbb Z_2$.
\item
$M=\mathbb D^2$, $|\mathscr P|=6$, $\mathscr P\subset\partial\mathbb D^2$. In this case, $\mathscr Z(M,\mathscr P)\cong\mathbb RP^3\times\mathbb Z_2$.
\item
$M=\mathbb T^2$, $|\mathscr P|=1$. In this case, $\mathscr Z(M,\mathscr P)\cong\mathbb R^2\times\mathbb S^1\times\mathbb Z_2$.
(This case was classified as sporadic and excluded from the main consideration in this paper, but this was only because the matrix presentation discussed
throughout the paper gives rise to the quotient group~$\mcg(\mathbb T^2,\{P_1\})/(\pm1)\cong PSL(2,\mathbb Z)$
instead of the mapping class group itself.)
\end{enumerate}
In cases~(ii) and~(iii), establishing this was far from being straightforward. However, these cases may still be too simple to reflect
the general picture.

\end{document}